\documentclass{article}
\usepackage{amssymb}
\usepackage{amsmath}
\usepackage{amsthm}
\usepackage{bbold}

\newcounter{enucount}
\newcounter{enuref}
\newcounter{enualph}
\newcounter{enunumb}
\renewcommand{\theenucount}{(\roman{enucount})}

\renewcommand{\theenunumb}{(\arabic{enunumb})}
\renewcommand{\theenualph}{(\alph{enualph})}

\theoremstyle{plain}
\newtheorem{theorem}{Theorem}
\newtheorem{question}{Question}
\newtheorem{conj}[question]{Conjecture}
\newtheorem{problem}[question]{Problem}
\newtheorem{lem}[theorem]{Lemma}
\newtheorem{mainlem}[theorem]{Main Lemma}
\newtheorem{basiclem}[theorem]{Basic Lemma}
\newtheorem{prelem}[theorem]{Preliminary Lemma}
\newtheorem{obs}[theorem]{Observation}
\newtheorem{fact}[theorem]{Fact}
\newtheorem{cor}[theorem]{Corollary}
\newtheorem{prop}[theorem]{Proposition}
\newtheorem{sclaim}{Claim}[theorem]

\newtheorem*{maintheoremA}{Main Theorem A}
\newtheorem*{maintheoremB}{Main Theorem B}

\theoremstyle{definition}
\newtheorem{defin}{Definition}

\theoremstyle{remark}
\newtheorem{exam}{Example}
\newtheorem{counterexam}[exam]{Counterexample}

\renewenvironment{enumerate}{\begin{list}{\rm \theenucount}{\usecounter{enucount}
    \setlength{\labelwidth}{1cm}}}
   {\end{list}}
\newenvironment{alphenumerate}{\begin{list}{\rm \theenualph}{\usecounter{enualph}
    \setlength{\labelwidth}{1cm}}}
   {\end{list}}

\newcommand{\forces}{\Vdash}
\newcommand{\FORCES}{\left\|\!-\right.}

\newcommand{\embed}{<\!\!\circ\;}
\newcommand{\Bool}{[\![}   
\newcommand{\Boor}{]\!]}   
\newcommand{\BOOL}{\left[\!\!\left[}   
\newcommand{\BOOR}{\right]\!\!\right]}   

\newcommand{\A}{{\cal A}}
\newcommand{\B}{{\cal B}}

\newcommand{\F}{{\cal F}}

\newcommand{\I}{{\cal I}}

\newcommand{\M}{{\cal M}}
\newcommand{\N}{{\cal N}}

\newcommand{\U}{{\cal U}}

\renewcommand{\AA}{{\mathbb A}}
\newcommand{\BB}{{\mathbb B}}
\newcommand{\CC}{{\mathbb C}}
\newcommand{\DD}{{\mathbb D}}
\newcommand{\EE}{{\mathbb E}}

\newcommand{\MM}{{\mathbb M}}

\newcommand{\PP}{{\mathbb P}}
\newcommand{\QQ}{{\mathbb Q}}
\newcommand{\RR}{{\mathbb R}}

\newcommand{\TT}{{\mathbb T}}

\renewcommand{\aa}{{\mathfrak a}}
\newcommand{\bb}{{\mathfrak b}}
\newcommand{\cc}{{\mathfrak c}}
\newcommand{\dd}{{\mathfrak d}}

\newcommand{\mm}{{\mathfrak m}}

\newcommand{\rr}{{\mathfrak r}}

\newcommand{\sss}{{\mathfrak s}}

\newcommand{\add}{{\mathsf{add}}}
\newcommand{\cov}{{\mathsf{cov}}}
\newcommand{\non}{{\mathsf{non}}}
\newcommand{\cof}{{\mathsf{cof}}}
\newcommand{\Add}[1]{{{\mathsf{add}}({\cal #1})}}
\newcommand{\Cov}[1]{{{\mathsf{cov}}({\cal #1})}}
\newcommand{\Non}[1]{{{\mathsf{non}}({\cal #1})}}

\newcommand{\dom}{{\mathrm{dom}}}
\newcommand{\amal}{{\mathrm{amal}}}

\newcommand{\ro}{{\mathrm{r.o.}}}
\newcommand{\dir}{{\mathrm{dir}}}
\newcommand{\pure}{{\mathrm{Pure}}}

\newcommand{\zero}{{\mathbb{0}}}
\newcommand{\one}{{\mathbb{1}}}

\newcommand{\sub}{\subseteq}
\newcommand{\sem}{\setminus}
\newcommand{\twoom}{2^\omega}
\newcommand{\twolom}{2^{<\omega}}
\newcommand{\omlom}{\omega^{<\omega}}
\newcommand{\omom}{\omega^\omega}
\newcommand{\omloms}{[\omega]^{<\omega}}
\newcommand{\omoms}{[\omega]^\omega}

\newcommand{\la}{\langle}
\newcommand{\ra}{\rangle}
\newcommand{\Loriar}{\mbox{$\Longrightarrow$}}

\newcommand{\noi}{\noindent}

\newcommand{\emp}{\emptyset}

\title{Shattered Iterations}

\author{J\"org Brendle\thanks{Partially supported by Grant-in-Aid for Scientific Research
   (C) 18K03398, Japan Society for the Promotion of Science. \newline
     \indent {\it 2020 Mathematics Subject Classification.} Primary 03E35; Secondary 03E17, 03E40 \newline
    \indent {\it Keywords.} Lebesgue measure, Baire category, cardinal invariants, Cohen forcing, random forcing, iterated forcing  }  \\ 
   Graduate School of System Informatics \\
   Kobe University \\
   Rokko-dai 1-1, Nada-ku \\
   Kobe 657-8501, Japan \\   
   email: {\sf brendle@kobe-u.ac.jp}}

\begin{document}
\maketitle

\begin{abstract}
\noi We develop iterated forcing constructions dual to finite support iterations in the sense that
they add random reals instead of Cohen reals in limit steps. In view of useful applications we focus in particular on two-dimensional
``random" iterations, which we call {\em shattered iterations}. As basic tools for such iterations we investigate several  concepts
that are interesting in their own right. Namely, we discuss correct diagrams, we introduce the amalgamated limit of cBa's, 
a construction generalizing both the direct limit and the two-step amalgamation of cBa's, we present a detailed account of cBa's
carrying finitely additive strictly positive measures, and we prove a general preservation theorem for  such cBa's in amalgamated limits. 
As application, we obtain new consistency results on cardinal invariants in Cicho\'n's diagram. For example,
we show the consistency of $\aleph_1 < \cov (\M) < \non (\M)$, thus answering an old question of A. Miller~\cite{Mi81}.
\end{abstract}



\section*{Introduction}

Lebesgue measurability and Baire category play a fundamental role across many areas of mathematics, the former providing the
foundation of integration and probability theory, the latter being crucial for the basics of functional analysis. Foundational problems about
the two concepts, which exhibit important similarities~\cite{Ox80}, have given a strong impetus to the development of set theory (see,
e.g., \cite{So70, Sh84, BJ95}), and it is aspects of such foundational issues we will be dealing with here. Let $\M$ and $\N$ denote the
ideals of meager and null sets on the real numbers $\RR$, respectively. The {\em covering number} $\cov (\M)$ of $\M$ is the least size of
a family of meager sets covering the real line, the {\em uniformity} $\non (\M)$ of $\M$ is the smallest cardinality of a non-meager set
of reals, and $\cov (\N)$ and $\non (\N)$ are defined analogously for the null ideal. All these cardinal invariants are uncountable and below
the size of the continuum $\cc$ (the statement $\cov (\M) \geq \aleph_1$ is just the classical Baire category theorem). Furthermore,
using {\em orthogonality} of $\M$ and $\N$ (that is, the fact that there is a comeager $G_\delta$ null set) as well as the
translation invariance of both ideals, the inequalities $\cov (\N) \leq \non (\M)$ and $\cov (\M) \leq \non (\N)$ are easy to establish
(this is Rothberger's Theorem, see~\cite[Theorem 2.1.7]{BJ95}). No further inequalities can be established in ZFC. 
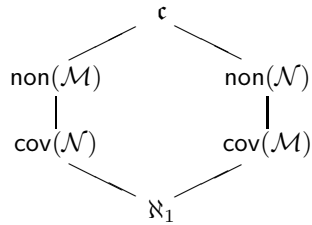
\[ 
\begin{picture}(80,90)(0,0)
\put(40,5){\makebox(0,0){$\aleph_1$}}
\put(0,30){\makebox(0,0){$\cov (\N)$}}
\put(80,30){\makebox(0,0){$\cov (\M)$}}
\put(0,55){\makebox(0,0){$\non (\M)$}}
\put(80,55){\makebox(0,0){$\non (\N)$}}
\put(40,80){\makebox(0,0){$\cc$}}
\put(0,36){\line(0,1){12}}
\put(80,36){\line(0,1){12}}
\put(5,23){\line(2,-1){25}}
\put(45,10){\line(2,1){25}}
\put(5,62){\line(2,1){25}}
\put(45,75){\line(2,-1){25}}
\end{picture}
\]
Indeed, letting $\B_\kappa$, $\M_\kappa$, and $\N_\kappa$ denote the Baire sets, the meager ideal, and the null ideal on $2^\kappa$, respectively,
where $\kappa$ is any cardinal, we define the {\em Cohen algebra} $\CC_\kappa =\B_\kappa / \M_\kappa$ and the {\em random algebra}
$\BB_\kappa = \B_\kappa / \N_\kappa$. For $\kappa = \omega$, we simply write $\B = \B_\omega$ for the Borel sets, and $\CC = \CC_\omega$ as well as $\BB = \BB_\omega$.
The impact of these forcing notions on measure and category was first investigated in pioneering work by Solovay~\cite{So70}, and it turns out that,
assuming $\kappa^\omega = \kappa$, forcing with $\CC_\kappa$ gives a model for $\cov (\N) = \non (\M) = \aleph_1$ and $\cov (\M) = \non (\N)
= \kappa = \cc$ while forcing with $\BB_\kappa$ yields $\cov (\M) = \non (\N) = \aleph_1$ and $\cov (\N) = \non (\M)   = \kappa = \cc$
(see~\cite{Ku84} for details). It is also known that $\cov (\N) < \non (\M)$ and $\cov (\M) < \non (\N)$ are consistent~\cite[Chapter 7]{BJ95}.

We will address the question to what extent these cardinal invariants may assume arbitrary values. This sort of
problem has gotten a lot of attention recently (see, in particular, ~\cite{GKS19, GKMSta1, KST19, KTT18}). Let $\kappa < \lambda$ be regular cardinals. Using a
finite support iteration (fsi, for short) of ccc forcing, a method introduced in seminal work of Solovay and Tennenbaum~\cite{ST71}, it is not difficult to
obtain a model with $\cov (\N) = \non (\M) = \kappa$ and $\cov (\M) = \non (\N) = \lambda$. Namely, assuming GCH in the ground model for simplicity, this is
achieved by performing a $\lambda$-stage fsi of subalgebras  of $\BB$ of size $< \kappa$ such that all such subalgebras are taken care of by a book-keeping argument
(see~\cite{Mi81} and~\cite[Model 7.5.1]{BJ95}). The dual consistency, namely CON($\aleph_1 < \kappa = \cov (\M) < \lambda = \non (\M)$), has been
an old problem of A. Miller~\cite[Problem (11), p. 348]{Mi81}. Why is this much much harder? --- The reason is that the standard iteration techniques for forcing, fsi mentioned
above, and countable support iteration of proper forcing, created by Shelah~\cite{Sh98}, won't work, the former because it adds Cohen reals in limit stages
and therefore necessarily forces $\non (\M) \leq \cov (\M)$, and the latter, because it can only give models for $\cc \leq \aleph_2$.

To appreciate how this problem might be solved, let us try to dualize the consistency proof of $\non (\M) = \kappa$ and $\cov (\M) = \lambda$ mentioned above,
by interchanging the roles of measure and category. The heuristic we obtain says we should set up an iterated forcing construction going through subalgebras of $\CC$ of size $< \kappa$ and adding
random reals instead of Cohen reals in limit stages. The former is clearly nonsense because $\CC$ is countable and so are all its subalgebras which are then
again isomorphic to $\CC$, and while it is a priori unclear what the latter should mean, the development of such ``random" iterations will turn out
to be one of the main points of our work. So let us shift our perspective. Notice that $\CC_\lambda$ completely embeds into the
fsi $\PP$ forcing $\kappa = \non (\M) < \cov (\M) = \lambda$. That is, we may write $\PP$ as a two-step iteration $\CC_\lambda \star \dot \QQ$ where 
the remainder forcing $\dot \QQ$ is itself an iteration adding partial random reals. So the revised heuristic says that we should start with $\BB_\lambda$
and iterate adjoining reals which are Cohen only over a small fragment of the $\BB_\lambda$-extension. This works indeed as we shall see:

\begin{maintheoremA}
Let $\kappa < \lambda$ be arbitrary regular uncountable cardinals. Then there is a ccc forcing extension satisfying $\cov (\M) = \non (\N) = \kappa$ and
$\cov (\N) = \non (\M) = \lambda$.
\end{maintheoremA}

The proof needs quite a bit of preliminary work, which is interesting in its own right. In particular, we shall discuss {\em correct diagrams}, a 
fundamental tool when dealing with nonlinear iterations, in detail in Section~\ref{correct}. Section~\ref{amalgamated} introduces and develops the basics
of the {\em amalgamated limit}, a novel limit construction for iterated forcing encompassing the direct limit and the two-step amalgamation as
special cases. In section~\ref{coherent} we use {\em finitely additive strictly positive measures} on complete Boolean algebras (cBa's) to establish
a general preservation theorem for the amalgamated limit construction. Equipped with these tools we prove Main Theorem A in Section~\ref{Cohen}.
Since the type of iteration we are using gradually destroys the homogeneity of the initial $\BB_\lambda$-extension by introducing more and more partial Cohen
(or other) generics, we call such iterations {\em shattered iterations}. We briefly discuss in Section~\ref{problems} to what extent shattered
iterations can be thought of as two-dimensional iterations adding random reals in limit stages.

There are many other cardinal invariants of interest in set theory and its applications, and two of the most prominent and useful are the
{\em unbounding number} $\bb$ and the {\em dominating number} $\dd$. For functions $f, g$ in the Baire space $\omom$ say that
$g$ {\em eventually dominates} $f$ ($f \leq^* g$ for short) if $f(n) \leq g(n)$ holds for all but finitely many $n$. A family $\F \sub \omom$
is {\em unbounded} if no $g \in \omom$ eventually dominates all members of $\F$, and $\F$ is {\em dominating} if for all $g \in \omom$
there is $f \in \F$ eventually dominating $g$. $\bb$ ($\dd$, respectively), is the least size of an unbounded family (of a dominating family, resp.).
It is easy to see that $\bb \leq \dd$ are uncountable cardinals below $\cc$, and since non-meager subsets of $\omom$ are unbounded,
$\bb \leq \non (\M)$ and $\cov (\M) \leq \dd$ follow (see~\cite[Proposition 5.5]{Bl10} for details). Unfortunately, we do not know the values of $\bb$ and $\dd$
in the model of Main Theorem A (though we strongly suspect them to be $\aleph_1$ and $\cc$, respectively, see Conjecture~\ref{b-d-conjecture}
and the discussion afterwards in Section~\ref{preservation}). By modifying the model of Main Theorem A, we prove in Section~\ref{Hechler}:

\begin{maintheoremB}
Let $\kappa < \lambda$ be arbitrary regular uncountable cardinals. Then there is a ccc forcing extension satisfying $\bb = \cov (\M) = \non (\N) = \kappa$ and
$\dd = \cov (\N) = \non (\M) = \lambda$.
\end{maintheoremB}

Section~\ref{preservation} contains a brief and still unsatisfactory discussion on preservation results for shattered iterations.

\bigskip

\noindent \underline{Notation}. We assume familiarity with basic techniques of set theory~\cite{Je03, Ku11} and, in particular,
with forcing theory. For cardinal invariants see~\cite{BJ95, Bl10}. For technical reasons, we shall formulate most of our forcing-theoretic
results in the language of complete Boolean algebras (cBa's, for short) $\la \AA, + , \cdot , \zero ,\one \ra$. In partial orders (p.o.'s for short) or cBa's,
$\bot$ denotes incompatibility of forcing conditions. For a p.o.  $\PP$, $\ro (\PP)$, the regular open algebra of $\PP$ is the complete
Boolean algebra associated with $\PP$ (so $\PP$ and $\ro (\PP)$ are forcing equivalent). In the forcing language of cBa's we use
$\Bool \varphi \Boor$ for Boolean values of formulas $\varphi$. Given cBa's $\AA_0$ and $\AA$, a mapping $e : \AA_0 \to \AA$ is a complete embedding
if it is an injection preserving infinite sums and products in addition to the Boolean operations. We shall usually assume $e$ is the identity $id$. In this
case we write $\AA_0 \embed \AA$ for ``$\AA_0$ is completely embedded in $\AA$". If $\AA_0$ is a cBa, and $\dot \AA_1$ is an $\AA_0$-name
for a cBa, $\AA : = \AA_0 \star \dot \AA_1$ denotes the two-step iteration. In this situation $\AA_0 \embed \AA$. On the other hand, if $\AA_0 \embed \AA$,
then there is an $\AA_0$-name $\dot \AA_1$ for a cBa such that $\AA = \AA_0 \star \dot \AA_1$; namely, letting $G_0$ be $\AA_0$-generic over the ground model  $V$,
we may define the quotient cBa (also referred to as remainder forcing) $\AA / G_0$ in $V [G_0]$ and then have that $\dot \AA_1 [G_0] = \AA / G_0$.
If $\la \AA_i : i \in I \ra$ is a system of cBa's with complete embeddings $id: \AA_i \embed \AA_j$ for $i < j$ where $\la I, \leq \ra$
is a directed p.o., we let $\lim \dir_{i \in I} \AA_i : = \ro (\bigcup_{i\in I} \AA_i )$ denote its direct limit.



\section{Correct diagrams}
\label{correct}

Let $I = \{ 0 \land 1 , 0, 1, 0 \lor 1 \}$ be the four-element lattice with top element
$0 \lor 1$, bottom element $0 \land 1$, and $0$ and $1$ in between. Assume we are
given cBa's $\AA_i$ for $i \in I$ and complete embeddings $e_{ij} : \AA_i \to \AA_j$
for $i \leq j$ from $I$, that is, we consider the diagram
\[ 
\begin{picture}(80,50)(0,0)
\put(48,5){\makebox(0,0){$\AA_{0\land 1}$}}
\put(7,25){\makebox(0,0){$\bar \AA : \;\;\; \AA_{0}$}}
\put(68,25){\makebox(0,0){$\AA_{1}$}}
\put(48,45){\makebox(0,0){$\AA_{0\lor 1}$}}
\put(26,28){\line(1,1){13}}
\put(28,20){\line(1,-1){11}}
\put(51,8){\line(1,1){13}}
\put(53,40){\line(1,-1){11}}
\end{picture}
\]
of cBa's. The $e_{ij}$ give rise to projection mappings $h_{ji} : \AA_j \to \AA_i$
naturally defined by
\[ h_{ji} (a_j) = \prod \{ a \in \AA_i : e_{ij} (a) \geq a_j \} \]
To make life easier let us assume that $\AA_i \subseteq \AA_j$ for $i \leq j$
and that $e_{ij}$ is the identity.  In this situation $\AA_{0 \land 1} = \AA_0 \cap \AA_1$ is a natural requirement,
but we shall need a stronger property.

\begin{defin}
The diagram $\bar \AA = ( \AA_i : i \in I )$ is {\em correct} if for all $a_0 \in \AA_0$ and $a_1 \in \AA_1$, if
$h_{0,0\land 1} (a_0) \cdot h_{1,0\land 1} (a_1) \neq \zero$ then $a_0 \cdot a_1 \neq \zero$ in $\AA_{0 \lor 1}$.
\hfill $\dashv$
\end{defin}

\begin{exam} ({\em product})   \label{exam1}
Let $\AA_{0 \land 1} = \{ \zero, \one \}$ be the trivial algebra and let $\AA_{0 \lor 1}$
be the cBa generated by the product forcing $( \AA_0 \sem \{ \zero \}) \times ( \AA_1 \sem \{ \zero \} )$.
Clearly the resulting diagram is correct. \hfill $\dashv$
\end{exam}

\begin{exam} ({\em amalgamation})   \label{exam2}
Let $\AA_{0\lor 1} = \amal_{\AA_{0\land 1}} (\AA_0 , \AA_1)$ be the {\em amalgamation of $\AA_0$
and $\AA_1$ over $\AA_{0 \land 1}$.} Conditions are pairs $(a_0 , a_1 ) \in  
( \AA_0 \sem \{ \zero \}) \times ( \AA_1 \sem \{ \zero \} )$ such that $h_{0, 0 \land 1} (a_0) \cdot
h_{1, 0\land 1} (a_1) \neq \zero$, equipped with the coordinatewise order,
$(a_0 ' , a_1 ') \leq (a_0 , a_1)$ if $a_0 ' \leq a_0$ and $a_1 ' \leq a_1$.
$\AA_{0 \lor 1 }$ is the completion of this set of conditions. Identify $\AA_0$ with 
$\{ (a_0 , \one) : a_0 \in \AA_0 \}$; similarly for $\AA_1$; thus $\AA_0 , \AA_1 \sub \AA_{0 \lor 1}$,
and it is well-known the embeddings are complete.
The diagram is correct by definition of the amalgamation: if $h_{0,0\land1} (a_0 ) \cdot
h_{1, 0\land 1} (a_1) \neq \zero$, then $(a_0 , a_1) \neq \zero$ is a common extension. 

Notice that the set of conditions $(a_0,a_1)$ thus defined is not separative.
By considering only pairs $(a_0,a_1)$ such that $h_{0, 0 \land 1} (a_0) = h_{1, 0 \land 1} (a_1)$, however,
we get an equivalent separative reformulation (we will put this fact to good use below, in Definition~\ref{defin3}; see
also Observation~\ref{basicobs} (iii)).

Note that the product is a special case of the amalgamation with $\AA_{0 \land 1} = \{ \zero, \one \}$.
In fact, there is another way of looking at the amalgamation: forcing with $\amal_{\AA_{0\land 1}} (\AA_0 , \AA_1)$
is equivalent to first forcing with $\AA_{0 \land 1}$ and then with the product of
the quotients $\AA_0 / G_{0\land 1} $ and $\AA_1 / G_{0\land 1}$ where $G_{0 \land 1}$ is the
$\AA_{0 \land 1}$-generic filter. \hfill $\dashv$
\end{exam}

Amalgamations have been little used so far in forcing theory, the most prominent example being 
those occurring in Shelah's consistency proof of the projective Baire property on the basis of the consistency of ZFC~\cite{Sh84}.

\begin{exam} ({\em two-step iteration of definable forcing})  \label{exam3}
Recall {\em Hechler forcing} $\DD$ consists of pairs
$(s,f) \in \omlom \times \omom$ such that $s \sub f$. The ordering is given by $(t,g) \leq (s,f)$
if $t \supseteq s$ and $g(n) \geq f(n)$ for all $n \in \omega$. Hechler forcing is $\sigma$-centered and generically adds
a real $d \in \omega^\omega$ which eventually dominates all ground model reals (that is, for all
$f \in \omom \cap V$, for all but finitely many $n$, $d(n) > f(n)$). Let $\AA_0 = \AA_1 =
\ro (\DD)$, the cBa generated by $\DD$, put $\AA_{0\lor 1} = \AA_0 \star \dot \AA_1$, the two-step iteration
(where $\dot \AA_1$ is an $\AA_0$-name for $\ro (\dot\DD)$, Hechler forcing in the sense of $V^{\AA_0}$),
and let $\AA_{0 \land 1} = \{ \zero, \one \}$ be trivial. As before identify $\AA_0$ with $\{ (a_0 , \one) :
a_0 \in \AA_0 \}$; similarly for $\AA_1$. 

Then all embeddings are complete, the only nontrivial case being $\AA_1 \embed \AA_{0 \lor 1}$:
we need to argue that every maximal antichain $A \sub \AA_1$ is still maximal in $\AA_0 \star \dot \AA_1$.
By the above identification, this is clearly equivalent to saying that every maximal
antichain $A \sub \DD$ in the ground model $V$ is still a maximal antichain of $\DD$ in the sense
of $V [G_0]$ where $G_0 \sub \AA_0$ is an arbitrary generic. (Note here that $\DD^V$ is a proper
subset of $\DD^{V[G_0]}$.) However, since $\DD$ is Borel ccc ($\DD$, the ordering and incompatibility
are Borel), being a maximal antichain in $\DD$ is a $\Pi^1_1$ statement and therefore absolute.
Correctness, then, is straightforward.\hfill $\dashv$
\end{exam}

Note that absoluteness was only relevant for the second step of the two-step iteration $\AA_{0\lor 1} = \AA_0 \star \dot \AA_1$ in this example.
That is, $\AA_0$ can be replaced by any cBa, while $\AA_1$ can be any sufficiently definable forcing, e.g.,
random forcing $\BB$ or a large random algebra $\BB_\kappa$.
Such absoluteness considerations will be used throughout our work for Borel ccc forcing notions. 

\begin{obs}[basic properties of correct diagrams]~   \label{correct-basic}
\begin{enumerate}
\item The diagram $\bar \AA$ is correct iff for all $a_1 \in \AA_1$, $h_{0 \lor 1 , 0} (a_1) = h_{1 , 0\land 1} (a_1)$.
\item $\bar \AA$ is correct iff $\bar \EE$ is correct for every $\AA_{0 \land 1}$-generic filter $G_{0\land 1}$ where $\EE_i = \AA_i / G_{0 \land 1}$ 
    (so $\EE_{0\land 1} = \{ \zero , \one \}$).
\item If $\bar \AA$ is correct, $\AA_{0\lor 1} \embed \AA_{0\lor 1}'$ and $\bar \AA'$ is obtained from $\bar \AA$ by replacing
   $\AA_{0\lor 1}$ by $ \AA_{0\lor 1}'$ and keeping the other cBa's, then $\bar \AA'$ is correct as well. 
\item If $\bar \AA$ is correct then $\AA_{0 \land 1} = \AA_0 \cap \AA_1$.
\end{enumerate}
\end{obs}

\begin{proof} (i) First assume $h_{0 \lor 1, 0} (a_1) < h_{ 1, 0 \land 1} (a_1)$ for some 
$a_1 \in \AA_1$. Then there is $a_0 \in \AA_0$ incompatible with $h_{0 \lor 1, 0} (a_1)$
such that $h_{0, 0 \land 1} (a_0) \leq h_{1, 0 \land 1} (a_1)$. Thus
$h_{0, 0 \land 1} (a_0) \cdot h_{1, 0 \land 1} (a_1) \neq \zero$, yet
$a_0 \cdot a_1 = \zero$ in $\AA_{0 \lor 1}$.
Suppose, on the other hand, that $h_{0 \lor 1, 0} (a_1) = h_{ 1, 0 \land 1} (a_1)$ for all $a_1 \in \AA_1$.
If $h_{0, 0 \land 1} (a_0) \cdot h_{1, 0 \land 1} (a_1) \neq \zero$
in $\AA_{0 \land 1}$, then $a_0 \cdot h_{0 \lor 1, 0 } (a_1) =
a_0 \cdot h_{1, 0 \land 1} (a_1) \neq \zero$ in $\AA_0$, and thus $a_0 \cdot a_1 \neq \zero$ in $\AA_{0 \lor 1}$.

(ii) Indeed if $a_0 \cdot a_1 = \zero$ but $p = h_{0,0\land 1} (a_0) \cdot h_{1, 0 \land 1} (a_1)
\neq \zero$, take a generic $G_{0\land 1}$ containing $p$ to see that $e_0 \cdot e_1 = \zero$
where $e_i = a_i / G_{0\land 1} \neq \zero$. 
If, on the other hand, there are a generic $G_{0 \land 1}$
and $e_0 \in \EE_0, e_1 \in \EE_1, e_i \neq \zero$ such that $e_0 \cdot e_1 = \zero$, then,
for some $p \in G_{0 \land 1}$, $p \forces_{\AA_{0 \land 1}} `` \dot e_i \neq \zero, \dot e_0 \cdot \dot e_1
= \zero "$. Thus, letting $a_i = (p, \dot e_i) \in \AA_i = \AA_{0\land 1} \star \dot \EE_i$, 
we get $h_{0,0\land 1} (a_0) = h_{1, 0 \land 1} (a_1) = p
\neq \zero$, yet $a_0 \cdot a_1 = \zero$.

(iii) Obvious.

(iv) If $\AA' $ is a strict complete subalgebra of $\AA_0 \cap \AA_1$, then there are incompatible $a , a' \in \AA_0 \cap \AA_1$
that project to the same element of $\AA '$.
\end{proof}

Naturally occurring diagrams satisfying (iv) are correct, but correctness is stronger than merely saying
$\AA_{0 \land 1} $ is $\AA_0 \cap \AA_1$.

\begin{counterexam} \label{counterexam1}
Let $\AA_{0 \lor 1}$ be the forcing which adds two pairwise disjoint
Cohen reals. Thus $\AA_{0\lor 1} = \ro (\PP)$ where $\PP$ consists of pairs $(s,t) \in 
(\twolom)^2$ such that $s^{-1} (\{ 1 \} ) \cap t^{-1} (\{ 1 \}) = \emp$. Set $\AA_0
= \AA_1 = \CC$ and $\AA_{0 \land 1} = \{ \zero, \one \}$.
Note that $\AA_0 \embed \AA_{0 \lor 1}$ and similarly for $\AA_1$. 
(To see this let $A \sub \AA_0$ be a maximal antichain, and let $(s,t) \in \AA_{0 \lor 1}$.
By extending if necessary, we may assume without loss of generality that $|s| \geq |t|$. Choose
$s_0 = (s_0 , \emp ) \in A$ compatible with $s$. Then $(s \cup s_0 , t)$ is a common extension
of $s_0$ and $(s,t)$.)

Next we argue $\AA_0 \cap \AA_1 = \{ \zero, \one \}$. For assume $a \in \AA_0 \cap \AA_1$ is
nontrivial. There are $s_n , t_n \in \twolom$ such that $a = \sum_n (s_n, \emp) = \sum_n
s_n = \sum_n (\emp , t_n) = \sum_n t_n$.
Choose $s \in \twolom$ incompatible with all $s_n$. So $s \in \AA_0$ is incompatible with all $t_n$, i.e.,
for all $n$ there is $i \in |s| \cap |t_n|$ such that $s(i) = t_n (i) = 1$. Thus defining
$t$ by $|t| = |s|$ and $t(i) = 0$ for all $i \in |t|$, $t \in \AA_1$ is incompatible with
all $t_n$. Yet, $t$ is compatible with all conditions in $\AA_0$, a contradiction.

Finally notice that, letting $t = \{ \la 0,1 \ra \} \in \AA_1$,
$h_{1, 0 \land 1} (t) = \one$ while $h_{0 \lor 1, 0} (t) < \one$ because $t$ is incompatible
with $s = \{ \la 0,1 \ra \} \in \AA_0$. So the diagram is not correct.
\hfill $\dashv$
\end{counterexam}

An important property of correctness is its preservation under iteration with definable forcing.
For example:

\begin{lem}[preservation of correctness in two-step iterations with Hechler forcing $\DD$, {\cite[Lemma 28]{survey}}]  \label{Hechler-correct}
Assume $\bar \AA$ is correct. Then so is $\bar \EE$ where $\EE_i = \AA_i \star \ro (\dot \DD)$.
\end{lem}

We omit the proof since we will not need this. We shall need, however, the following:

\begin{lem}[preservation of correctness in two-step iterations with random forcing $\BB$]   \label{random-correct}
Assume $\bar \AA$ is correct. Then so is $\bar \EE$ where $\EE_i = \AA_i \star \dot \BB$
or $\EE_i = \AA_i \star \dot \BB_\kappa$.  
\end{lem}

\begin{proof} Again assume that $\AA_{0 \land 1} = \{ \zero, \one \}$ (Observation~\ref{correct-basic} (ii)).
We only show the first part (in fact, there is no need to care about $\kappa$ at all because
any condition in $\BB_\kappa$ is a condition of $\BB_C$ for some countable $C \sub \kappa$).
Thus $\EE_{0\land 1} = \BB$. 

Let $e_0 = (p_0 , \dot b_0) \in \EE_0 = \AA_0 \star \dot \BB$ and $e_1 = (p_1 , \dot b_1) \in \EE_1 = \AA_1 \star \dot \BB$
such that $h^\EE_{0, 0\land 1} (e_0) \cdot h^\EE_{1, 0\land 1} (e_1) \neq \zero$. Strengthening $e_0$, we may assume 
that $h^\EE_{0, 0\land 1} (e_0) \leq  h^\EE_{1, 0\land 1} (e_1)$. 
By Lebesgue density and a further strengthening of $e_0$ we may assume that for
$b_{0\land 1} : = h^\EE_{0, 0\land 1} (e_0)$, $p_0$ forces in $\AA_0$ that $\mu (\dot b_0 \cap b_{0\land 1} ) > {3 \over 4}
\mu (b_{0\land 1})$. 

Again by Lebesgue density we can find a maximal antichain of conditions $c^n \in \BB = \EE_{0\land 1}$,
$n\in\omega$, below $b_{0\land 1}$ and conditions $p^n_1$, $n\in\omega$, below $p_1$ in $\AA_1$ such that
$p_1^n$ forces in $\AA_1$ that  $\mu (\dot b_1 \cap c^n) \geq {3 \over 4} \mu (c^n)$.
Note there must be $n \in\omega$ and $q_0 \leq p_0$ such that
$q_0 \forces_{\AA_0} \mu (\dot b_0 \cap c^n) \geq {3 \over 4} \mu (c^n)$.
For if this were not the case, $p_0 \forces_{\AA_0} \mu (\dot b_0 \cap  c^n) < {3 \over 4} \mu ( c^n)$ for all $n$.
So $p_0 \forces_{\AA_0} \mu (\dot b_0 \cap b_{0 \land 1} ) = \sum_n \mu (\dot b_0 \cap  c^n)
< \sum_n  {3 \over 4} \mu (c^n) =  {3 \over 4} \mu ( b_{0 \land 1})$, a contradiction.

Let $q_1 = p_1^n$. By correctness of the diagram $\bar \AA$, $q_0$ and $q_1$ are compatible in $\AA_{0 \lor 1}$
with common extension $q_{0\lor 1} = q_0 \cdot q_1$. Let $\dot c_{ 0 \lor 1} = c^n \cap \dot b_0 \cap \dot b_1$,
an $\AA_{0 \lor 1}$-name for a condition in $\dot \BB$. We easily see that $q_{0\lor 1} \forces_{\AA_{0\lor 1}} 
\mu (\dot c_{0 \lor 1}) \geq {1 \over 2} \mu (c^n) > 0$. So $e_{0 \lor 1} = (q_{0 \lor 1} , \dot c_{0 \lor 1})$ belongs to 
$\EE_{0 \lor 1} \setminus \{ \zero \}$ and is a common extension of $e_0$ and $e_1$.
\end{proof}

We actually believe this is true for a large class of forcing notions.

\begin{conj}[{\cite[Conjecture 1.5]{luminy}}] Let $\PP$ be a Borel (or Suslin) ccc forcing and assume $\bar \AA$ is correct.
Then so is $\bar \EE$ where $\EE_i = \AA_i \star \ro (\dot \PP)$.
\end{conj}

An important application of correctness is:

\begin{lem}[embeddability of direct limits, {\cite[Lemma 3]{survey}}] \label{correct-limit}
Let $K$ be a directed index set.
Assume $\la \AA_k : k \in K \ra$ and $\la \EE_k : k \in K \ra$ are systems of cBa's such that
$\AA_k \embed \AA_\ell$, $\EE_k \embed \EE_\ell$ and $\AA_k \embed \EE_k$ for any
$k \leq \ell$. Assume further all diagrams
\begin{picture}(60,50)(0,0)
\put(28,3){\makebox(0,0){$\AA_k$}}
\put(5,23){\makebox(0,0){$ \EE_k$}}
\put(48,23){\makebox(0,0){$\AA_\ell$}}
\put(28,43){\makebox(0,0){$\EE_\ell$}}
\put(8,26){\line(1,1){13}}
\put(10,18){\line(1,-1){11}}
\put(31,6){\line(1,1){13}}
\put(33,38){\line(1,-1){11}}
\end{picture}
are correct for $k \leq \ell$. Then $\lim \dir_k \AA_k \embed \lim \dir_k \EE_k$.
\end{lem}

\begin{proof}
Let $A \sub \bigcup_k \AA_k$ be a maximal antichain in $\AA : = \lim \dir_k \AA_k$. Choose
$b \in \EE : = \lim \dir_k \EE_k$, i.e. $b \in \EE_k$ for some $k$. By maximality of $A$
there is $a \in A$ such that $h_{\EE_k , \AA_k} (b) \cdot h_{\AA, \AA_k} (a) \neq \zero$.
Find $\ell \geq k$ such that $a \in \AA_\ell$. By correctness, $b \cdot a \neq \zero$
in $\EE_\ell$ as required.
\end{proof}

Again, embeddability of direct limits is far from being true in general.

\begin{counterexam}
Let $\AA_n$ add $n$ Cohen reals and let $\EE_n$ add $n$ Cohen reals and an additional one
which is disjoint from all of them. Clearly $\AA_m \embed \AA_n , \EE_m \embed \EE_n$, for
$m \leq n < \omega$, and $\AA_n \embed \EE_n$ by Counterexample~\ref{counterexam1}.
On the other hand, $\AA : = \lim \dir_n \AA_n \not{\!\!\embed} \lim \dir_n \EE_n = : \EE$
because, denoting the $\omega$ Cohen reals added by the sequence of $\AA_n$'s by $\la c_n : n\in \omega \ra$,
$\bigcup_n c_n = \omega$ in the $\AA$-generic extension while $\bigcup_n c_n$
is coinfinite in the $\EE$-generic extension. 
\hfill $\dashv$
\end{counterexam}

While it is not very relevant for linear iterations, correctness is a crucial notion
when it comes to building nonlinear iterations. In fact, the author originally introduced this notion in~\cite{luminy}
when trying to grasp Shelah's technique of iteration along templates
(see~\cite{Sh700}, \cite{madit} and~\cite{adcc}). For a rephrasing of this technique in the
language of complete embeddings with correct projections, using Lemmata~\ref{Hechler-correct}
and~\ref{correct-limit} above, see either~\cite{luminy} or the recent survey~\cite{survey}.
In the template context, correctness has been used in other work since, see e.g.
the detailed account in~\cite[Section 2]{FM17}.

It should be pointed out, however, that correctness can be useful in classical finite support iteration too.
For example, if all iterands $\dot \QQ_\alpha$ are Suslin ccc and coded in the ground model (see~\cite[Section 3.6]{BJ95}), 
we obtain from Lemma~\ref{correct-limit} immediately that  any subiteration using only some of the $\dot \QQ_\alpha$
and the trivial forcing elsewhere completely embeds into the whole iteration.



\section{The amalgamated limit}  
\label{amalgamated}

The purpose of this section is to present an iterated forcing construction which generalizes both the
direct limit and the (two-step) amalgamation described in Example~\ref{exam2} above.

\begin{defin}   \label{defin2}
Assume $\la I, \leq \ra$ is a partial order. Call $I$ a {\em distributive almost-lattice} if
\begin{enumerate}
\item any two elements $i,j \in I$ have a greatest lower bound $i \land j$, the {\em meet}
   of $i$ and $j$,
\item if $i$ and $j$ have an upper bound, then they have a least upper bound $i \lor j$, the
   {\em join} of $i$ and $j$,
\item if $i_0, i_1, i_2 \in I$ then there are $j \neq k$ such that $i_j $ and $i_k$ have an upper bound
   (so that $i_j \lor i_k$ exists), that is, given any three elements, two of them have an upper bound,
\item the distributive laws hold for $\land$ and $\lor$ (as long as both sides of the law in question 
   exist\footnote{One may have, e.g., that $i$ and $j$ have no upper bound, so the expression
   $(i \lor j) \land k$ makes no sense while $(i \land k) \lor (j \land k)$ is below $k$ and thus always exists.}),
\item if $i$ and $j$ have no upper bound and $j'$ is arbitrary, then \underline{either} neither $i, j'$ nor 
   $i, j\land j'$ have an upper bound \underline{or} $i \lor j' = i \lor (j \land j')$, and
\item if $j,j'$ have no upper bound, then $(i \land j) \lor (i \land j') = i$ for any $i$.  \hfill $\dashv$
\end{enumerate}
\end{defin}

The only reason for items (v) and (vi) in this definition is the following:

\begin{obs}  \label{almost-lattice-basic}
Assume $\la I,  \leq \ra$ satisfies (i) to (iv). Also assume $I$ has no maximal element, let $\ell \notin I$, $J = I \cup \{ \ell \}$, and
stipulate $i \leq \ell$ for all $i \in I$. Then the following are equivalent:
\begin{itemize}
\item $I$ satisfies (v) and (vi).
\item The distributive laws hold in $J$, that is, $J$ is a lattice.
\end{itemize}
\end{obs}

\begin{proof}
(v) clearly guarantees the distributive law $i \lor (j \land j') = (i \lor j) \land (i \lor j')$ $\;(\star)$. 
On the other hand, if (v) fails, then \underline{either} $i,j'$ has no upper bound in $I$ while $i, j\land j'$ does in $I$
in which case the left-hand side of $(\star)$ belongs to $I$ while the right-hand side is $\ell$,
\underline{or} $i \lor (j\land j') < i \lor j'$ in $I$, in which case the left-hand side of $(\star)$ is
still strictly smaller than the right-hand side in $J$.

(vi) clearly guarantees the distributive law $i \land (j \lor j') = (i \land j) \lor (i \land j')$ $\; (\star\star)$.
On the other hand, if $(i \land j) \lor (i \land j') < i$ for $j , j'$ without upper bound in $I$, then the right-hand side is strictly smaller than
the left-hand side in $(\star\star)$.
\end{proof}

Typical examples are distributive lattices or the three-element almost-lattice $\{ 0, 1, 0 \land 1 \}$.
We shall encounter a more interesting example in Section~\ref{Cohen}. The next result describes the basic structure of distributive
almost lattices\footnote{This structural result simplifies the definition of, and several proofs about, the amalgamated limit below. It was not
part of the first version of this paper nor of the survey~\cite{amallimit}, and this accounts for the differences between the latter work and the present one.}.

\begin{basiclem}   \label{basic-lemma}
Let $\la I, \leq \ra$ be a distributive almost-lattice. Then there are the {\em kernel} $\ker (I)$ as well as the {\em left part} $L = L(I)$ and
the {\em right part} $R = R(I)$ such that 
\begin{enumerate}
\item $i \in \ker (I)$ iff $i \lor j$ exists for all $j \in I$,
\item $L$ and $R$ are downward closed, i.e., $i \in L$ and $j \leq i$ imply $j \in L$ (same for $R$),
\item if $i,j \in L$ then $i \lor j$ exists and belongs to $L$ (same for $R$),
\item $I = L \cup R$, $\ker (I) = L \cap R$, and $L \sem \ker (I) \neq \emptyset $ iff $R \sem \ker (I) \neq \emptyset $,
\item if $i \in L \sem \ker (I)$ and $j \in R \sem \ker (I)$, then $i \lor j$ does not exist,
\item if $i \in L$ and $j \in R$ then $i \land j \in \ker (I)$.
\end{enumerate}
\end{basiclem}

\begin{proof}
Define $\ker (I) := \{ i : \forall j \in I \; ( i \lor j$ exists$)\}$. Then (i) clearly holds. If $\ker (I) = I$ let $L = R = I$, and the remaining clauses
are trivially satisfied.

So assume $\ker (I) \neq I$, and choose $i_0 \in I \sem \ker (I)$. Let $L : = \{ j \in I : i_0 \lor j$ exists$\}$. Since $i_0 \notin \ker (I)$,
there is $j_0 \in I$ such that $i_0 \lor j_0$ does not exist. So $j_0 \notin L$. Let $R : = \{ k \in I : j_0 \lor k$ exists$\}$. 
Assume now $i \in L$ and $j \leq i$. Then $i_0 \lor i$ is an upper bound of $i_0$ and $j$ and therefore $i_0 \lor j$ exists so that $j \in L$.
This shows (ii).

To see (iii), take $i,j \in L$. So both $i_0 \lor i$ and $i_0 \lor j$ exist. On the other hand, both $(i_0 \lor i) \lor j_0$ and $(i_0 \lor j) \lor j_0$
cannot exist. Therefore $i_0 \lor i \lor j$ must exist and so does $i \lor j$. Similarly for $R$.

Assume $j \notin L$. So $i_0 \lor j$ does not exist. Therefore $j_0 \lor j$ must exist and $j$ belongs to $R$. Thus $I = L \cup R$. Furthermore
$\ker (I) \sub L \cap R$ is obvious. Assume $k \in L \cap R$, and let $i \in I$. If $i \in L$, $k \lor i$ exists by (iii), and similarly if $i \in R$. Hence
$k \in \ker (I)$, and (iv) is proved.

Let $i \in L \sem \ker (I)$ and $j \in R \sem \ker (I)$ and suppose $i \lor j$ exists. Without loss of generality (by (iv)) $i \lor j \in L$. By (ii),
$j \in L$ follows, a contradiction. Thus (v) holds, and (vi) is clear by (ii) and (iv). 
\end{proof}

\begin{defin}   \label{defin3}
Given a distributive almost-lattice $\la I , \leq \ra$ with minimum $0$ and a system $\la \AA_i : i \in I \ra$
of cBa's with complete embeddings $id : \AA_i \embed \AA_j$ for $i < j$ and with $\AA_0 = \{ \zero, \one \}$ such that
all diagrams of the form
\begin{picture}(50,50)(0,0)
\put(28,3){\makebox(0,0){$\AA_{i\land j}$}}
\put(-7,23){\makebox(0,0){$ \bar \AA : \;\;\; \AA_i $}}
\put(48,23){\makebox(0,0){$\AA_j$}}
\put(28,43){\makebox(0,0){$\AA_{i\lor j}$}}
\put(6,26){\line(1,1){13}}
\put(8,18){\line(1,-1){11}}
\put(31,6){\line(1,1){13}}
\put(33,38){\line(1,-1){11}}
\end{picture}
are correct, we define the {\em amalgamated limit} $\AA_\ell = \AA_\amal = \lim \amal_{i\in I} \AA_i$ as follows. First
set $\AA = \bigcup_{i\in I} \AA_i \sem \{ \zero \}$.
\begin{itemize}
\item The set $D$ of nonzero conditions consists of ordered pairs $(p,q) \in \AA \times \AA$
   such that there are $i\in L$ and $j \in R$ with $p \in \AA_i$, $q \in \AA_j$ and
   $h_{i, i \land j} (p) = h_{j, i\land j} (q)$ (note that $i\land j \in \ker(I)$ by (vi) of the basic lemma).
\item The order is given by stipulating $(p',q') \leq (p,q)$ if 
   $h_{i' , i' \land i} (p')  \leq_{i} p$ and $h_{j', j' \land j} (q') \leq_{j} q$ 
   where $(i',j')$ witnesses $(p',q') \in D$.
\end{itemize}
$\AA_\ell$ is the cBa generated by $D$, i.e., $\AA_\ell = \ro (D)$.  \hfill $\dashv$
\end{defin}

We usually omit the subscript when referring to the ordering and write $\leq$ instead of $\leq_i$
since it will always be understood in which cBa we work. Notice, however, that while we do not
require $i \leq i'$ ($j \leq j'$, resp.) in the definition of the partial order above, the latter assumption is
natural, often true in concrete situations, and makes the definition of the order simpler: 
\begin{itemize}
\item $(p',q') \leq (p,q)$ if $p' \leq_{i'} p$ and $q' \leq_{j'} q$
\end{itemize}
because $h_{i'i} (p') \leq_i p$ and $p' \leq_{i'} p$  are clearly equivalent. 
Furthermore, note that the order on $D$ is not separative, that is, distinct
members of $D$ are equivalent in the completion $\AA_\ell$. This will be further clarified in 
Observation~\ref{basicobs} below.
Moreover, we identify each $\AA_i$, $i \in L$, with a
subset of $D$ (and thus with a subalgebra of $\AA_\ell$) via $p \mapsto (p,\one)$ where $\one \in \AA_0$
(that is, $(i,0)$ witnesses $(p,\one) \in D$), and similarly for $\AA_j$, $j \in R$. Note that we always have
$(p,q) \leq (p,\one) = p$ and $(p,q) \leq (\one,q) = q$. 
Finally, to avoid verbosity, we shall in future state the assumption on $\la \AA_i : i \in I \ra$ in Definition~\ref{defin3}
simply as ``$\la \AA_i : i \in I \ra$ is a system of cBa's with complete embeddings and correct projections".
Let us first check that $\leq$ on $D$ is indeed a p.o.

\begin{obs}
$\leq$ is transitive on $D$.
\end{obs}

\begin{proof}
Assume $(p'',q'') \leq (p',q') \leq (p,q)$ where $(i'',j'')$ witnesses 
$(p'',q'') \in D$ and we have $h_{i'' , i'' \land i'} (p'')  \leq_{i'} p'$ and $h_{i' , i' \land i} (p')  \leq_{i} p$ and similarly for the $j$'s and the $q$'s.
This clearly implies 
\[ h_{i'', i'' \land i} (p'') \leq_{i'' \land i} h_{i'', i'' \land i' \land i} (p'') = h_{i' , i' \land i} ( h_{i'', i'' \land i'} (p''))  \leq_{i' \land i} h_{i', i' \land i} (p') \leq_i p, \]
where the equality holds by correctness, and $h_{i'', i'' \land i} (p'')\leq_i p $ then follows from transitivity of $\leq_i$.
\end{proof}

\begin{exam}
Let $I = \{ 0, 1, 0 \land 1 \}$, the three-element almost-lattice. Then clearly
$\AA_\ell = \lim \amal_{i\in I} \AA_i = \amal_{\AA_{0\land 1}} ( \AA_0 , \AA_1)$,
the amalgamation of $\AA_0$ and $\AA_1$ over $\AA_{0 \land 1} $ studied in
Example~\ref{exam2}. As mentioned earlier, the product construction (Example~\ref{exam1})
is a special case of this for $\AA_{0 \land 1} = \{ \zero, \one \}$.
\hfill $\dashv$
\end{exam}

\begin{exam}
Let $I$ be a distributive lattice. Then $\AA_\ell = \lim \amal_{i \in I} \AA_i  =  \lim \dir_{i \in I}
\AA_i$ is the usual direct limit. If $I$ has a maximum $i_0$, we get $\AA_\ell = \AA_{i_0}$
as a special case. \hfill $\dashv$
\end{exam}

The point of Definition~\ref{defin3} is, of course, that it allows us to define an iteration. For this we first
need to establish that the initial  steps of the iteration completely embed into the limit, and this constitutes the
main result of this section.

\begin{mainlem}[complete embeddability into the amalgamated limit]  \label{amal-embed}
All $\AA_k$, $k \in I$, completely embed into $\AA_\ell = \AA_\amal$.
\end{mainlem}

\begin{proof}
Fix $k \in I$. Let $(p,q) \in D$ as witnessed by $(i,j)$ where $i \in L$ and $j \in R$.
So $p \in \AA_i$, $q \in \AA_j$, and $h_{ i, i\land j} (p) = h_{ j ,i\land j} (q)$.
We need to define $h_{\ell k} (p,q)$.

By part (iv) of Basic Lemma~\ref{basic-lemma}, we may assume without loss of generality that $k \in L$. 
(The case $k\in R$ is analogous.) By part (iii) of the same lemma, then, $i \lor k$ exists and belongs to $L$ as well. Then $p \in \AA_{i \lor k}$ and, by correctness of projections,
$h_{i\lor k , (i \lor k) \land j} (p) = h_{i, i \land j} (p) = h_{j, i \land j} (q) \geq h_{j, (i \lor k) \land j} (q)$.
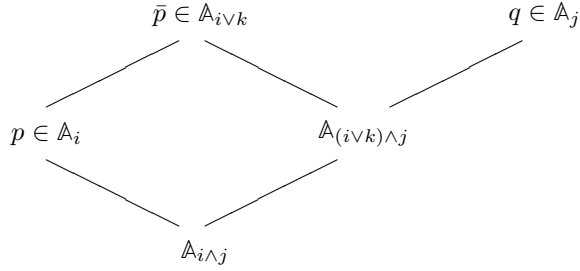
\[ 
\begin{picture}(190,100)(0,0)
\put(60,5){\makebox(0,0){$ \AA_{i \land j} $}}
\put(0,50){\makebox(0,0){$p \in \AA_{i}$}}
\put(60,95){\makebox(0,0){$ \bar p \in \AA_{i \lor k}$ }}
\put(120,50){\makebox(0,0){$ \AA_{(i \lor k) \land j} $}}
\put(190,95){\makebox(0,0){$q \in \AA_j $ }}
\put(0,40){\line(2,-1){50}}
\put(60,15){\line(2,1){50}}
\put(60,85){\line(2,-1){50}}
\put(0,60){\line(2,1){50}}
\put(130,60){\line(2,1){50}}
\end{picture}
\]
Letting $\bar p = p \cdot h_{j, (i \lor k) \land j} (q) \in \AA_{i \lor k}$, we obtain $h_{i \lor k, (i \lor k) \land j} (\bar p) =
h_{j, (i \lor k) \land j} (q)$, and therefore $( i \lor k, j)$ witnesses that $(\bar p, q)$ belongs to $D$.
(In fact, as proved below in Observation~\ref{basicobs} (iv), $(p,q)$ and $(\bar p, q)$ are equivalent conditions in $\AA_\ell$.)
Now let $h_{\ell k} (p,q) = h_{i \lor k, k} (\bar p)$.

To check this works, let $r \in \AA_k$, $r \leq h_{\ell k} (p,q) $.
We need to show that $r = (r, \one)$ and $(p,q)$ are compatible in $\AA_\ell$.
Since $r \leq h_{i \lor k, k} (\bar p)$, $r$ and $\bar p$ are compatible in $\AA_{i \lor k}$ with common extension $p' = \bar p \cdot r$.
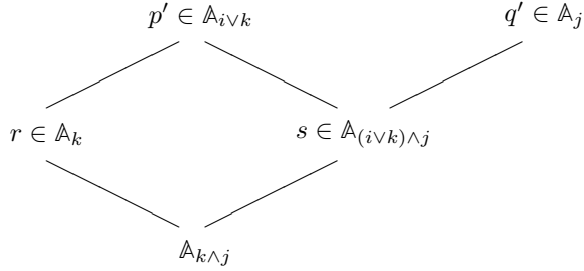
\[ 
\begin{picture}(190,100)(0,0)
\put(60,5){\makebox(0,0){$ \AA_{k \land j} $}}
\put(0,50){\makebox(0,0){$r \in \AA_{k}$}}
\put(60,95){\makebox(0,0){$  p ' \in \AA_{i \lor k}$ }}
\put(120,50){\makebox(0,0){$ s \in \AA_{(i \lor k) \land j} $}}
\put(190,95){\makebox(0,0){$q ' \in \AA_j $ }}
\put(0,40){\line(2,-1){50}}
\put(60,15){\line(2,1){50}}
\put(60,85){\line(2,-1){50}}
\put(0,60){\line(2,1){50}}
\put(130,60){\line(2,1){50}}
\end{picture}
\]
Then $h_{i \lor k, ( i \lor k) \land j} (p') \leq h_{i \lor k, (i \lor k) \land j} (\bar p) = h_{j, (i \lor k) \land j} (q) $ and, letting
$s = h_{i \lor k, (i \lor k) \land j} (p')$ and $q' = q \cdot s$, we see that $h_{j, (i \lor k) \land j} (q') = s$, so that $(i \lor k, j)$
witnesses that $(p',q') \in D$. Since $p' \leq r$, $p' \leq \bar p \leq p$ and $q' \leq q$, $(p',q')$ is a common extension 
of $(p,q)$ and $r = (r,\one)$. 
This completes the proof.
\end{proof}

We saw that correctness was crucial for this proof, and since we will later (in Section~\ref{Cohen})
build an iteration using the amalgamated limit, we had better check that correctness is preserved by the
latter.

\begin{lem}   \label{amal-correct}
Assume $\AA_\ell = \AA_\amal$ is the amalgamated limit of $\la \AA_i : i \in I \ra$.
Then $\la \AA_i : i \in I \cup \{ \ell \} \ra$ still has correct projections, where $\ell \geq i$ for all $i \in I$.
\end{lem}

Note that $I \cup \{ \ell \}$ is a lattice by Observation~\ref{almost-lattice-basic}.

\begin{proof}
We need to check that if $i\in L$ and $j \in R$ have no upper bound in $I$, then $h_{\ell  i} (q) = h_{j , i \land j} (q)$ for $q \in \AA_j$.
That is, it suffices to argue that any given $p \in \AA_i$ with $p \leq h_{j, i \land j} (q)$
is compatible with $q$ in $\AA_\ell$. This, however, is easy.
Indeed,  $h_{i, i\land j} (p) \leq h_{j ,i\land j} (q)$ so that,
letting $q' = q \cdot h_{i, i\land j} (p)$ we have $h_{j ,i\land j} (q') = h_{i, i\land j} (p)$.
This means $(p, q') \in \AA_\ell$ and $(p, q') \leq p, q$ is obvious.
\end{proof}

One may wonder why we require the stringent condition that out of three elements of $I$, two have
a common upper bound. Of course many cBa's can be amalgamated over a {\em  single} common
subalgebra:
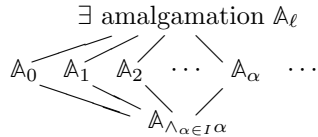
\[
\begin{picture}(110,50)(0,0)
\put(65,3){\makebox(0,0){$\AA_{\land_{\alpha \in I} \alpha}$}}
\put(3,23){\makebox(0,0){$ \AA_0 $}}
\put(23,23){\makebox(0,0){$\AA_1$}}
\put(43,23){\makebox(0,0){$\AA_2$}}
\put(65,23){\makebox(0,0){$\cdots$}}
\put(87,23){\makebox(0,0){$\AA_\alpha$}}
\put(109,23){\makebox(0,0){$\cdots$}}
\put(65,43){\makebox(0,0){$\exists$ amalgamation $\AA_\ell$}}
\put(6,26){\line(3,1){30}}
\put(26,26){\line(2,1){20}}
\put(46,26){\line(1,1){10}}
\put(68,36){\line(1,-1){10}}
\put(9,17){\line(3,-1){34}}
\put(29,17){\line(2,-1){18}}
\put(49,17){\line(1,-1){10}}
\put(68,6){\line(1,1){10}}
\end{picture}
\]
Indeed, $\AA_\ell$ is simply the two-step iteration where we first force with $\AA_{\land_{\alpha \in I} \alpha}$ and then with the
product of the quotients $\AA_\beta / G_{\land_{\alpha \in I} \alpha}$, $\beta \in I$. But what about the following diagram?
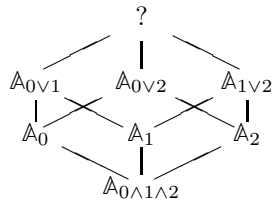
\[
\begin{picture}(90,70)(0,0)
\put(43,3){\makebox(0,0){$ \AA_{0\land 1 \land 2} $}}
\put(3,23){\makebox(0,0){$ \AA_0 $}}
\put(43,23){\makebox(0,0){$ \AA_1 $}}
\put(83,23){\makebox(0,0){$ \AA_2 $}}
\put(3,43){\makebox(0,0){$ \AA_{0 \lor 1} $}}
\put(43,43){\makebox(0,0){$ \AA_{0 \lor 2} $}}
\put(83,43){\makebox(0,0){$ \AA_{1 \lor 2} $}}
\put(43,68){\makebox(0,0){?}}
\put(6,26){\line(2,1){22}}
\put(6,48){\line(2,1){25}}
\put(48,6){\line(2,1){25}}
\put(48,26){\line(2,1){22}}
\put(38,6){\line(-2,1){25}}
\put(36,26){\line(-2,1){22}}
\put(78,26){\line(-2,1){22}}
\put(80,48){\line(-2,1){25}}
\put(3,28){\line(0,1){8}}
\put(83,28){\line(0,1){8}}
\put(43,48){\line(0,1){10}}
\put(43,8){\line(0,1){10}}
\end{picture}
\]
As the next example shows, the amalgamation does not exist in general.

\begin{counterexam}   \label{counterexam3}
Let $\AA_{0 \land 1 \land 2} = \{ \zero, \one \}$, $\AA_0 = \AA_1 = \AA_2 = \ro (\DD)$,
that is, Hechler forcing, $\AA_{0 \lor 1} = \AA_0 \star \dot \AA_1$,
$\AA_{1 \lor 2 } = \AA_1 \star \dot \AA_2$, and $\AA_{0 \lor 2} = \AA_2 \star \dot \AA_0$. By
Example~\ref{exam3} we know all embeddings are complete and the corresponding
projections are correct. Suppose there was an algebra $\AA_\ell$ with
$\AA_{0 \lor 1 } \embed \AA_\ell$, $\AA_{0 \lor 2} \embed \AA_\ell$, and
$\AA_{1 \lor 2} \embed \AA_\ell$. Then $\AA_\ell$ would adjoin three reals
$d_0$, $d_1$ and $d_2$, added by $\AA_0$, $\AA_1$ and $\AA_2$, respectively,
such that $d_1$ eventually dominates $d_0$, $d_2$ eventually dominates $d_1$,
and $d_0$ eventually dominates $d_2$, which is absurd. \hfill $\dashv$
\end{counterexam}

Thus Main Lemma~\ref{amal-embed} is optimal and this accounts for the fact that we do not
know how to build ``three-dimensional" shattered iterations, see Question~\ref{threedimensions}
in Section~\ref{problems}.

Another natural question is whether Lemma~\ref{correct-limit}
can be extended to the context of the amalgamated limit. Again this is
false in general. This somewhat unfortunate state of affairs means that preservation
results are not as easy as for, say, finite support iteration of Suslin ccc forcing.
We shall come back to this issue in Section~\ref{preservation}.

\begin{counterexam}   \label{counterexam4}
Let $\AA_{0 \land 1} = \{ \zero, \one \}$, $\AA_0 = \AA_1 = \EE_{0 \land 1} = \ro (\DD)$,
$\EE_0 = \AA_0 \star \dot \EE_{0 \land 1}$ and $\EE_1 = \EE_{0 \land 1} \star \dot \AA_1$,
and let $d_0$, $d_1$ and $d$ denote the corresponding Hechler generics.
Then all diagrams of the form 
\begin{picture}(60,50)(0,0)
\put(28,3){\makebox(0,0){$\AA_i$}}
\put(5,23){\makebox(0,0){$ \EE_i$}}
\put(48,23){\makebox(0,0){$\AA_j$}}
\put(28,43){\makebox(0,0){$\EE_j$}}
\put(8,26){\line(1,1){13}}
\put(10,18){\line(1,-1){11}}
\put(31,6){\line(1,1){13}}
\put(33,38){\line(1,-1){11}}
\end{picture}
are correct. Put $\AA_{0 \lor 1} = \amal_{\AA_{0 \land 1}} (\AA_0 , \AA_1) = \AA_0 \times \AA_1$
and $\EE_{0 \lor 1} = \amal_{\EE_{0 \land 1}} (\EE_0 , \EE_1)$. Thus $d_0$ and $d_1$ are independent
in the $\AA_{0 \lor 1}$-generic extension while $d_0 <^* d <^* d_1$ in the $\EE_{0 \lor 1}$-generic
extension. A fortiori $\AA_{0 \lor 1} \not{\!\!\embed}  \EE_{0 \lor 1}$. \hfill $\dashv$
\end{counterexam}

\begin{problem}    \label{amal3-problem}
Find sufficient conditions which guarantee the existence of the amalgamation and/or the amalgamated
limit for structures more general than distributive almost-lattices.
\end{problem}

\begin{problem}   \label{amalembed-problem}
Find sufficient conditions which guarantee complete embeddability between amalgamated limits
(\`a la Lemma~\ref{correct-limit}).
\end{problem}

We collect a number of basic facts about the amalgamated limit which we will
use later on, in particular in Section~\ref{coherent}. While quite technical, all these properties should be seen
as straightforward. For $(p,q) \in D$, say that {\em $p \cdot q$ exists in some $\AA_k$} if both $p$ and $q$ belong to
$\AA_k$ so that their product $p \cdot q$ is defined in $\AA_k$. If $(p,q) \in D$ is witnessed by $(i,j)$, this means that
$p \in \AA_{i \land k}$ and $q \in \AA_{j \land k}$\footnote{This notion, as well as the last two items of Observation~\ref{basicobs}
and Observation~\ref{basicobsquotient} (iv), are not really needed later on (in Section~\ref{coherent}); we include them since we consider them helpful for better understanding 
the amalgamated limit.}.

\begin{obs}   \label{basicobs}
\begin{enumerate}
\item If $(p,q) \in D$ is witnessed by $(i,j)$, $i' \leq i, j' \leq j$,  $p \in \AA_{i'}$
   and $q \in \AA_{j'}$, then $(i',j')$ also witnesses $(p,q) \in D$.
\item If $(p,q) \in D$ is witnessed by $(i,j)$, $i' \geq i$ and $j' \geq j$, then clearly $p \in \AA_{i'}$,
   $q \in \AA_{j'}$, and, setting $p' = p \cdot h_{j, i' \land j} (q)$ and $q' = q \cdot h_{i,i\land j'} (p)$,
   $(i',j')$ witnesses $(p ',q') \in D$.
\item If $(p,q), (p',q') \in D$ are both witnessed by $(i,j)$, and $(p,q)$ and $(p',q')$ are equivalent in $\AA_\ell$,
   then $p=p'$ and $q = q'$.
\item If $(p,q) \in D$ is witnessed by $(i,j)$, $i' \geq i$ and $j' \geq j$, then there is
   a unique $(p',q') \in D$ with witness $(i',j')$ such that $(p,q)$ and $(p',q')$ are
   equivalent in $\AA_\ell$ and $(p',q') \leq (p,q)$. 
   Furthermore $p' = p \cdot h_{j, i' \land j} (q)$ and $q' = q \cdot h_{i,i\land j'} (p)$.
\item If $(p,q)$ and $(p',q')$ belong to  $D$, as witnessed by $(i,j)$ and $(i',j')$ respectively, and
   are equivalent in $\AA_\ell$, 
   then there is $(p'',q'') \in D$ with witness $(i \lor i', j \lor j')$,
   $(p'',q'') \leq (p,q) , (p',q')$, equivalent to both. Furthermore $p'' = p \cdot
   h_{j,j \land (i \lor i')} (q) = p' \cdot h_{j', j' \land (i \lor i')} (q')$ and $q'' = q \cdot h_{i,
   i \land (j \lor j')} (p) = q' \cdot h_{i' , i' \land (j \lor j')} (p')$.
\item If $(p_m, q_m) \in D$ as witnessed by $(i_m, j_m)$, $m < n$, then there are $(i',j')$ and
   $(p_m ' , q_m ')$, $m < n$, such that $(p_m ' , q_m ') \in D$ is witnessed by $(i',j')$ and
   $(p_m ', q_m ')$ is equivalent to $(p_m , q_m)$ in $\AA_\ell$.
\item If $(p,q) \in D$ is witnessed by $(i,j)$ and $p \cdot q$ exists in some $\AA_k$, $k\in I$,
   then $(p,q)$ and $p \cdot q $ are equivalent in $\AA_\ell$.
\item If $(p,q)$, $(p',q') \in D$ are equivalent in $\AA_\ell$ and $p \cdot q$ exists,
   then so does $p' \cdot q'$ and $p \cdot q = p' \cdot q'$ is equivalent to both.
\end{enumerate}
\end{obs}

\begin{proof}
(i) Straightforward.

(ii) By correctness, $h_{i',i' \land j} (p) = h_{i,i\land j} (p) = h_{j, i \land j} (q) \geq h_{j , i' \land j} (q)$
so that $h_{i', i' \land j} (p') = h_{j , i' \land j} (q)$. Again by correctness, $h_{j' , i' \land j'} (q)
= h_{j, i' \land j} (q)$ so that $h_{j', i' \land j'} (q) \geq h_{i' , i' \land j'} (p')$.
As evidently $h_{i, i \land j'} (p) \geq h_{i' , i' \land j'} (p')$,
$h_{j' , i' \land j'} (q') \geq h_{i', i' \land j'} (p')$ follows.
By symmetry, they must be equal and $(i',j')$ witnesses $(p',q') \in D$.

(iii) Assume not. Then, without loss of generality, there is $p'' \in \AA_i$, $p'' \leq p$,
$p''$ incompatible with $p'$. Set $q'' = q \cdot h_{i,i\land j} (p'')$. Thus $h_{j, i \land j} (q'') =
h_{i , i \land j} (p'')$, and $(p'',q'') \in D$ is witnessed by $(i,j)$.
Since $(p'',q'') \leq (p,q)$, yet $(p'',q'')$ and $(p',q')$ are incompatible,
we reach a contradiction.

(iv) By (ii), we know $(i',j')$ witnesses $(p',q') \in D$. $(p',q') \leq (p,q)$ is immediate by definition.
Uniqueness follows from (iii). So it suffices to prove equivalence. To this end let
$(p'',q'') \in D$ as witnessed by $(i'',j'')$ with $(p'',q'') \leq (p,q)$. We need to show
$(p'',q'')$ and $(p',q')$ are compatible.

Replacing $(i'',j'')$ by $(i'' \lor i', j'' \lor j')$ and $(p'',q'')$ by $(p'' \cdot h_{j'', (i'' \lor i') \land j''} (q''), q'' \cdot h_{i'' , i'' \land ( j'' \lor j')} (p'') )$ using (ii) if necessary,
we see that we may assume $i'' \geq i'$ and $j'' \geq j'$. By assumption $p'' \leq p$. Furthermore, by correctness, $p'' \leq h_{i'' , i' \land j} (p'') = h_{j'' , i' \land j} (q'') 
\leq h_{j, i' \land j} (q)$. Thus $p'' \leq p \cdot h_{j, i' \land j} (q) = p'$. Similarly we see that $q'' \leq q'$, so that $(p'',q'') \leq (p',q')$ follows as required.

(v) This follows from (iv).

(vi) Simply let $i' = \bigvee_{m<n} i_m$, $j' = \bigvee_{m<n} j_m$, and use (iv).

(vii) Without loss of generality assume that $k \in L$. Then $i \lor k$ exists and also belongs to $L$ while $j \land k \in \ker (I)$. The condition $p \cdot q \in \AA_k$
is identified with $(p \cdot q , \one) \in D$, and we need to show it is equivalent with $(p,q)$.

First assume $(p',q')$ belongs to $D$ as witnessed by $(i',j')$ with $(p',q') \leq (p \cdot q, \one)$. As in the proof of (iv), using (ii), we may assume without loss
of generality that $i' \geq i \lor k$ and $j' \geq j$. Since $p' \leq p \cdot q$, $p' \leq p$ is straightforward. By correctness we see that $q' \leq
h_{j', j \land k} (q') = h_{i', j \land k} (p') \leq h_{k, j \land k} ( p \cdot q) \leq h_{k, j \land k} (q) = q$ (where the last equality holds because $q \in \AA_{j\land k}$).
Therefore $(p',q') \leq (p,q)$.

On the other hand, assume $(p',q') \leq (p,q)$. Arguing as in the previous paragraph and using correctness again, we see that $p' \leq p$ and $p' \leq h_{i', j \land k} (p') = 
h_{j', j \land k} (q') \leq q$, so that $p' \leq p \cdot q$ and $(p' , q') \leq (p \cdot q , \one)$ follow, as required.

(viii) Assume again that $p \cdot q$ exists in $\AA_k$ where $k \in L$. Let $(i',j')$ witness $(p',q') \in D$. By (v) (see also (ii)) we may assume $i' \geq i \lor k$ and
$j' \geq j$. We then know that $p' = p \cdot h_{j\land k, i' \land j \land k} (q) \in \AA_k$ and $q' = q \cdot h_{i\land k, i \land j' \land k} (p) \in \AA_k$ so that, in
$\AA_k$, $p' \cdot q' = p \cdot h_{j\land k, i' \land j \land k} (q) \cdot q \cdot h_{i\land k, i \land j' \land k} (p) = p \cdot q$, as required. Equivalence follows from (vii).
\end{proof}

In Section~\ref{coherent} we will need a description of the amalgamated limit as
a two-step iteration of the form $\AA_\ell \cong \AA_{i_0} \star \dot\AA_{i_0 \ell}$ where
$i_0 \in I$ and $\dot \AA_{i_0 \ell}$ is the name for the remainder forcing leading from the
$i_0$-generic extension to the $\ell$-generic extension. 

Assume without loss of generality that $i_0 \in L$ (the case $i_0 \in R$ is analogous).
Let $G_{i_0}$ be $\AA_{i_0}$-generic over $V$. Work in $V [G_{i_0}]$ for the moment.
Let $D_{i_0}$ be the collection of (ordered pairs) $(p,q) \in \AA \times \AA$ such that
for some $i \in L$ with $i \geq i_0$ and $j\in R$, we have $p \in \AA_i$, $q \in \AA_j$, $h_{i, i \land j} (p) =
h_{j , i \land j} (q)$, and $h_{ii_0} (p) \in G_{i_0}$. The order is given by 
$(p',q') \leq_{i_0 \ell} (p,q)$ if for some $r \in G_{i_0}$ with $r \leq_{i_0} h_{i'  i_0} (p')$ we have
$h_{i' , i'\land i } (p') \cdot r \leq_{i} p$ and $h_{j', j' \land j} (q')  \cdot h_{i_0, i_0 \land j} (r) \leq_{j} q$ 
where $(i',j')$ witnesses $(p',q') \in D_{i_0}$. $\AA_{i_0 \ell}$ is the completion of $D_{i_0}$,
$\AA_{i_0 \ell} = \ro (D_{i_0})$. Note that $D_{i_0 } \sub D$ in $V[G_{i_0}]$. On the other hand,
$(p,q) \in D \sem D_{i_0}$ can be identified with $\zero$ in $\AA_{i_0 \ell}$.

\begin{lem}  \label{amal-remainder}
$\AA_\ell \cong \AA_{i_0} \star \dot\AA_{i_0 \ell}$.
\end{lem}

\begin{proof}
We define a dense embedding $e$ from $D$ into the two-step iteration $\AA_{i_0} \star
\dot D_{i_0}$. Let $(p,q) \in D$ be witnessed by $(i,j)$. Then $i \lor i_0$ exists and belongs to $L$.
Increasing $i$, if necessary (see Observation~\ref{basicobs} (iv)), we may assume $i \lor i_0 = i$,
i.e., $i \geq i_0$. Let $p_0 = h_{i i_0} (p)$. Then $p_0 \forces (p,q) \in \dot D_{i_0}$. Thus
let $e$ map $(p,q)$ to $(p_0 , (p,q))$. 

It is easy (though tedious, because there are several cases) to check that $e$ is well-defined (i.e., $e$ maps
equivalent conditions in $D$ to equivalent conditions in $\AA_{i_0} \star \dot D_{i_0}$) and that it preserves
order and incompatibility. 

To see that $e$ is dense, assume $p_0 \forces ``(p,q) \in \dot D_{i_0}$, as witnessed
by $i$ and $j$". Thus $p_0 \leq h_{ii_0} (p)$. Let $p' = p_0 \cdot p$. Then $h_{ii_0} (p')
= p_0$. Next let $q' = h_{i,i\land j} (p')\cdot q$. Then $h_{j,i\land j} (q') = h_{i,i\land j} (p')$.
Thus $e$ maps $(p',q')$ to $(p_0, (p',q'))$ and clearly $(p_0, (p',q')) \leq (p_0, (p,q))$.
\end{proof}

Note that this lemma subsumes Main Lemma~\ref{amal-embed}. We may think of 
$\dot \AA_{i_0 \ell}$ as the amalgamated limit of the remainders. In particular,
for $i \geq i_0$ from $L$, $\dot \AA_{i_0 i}$ is forced to canonically embed into $\dot \AA_{i_0 \ell}$:

\begin{cor}
$\forces_{i_0} \dot \AA_{i_0 i} \embed \dot \AA_{i_0 \ell}$ for $i_0 \leq i$.
\end{cor}

\begin{proof}
Applying the previous lemma twice we see that $\AA_{i_0} \star \dot \AA_{i_0 \ell}
\cong \AA_\ell \cong \AA_i \star \dot \AA_{i\ell} \cong \AA_{i_0} \star \dot\AA_{i_0 i} \star \dot\AA_{i\ell}$.
Hence $\AA_{i_0}$ forces that $\dot \AA_{i_0 \ell} \cong \dot \AA_{i_0 i} \star \dot \AA_{i \ell}$.
In particular, it forces $\dot \AA_{i_0 i} \embed \dot \AA_{i_0 \ell}$.
\end{proof}

We conclude this section with some basic properties of the $\AA_{i_0 \ell}$, which we will need in 
Section~\ref{coherent}. For $(p,q) \in D_{i_0}$, say that {\em $p \cdot q$ exists in $D_{i_0}$} if there are $k \in I$ and $r \in G_{i_0}$
such that $r \cdot p$ and $r \cdot q$ both belong to $\AA_k$ and therefore we may form the product $r \cdot p \cdot q \in \AA_k$.
Note that if $(p,q) \in D_{i_0}$ and $p \cdot q$ exists in $D$, as defined earlier (before Observation~\ref{basicobs}), then 
$p \cdot q$ exists in $D_{i_0}$. The converse, however, need not hold. (One may have, e.g., that $i_0 \leq i \in L, i_0 \leq j \in R$, $i \lor j$ does not
exist, $(p,q)$ belongs to $D$ as witnessed by $(i,j)$, $p \cdot q$ does not exist either, $h_{ii_0} (p) = h_{ji_0} (q) \in G_{i_0}$, and for some
$r \in G_{i_0}$ strictly stronger than $h_{ii_0} (p)$, $r = r \cdot p = r \cdot q$ belongs to $\AA_{i_0}$.)

\begin{obs} \label{basicobsquotient}
Letting $i_0 \in L$ and $G_{i_0}$ an $\AA_{i_0}$-generic filter over $V$, the following hold in $V[G_{i_0}]$:
\begin{enumerate}
\item Suppose $(p,q), (p',q') \in D_{i_0}$, both witnessed by $(i,j)$ with $i \geq i_0$, are equivalent in $\AA_{i_0\ell}$. Then there is
$r \in G_{i_0}$ such that $p \cdot r = p' \cdot r$ and $q \cdot h_{i_0, i_0 \land j} (r)  = q' \cdot h_{i_0, i_0 \land j} (r)  $.
\item Suppose $i' \geq i \geq i_0$ and $j' \geq j$, $(p,q) \in D_{i_0}$ and $(p',q') \in D_{i_0}$ are witnessed by $(i,j)$
and $(i',j')$ respectively, and the two conditions are equivalent in $\AA_{i_0\ell}$. Then for some $r \in G_{i_0}$, we have
$p'  \cdot r = p \cdot h_{j,i' \land j} (q) \cdot r$ and $q' \cdot h_{i_0, i_0 \land j'} (r) = q \cdot h_{i, i \land j'} (p) \cdot h_{i_0, i_0 \land j'} (r)$.
\item Suppose $(p,q), (p', q') \in D_{i_0}$, as witnessed by $(i,j)$ and $(i',j')$ with $i,i' \geq i_0$, respectively, are equivalent in $\AA_{i_0\ell}$. 
Then there are $r \in G_{i_0}$ and $(p'',q'') \in D_{i_0}$ witnessed by $(i \lor i', j \lor j')$ equivalent to both in $\AA_{i_0\ell}$ 
such that $p'' \cdot r = p \cdot h_{j, j \land (i \lor i')} (q) \cdot r = p' \cdot h_{j', j '\land (i \lor i')} (q') \cdot r $
and $q'' \cdot h_{i_0, i_0 \land (j \lor j')} (r) = q \cdot h_{i, i \land (j \lor j')} (p) \cdot h_{i_0, i_0 \land (j \lor j')} (r) = q' \cdot h_{i', i '\land (j \lor j')} (p') \cdot h_{i_0, i_0 \land (j \lor j')} (r)$.
\item Suppose $(p,q), (p',q') \in D_{i_0}$ are equivalent in $\AA_{i_0\ell}$, the latter witnessed by $(i',j')$ with $i' \geq i_0$, and
$p \cdot q$ exists in $D_{i_0}$ as witnessed by $k \in I$ and $r \in G_{i_0}$. Then for some $r' \in G_{i_0}$ with $r' \leq r$, we have
   $r' \cdot p \cdot q = r' \cdot p' \cdot q'$.
\end{enumerate}
\end{obs}

Note that the conclusions in the various parts of this observation are really characterizations of equivalence in $\AA_{i_0 \ell}$. For example, in (i), it
is obvious that if $p \cdot r = p' \cdot r$ and $q \cdot h_{i_0, i_0 \land j} (r)  = q' \cdot h_{i_0, i_0 \land j} (r)  $ for some $r \in G_{i_0}$ then
$(p,q)$ and $(p',q')$ are equivalent as conditions in $\AA_{i_0\ell}$, and similarly for the other items.

\begin{proof}
(i) First consider $p$ and $p'$. In the ground model $V$, let $E = \{ r \leq_{i_0}  h_{i i_0} (p)  : r \cdot p \leq_i p'$ or there is $p'' \leq_i p$ such that
$h_{ii_0} (p'') = r$ and $p'' \bot p' \}$ and note that this set is predense in $\AA_{i_0}$ below $h_{i i_0} (p)$. Therefore $r \in E$ for some $r \in G_{i_0}$.
If $h_{i i_0} (p'') = r$ and $p'' \bot p'$, then, letting $q'' = q \cdot h_{i, i \land j} (p'')$, we see $(p'' , q'') \in D_{i_0}$ is below $(p,q)$ and
incompatible with $(p',q')$ in $\AA_{i_0\ell}$, a contradiction. Thus $r \cdot p \leq p'$. By symmetry, $r \cdot p' \leq p$ for some $r \in
G_{i_0}$, and since $G_{i_0}$ is a filter we obtain $r \in G_{i_0}$ such that $r \cdot p = r \cdot p'$, as required.

Next consider $q$ and $q'$. Again in $V$, let $E = \{ r \leq_{i_0}  h_{i i_0} (p) : h_{i_0, i_0 \land j} (r) \cdot q \leq_j q'$ or there is $q'' \leq_j q$ such that
$h_{j, i_0 \land j} (q'') = h_{i_0, i_0 \land j} (r)$ and $q '' \bot q' \}$. Again we find $r \in G_{i_0} \cap E$, and if the second option holds, 
letting $p'' = p \cdot h_{j, i \land j} (q'')$, $(p'',q'') \in D_{i_0}$ is below $(p,q)$ and incompatible with $(p',q')$, a contradiction. So $h_{i_0, i_0 \land j} (r) 
\cdot q \leq q'$. The converse
is again by symmetry.

(ii) This is very much like the proof of (i). There are four predense sets to consider. We do the first one and leave the rest to the reader.
In the ground model $V$, let $E = \{ r \leq_{i_0}  h_{i i_0} (p) : r \cdot p' \leq_{i'} p \cdot h_{j, i' \land j} (q)$ or there is $p'' \leq_{i'} p'$ such that
$h_{i' i_0} (p'') = r$ and $p'' \bot p \cdot h_{j, i'\land j} (q) \}$. Again this set is predense in $\AA_{i_0}$ below $h_{i i_0} (p)$, and we may take $r \in G_{i_0} \cap E$.
If the second option holds, setting $q'' = q ' \cdot h_{i', i' \land j'} (p'')$, we see that $(p'',q'')$ belongs to $D_{i_0} $, is below $(p',q')$ and incompatible 
with $(p,q)$ in $\AA_{i_0\ell}$, a contradiction. Thus $r \cdot p' \leq_{i'} r \cdot p \cdot h_{j, i' \land j} (q)$. The rest is similar.

(iii) Letting $p''  = p \cdot h_{j, j \land (i \lor i')} (q)$ and $q''  = q \cdot h_{i, i \land (j \lor j')} (p)$, we know by Observation~\ref{basicobs} (iv) (see also (v)) that
$(p'',q'')$ is equivalent to $(p,q)$ in $\AA_\ell$ and, a fortiori, in $\AA_{i_0\ell}$. The rest follows by (ii).

(iv) This is again done with predense sets using part (viii) of Observation~\ref{basicobs}. We leave the details to the reader.
\end{proof}

While special cases of the amalgamated limit have been used implicitly in other constructions, see e.g. the discussion in~\cite[Section 1]{survey}, in particular
Lemma 13, it has never been considered in the generality in which we discuss it here. For a survey on applications of the amalgamated limit and on its topological
interpretation in the language of compact Hausdorff spaces see~\cite{amallimit}.



\section{Coherent systems of measures}   \label{coherent}

The principal technical problem in the proof of the main result is to show the forcing is
ccc and, thus, preserves cardinals. To appreciate this difficulty, recall that while
the countable chain condition is preserved under two-step iterations as well as in 
limit stages of finite support iterations,
it is in general not preserved under products, and strengthenings of the ccc which are
productive (like Knaster's condition) need not be preserved under amalgamation.

\begin{counterexam}   \label{counterexam5}
Let $\TT_{\dot c}$ be the (canonical) $\CC$-name for the Suslin tree added by $\CC$~\cite{Sh84,To87}.
It is well-known $\PP = \CC \star \TT_{\dot c}$ satisfies Knaster's condition (see the argument in~\cite[pp. 292-293]{To87}).
Furthermore, if MA + $\neg$ CH holds, then $\PP$ is even $\sigma$-centered.
However, the amalgamation of two copies of $\PP$ over $\CC$, namely $\CC
\star (\TT_{\dot c} \times \TT_{\dot c})$, is not ccc.   \hfill $\dashv$
\end{counterexam}

This example shows we need to guarantee the quotient algebras arising in our iteration are
nice enough. Knaster's condition may be a plausible choice, but we shall opt for the
stronger finitely additive measurability.
Recall that given a Boolean algebra $\AA$, $\mu : \AA \to [0,1]$ is a {\em finitely additive
measure} if
\begin{itemize}
\item $\mu (\zero) = 0, \mu (\one) = 1$,
\item $\mu (a + b) = \mu (a) + \mu (b)$ for incompatible $a,b \in \AA$ (i.e. $a \cdot b = \zero$).
\end{itemize}
$\mu$ is {\em strictly positive} if also
\begin{itemize}
\item $\mu (a) > 0$ whenever $a \in \AA \setminus \{ \zero \}$.
\end{itemize}
All measures we shall be dealing with will be finitely additive and strictly positive
so that we shall mostly refrain from stating this explicitly. Note that algebras carrying such
a measure are ccc. Such algebras have been studied in detail by Kamburelis~\cite{Ka89}
and we first review some of his work.

First note that any random algebra $\BB_\kappa$ carries a(n even countably additive) strictly positive measure.
Next, every $\sigma$-centered cBa $\AA$ carries such a measure, for there are ultrafilters $\U_n$, $ n \in \omega$,
on $\AA$ such that $\AA \sem \{ \zero \} = \bigcup_n \U_n$, and we may define $\mu (a) = \sum \{  2^{- (n+1)} : a \in \U_n \}$.
Strict positivity is obvious and finite additivity follows from the fact that the $\U_n$ are ultrafilters. Next we have:

\begin{lem}[Kamburelis {\cite[Proposition 2.1]{Ka89}}]   \label{secondextension}
Assume the Boolean algebra $\AA$ carries a finitely additive strictly positive measure $\mu$. Then $\mu$
can be extended to a finitely additive strictly positive measure on $\ro (\AA)$. 
\end{lem}

Kamburelis' proof of this result is rather indirect: canonically embed $\AA$ into the $\sigma$-algebra $\tilde \AA$ generated by the clopen
subsets of the Stone space $X$ of $\AA$; note that the measure on $\AA$ induces a $\sigma$-additive measure on the clopen subsets of $X$, which can
be extended to a $\sigma$-additive measure on $\tilde \AA$. By Sikorski's extension theorem~\cite[Theorem 5.9]{Ko89},
the embedding of $\AA$ into $\tilde \AA$ can be extended to an (in general not complete) embedding of $\ro (\AA)$ into $\tilde \AA$,
and this embedding induces a measure on $\ro (\AA)$. Since the extension is (in general) not unique, neither is the measure.
We will later use a different and more direct argument (proof of Main Lemma~\ref{amal-coherent}), which we include for the sake of completeness.
The point is that we need that various extended measures cohere, and for this a non-arbitrary uniform way of measure extension is necessary. 

\begin{proof}
By a standard application of Zorn's Lemma, it suffices to prove that given $b \in \ro (\AA) \sem \AA$, there is a finitely additive strictly
positive measure $\mu '$ on $\AA' = \la \AA \cup \{ b \} \ra$, the Boolean algebra generated by $\AA$ and $b$, extending $\mu$.

Let $\varepsilon = \pm 1$, and $(+1) b = b, (-1) b = -b$. For $a \in \AA$ define
\[ \tilde \mu ( a \cdot \varepsilon b ) = \sup \{ \mu ( a') : a' \in \AA, a' \leq a \cdot \varepsilon b \} . \]
Note that since $\AA$ is dense in $\ro (\AA)$, $\tilde \mu ( a \cdot \varepsilon b ) = 0$ iff $ a \cdot \varepsilon b = \zero$. Next let
\[ \mu ' (a \cdot \varepsilon b) = {1 \over 2} \Big( \mu (a) +  \tilde \mu ( a \cdot \varepsilon b ) -  \tilde \mu ( a \cdot (- \varepsilon b) )  \Big). \]
Again, $\mu ' ( a \cdot \varepsilon b ) = 0$ iff $ a \cdot \varepsilon b = \zero$ is immediate. Furthermore, if
$a \cdot \varepsilon b \in \AA$, then $\mu ' (a \cdot \varepsilon b) = \tilde \mu (a \cdot \varepsilon b) = \mu (a \cdot \varepsilon b)$.
An arbitrary element of $\AA'$ is of the form $a_0 + a_1 b + a_2 (-b)$ where the $a_i \in \AA$ are pairwise incompatible.
Define 
\[ \mu ' ( a_0 + a_1 b + a_2 (-b) ) = \mu (a_0) + \mu ' (a_1 b) + \mu ' (a_2 (-b)). \]
If $a_0 + a_1 b + a_2 (-b) \in \AA$, then $a_1 b$ and $a_2 (-b)$ also belong to $\AA$, and $\mu ' ( a_0 + a_1 b + a_2 (-b) )  = \mu (a_0) + \mu (a_1 b) + \mu (a_2 (-b))
= \mu  ( a_0 + a_1 b + a_2 (-b) )$ follows, so that $\mu'$ indeed extends $\mu$. We already saw that $\mu'$ is strictly positive,
and since we have 
\[ \mu' (a b) + \mu ' (a (-b)) = \mu (a) \]
as well as 
\[ \tilde \mu ( a \cdot \varepsilon b ) + \tilde \mu ( \bar a \cdot \varepsilon b ) = \tilde \mu ( (a + \bar a)\cdot \varepsilon b ) \;\;\;\mbox{ and }\;\;\;
 \mu '( a \cdot \varepsilon b ) +  \mu ' ( \bar a \cdot \varepsilon b ) = \mu ' ( (a + \bar a)\cdot \varepsilon b ) \]
for incompatible $a, \bar a \in \AA$, a simple computation shows finite additivity.
\end{proof}

We will also need some basic facts concerning the possibility of extending
measures on dense subsets to measures on the whole algebra.

\begin{prelem}   \label{firstextension}
Assume the (not necessarily complete) Boolean algebra $\AA$ is generated by the dense set $D$
and $\mu : D \to (0,1]$ is given such that
\begin{enumerate}
\item $\one \in D$ and $\mu (\one) = 1$, 
\item $D$ is closed under products (i.e., if $a,b \in D$ and $a \cdot b > \zero$ in $\AA$, then $a \cdot b \in D$),
\item if $a, b \in D$, then there are pairwise incompatible $c_j \in D, j < n,$ such that $a - b
   = \sum_{j<n} c_j$ (where we stipulate $\sum_{j<n} c_j = \zero$ if $n=0$),
\item $\mu$ is finitely additive, that is, given $a_i \in D, i<n,$ pairwise incompatible
   such that $a = \sum_{i<n} a_i \in D$, we have $\mu (a) = \sum_{i<n} \mu (a_i)$.
\end{enumerate}
Then there is a (unique) finitely additive strictly positive measure $\bar\mu$ on $\AA$ extending $\mu$.
\end{prelem}

\begin{proof}
Note first that (ii) and (iii) imply that $\AA$ is the closure of $D \cup \{ \zero \}$ under disjoint
(= incompatible) sums. For indeed, let $C$ denote this closure. Since
$(\sum_{i<m} a_i) (\sum_{j<n} b_j) = \sum_{i,j} a_i b_j$ for $a_i, b_j \in D$ and the last sum
is disjoint if the first two are, $C$ is closed under products by (ii).
Next, given pairwise incompatible $a_i, i<m,$ and pairwise incompatible $b_{i'}, i' < m',$
from $D$, we have
\begin{eqnarray}
\sum_{i<m} a_i - \sum_{i' < m'} b_{i'} & = & \sum_{i<m} \left( a_i - \sum_{i' < m'} b_{i'} \right) = \sum_{i<m} \prod_{i' < m'} 
   \left( a_i - b_{i'} \right)  \nonumber \\
& = & \sum_{i<m} \prod_{i' < m'} \sum_{j < n_{ii'}} c_j^{ii'} = \sum_{i<m} \sum_f \prod_{i' < m'} c^{ii'}_{
   f(i')}
\end{eqnarray}
where the $c_j^{ii'} \in D, j < n_{ii'},$ are pairwise incompatible as stipulated by (iii),
and the sum $\sum_f$ runs over all functions $f$ with $\dom (f) = m'$ and $f(i') < n_{ii'}$
(for fixed $i$). By (ii), $\prod_{i' < m'} c_{f(i')}^{ii'} \in D$, and the double sum
in (1) is pairwise disjoint because both the $c_j^{ii'}$ $(j< n_{ii'})$ and the
$a_i$ $(i<m)$ are. Thus $\sum_{i<m} a_i - \sum_{i' < m'} b_{i'} \in C$, and $C$ is closed under complements.
Finally, closure under products, complements and disjoint sums implies closure
under arbitrary sums as well. Hence $C$ is a Boolean algebra, and $C = \AA$ follows.

Let $\bar\mu (\zero) =0$. For $a \in \AA \sem \{ \zero \}$, let $\bar\mu (a) = \sup \{
\sum_{i<n} \mu (a_i) :  \{ a_i : i<n\} \sub D$ is an antichain below $a\}$. For $a\in D$,
$\mu (a) \leq \bar\mu (a)$ is immediate. To see the other inequality, let
$\{ a_i : i<n \} \sub D$ be an antichain below $a$. By (iii) and (ii), this antichain is complemented in $a$
(see (1) above), and
by (iv), $\sum_{i<n} \mu (a_i) \leq \mu (a)$ follows.
Since $D$ is dense, $\bar\mu$ clearly is strictly positive.
To see finite additivity, take incompatible $a,b \in \AA$.
$\bar \mu (a+b) \geq \bar\mu (a) + \bar\mu (b)$ is immediate from the definition.
For the converse note that if $\{ a_i   \in D : i < n \}$ is a finite antichain
below $a + b$, then each $a_i \cdot a$ is a finite disjoint sum of $a_{ij}^a \in D, j < n_i^a$.
Similarly for $a_i \cdot b$. Thus $\mu (a_i) = \sum_{j < n_i^a} \mu (a_{ij}^a )
+ \sum_{j < n_i^b} \mu (a_{ij}^b)$ by (iv). This easily entails $\bar\mu (a + b) \leq
\bar\mu (a) + \bar\mu (b)$.

To see uniqueness note that any extension $\mu '$ must satisfy $\mu ' (a) \geq \bar \mu (a)$ for
$a \in \AA$. Thus $\mu' (a) = \bar \mu (a)$ since $\AA$ is closed under complements.
\end{proof}

\begin{prelem}   \label{properties-of-D}
Let $D$ be the dense set in the definition of the amalgamated limit.
\begin{enumerate}
\item $D$ is closed under products, i.e., if $(p,q),  (p',q') \in D$ are compatible,
   then there is $(p'',q'') \in D$ which is equivalent to $(p,q) \cdot (p',q') \in \AA_\ell$.
\item Given $(p,q), (p',q') \in D$ there are (at most two) $(p_0,q_0), (p_1,q_1) \in D$
   such that $(p,q) - (p',q') $ and $(p_0,q_0) + (p_1,q_1)$ are equivalent in $\AA_\ell$
   and $(p_0,q_0)$ and $(p_1,q_1)$ are incompatible.
\item Let $i_0 \in L$ (or $R$), let $G_{i_0}$ be $\AA_{i_0}$-generic over $V$, and let $D_{i_0}$ be the dense subset of the remainder forcing
   $\AA_{i_0\ell}$ (see before Lemma~\ref{amal-remainder}) in $V [G_{i_0}]$. Then, in $V [G_{i_0}]$, $D_{i_0}$ also satisfies (i) and (ii).
\end{enumerate}
\end{prelem}

\begin{proof}
We present the construction of the required conditions for (i) and (ii), but leave the details to the
reader since the arguments are similar to those in Observation~\ref{basicobs}. Let $(p,q), (p',q') \in D$.
By Observation~\ref{basicobs} (vi) we may assume that this is witnessed by the same pair $(i,j)$.

(i) Assume  additionally that $(p,q)$ and $(p',q')$ are compatible. Letting $p'' = p\cdot p'\cdot h_{j, i\land j }
(q \cdot q') \in \AA_{i }$, $q'' = q \cdot q' \cdot h_{i, i \land j} (p \cdot p') \in \AA_j$, we see that $(p'',q'')$ belongs to $D$ as witnessed by $(i,j)$
because $h_{i,i\land j} (p'') = h_{i, i \land j} (p \cdot p') \cdot h_{j, i \land j} (q \cdot q') = h_{j, i \land j} (q'')$, 
and check that $(p'',q'')$ is equivalent to $(p,q) \cdot (p',q')$ in $\AA_\ell$. (In particular,  note that $p'' \neq \zero$ and $q'' \neq \zero$
by compatibility of $(p,q)$ and $(p',q')$.)

(ii)  Let $p_0 = p - p' \in \AA_{i}$, $q_0 = q \cdot h_{i, i \land j} (p_0) \in \AA_{j}$, $q_1 = q \cdot h_{i, i \land j} (p
\cdot p') - q' \in \AA_{j }$, and $p_1 = p \cdot p' \cdot h_{j , i \land j} (q_1) \in \AA_{i}$. Now check $(p_0,q_0), (p_1,q_1)$ either belong to
$D$ as witnessed by $(i , j )$ or are $\zero$, that $(p,q) - (p',q')$ is equivalent
to $(p_0,q_0) + (p_1,q_1)$ in $\AA_\ell$, and that $(p_0,q_0)$ and $(p_1,q_1)$ are incompatible.

(iii) To see (i) for $D_{i_0}$, note that if $(p,q),  (p',q') \in D_{i_0}$ are compatible in $D_{i_0}$, then they are also compatible in $D$ in $V$ (recall
$D_{i_0} \sub D$). So by (i) there is $(p'',q'') \in D$ equivalent to $(p,q) \cdot (p',q') \in \AA_\ell$. Clearly $(p'',q'')$ is still equivalent to $(p,q) \cdot (p',q') \in \AA_{i_0 \ell}$
in $V[G_{i_0}]$ and -- since this product is non-zero -- we must have $(p'',q'') \in D_{i_0}$. (ii) is proved similarly.
\end{proof}

Going back to Kamburelis' work, we recall that the class of cBa's with a measure is closed under two-step iterations, under products,
under limit stages of finite support iterations and (trivially) under complete subalgebras (see~\cite[2.4, 2.5, 2.6]{Ka89}).
Finally, it turns out that this class is generated by the $\BB_\kappa$ and the $\sigma$-centered cBa's in a rather strong
sense: namely, every cBa with a measure completely embeds into a two-step iteration of some $\BB_\kappa$ with a
$\sigma$-centered cBa~\cite[Proposition 3.7]{Ka89}. 

For later use, let us review the construction of the measure on the two-step iteration in some detail. Again, our presentation 
is more direct than Kamburelis' original one. First we introduce an integral on
a cBa $\AA$ carrying a finitely additive measure $\mu$. Namely, let $a \in \AA$ and let
$\dot x$ be an $\AA$-name for a real such that $a \forces \dot x \in [0,1]$. Then set
\begin{eqnarray*} 
\int_a \dot x \; d\mu & = & \lim_{m \to \infty} \sum_{k = 0}^{2^m -1} \left[ \mu \left(a \cdot \BOOL
{k + 1 \over 2^m} \geq \dot x > {k \over 2^m} \BOOR \right) \cdot {k \over 2^m} \right] \\   & = &
\lim_{m \to \infty} \sum_{k = 0}^{2^m -1} \left[ \mu \left(a \cdot \BOOL
{k + 1 \over 2^m} \geq \dot x > {k \over 2^m} \BOOR \right) \cdot {k + 1 \over 2^m}  \right] .
\end{eqnarray*}
It is clear that these limits exist and are equal, and, in fact, the first limit is a supremum and the second, an infimum.
Now assume $\AA = \AA_0 \star \dot\AA_1$ are cBa's, $\AA_0$ carries a finitely
additive measure $\mu_0$ and $\forces_{\AA_0} ``\dot \AA_1$ carries a finitely additive measure $\dot \mu_1"$.
Then we canonically get a finitely additive measure $\mu$ on $\AA$ by stipulating
\[ \mu (a) = \mu (a_0 , \dot a_1) = \int_{a_0} \dot\mu_1 (\dot a_1) d \mu_0 \]
for $a = (a_0, \dot a_1) \in \AA$. See~\cite[Section 2]{Ka89} for details.

\begin{defin} \label{defin3a}
(1) A distributive almost-lattice $I$ is {\em nice} if whenever $i < i'$ both belong to $I \sem \ker (I)$ then there is $j \in I$
incomparable with $i$ such that $i' = i \lor j$.

(2) Let $I$ be a nice distributive almost-lattice. The set $\pure (I)$ of {\em pure pairs} is the smallest subset of
$\{ (i,j) : i \leq j, i,j \in I \}$ such that
\begin{alphenumerate}
\item $(i,i) \in \pure (I)$ for all $i \in I$,
\item if $i,j \in I$ are incomparable then $(i \land j, i)$ and $(i\land j, j)$ belong to $\pure (I)$, and if additionally $i \lor j$ exists
then also $(i , i \lor j)$ and $(j , i \lor j)$ belong to $\pure (I)$,
\item if $i \leq j \leq j' \leq k$ are in $I$ and $(i,k) \in \pure (I)$, then $(j, j') \in \pure (I)$. \hfill $\dashv$
\end{alphenumerate}
\end{defin}

Note that in (1) if $i < i'$ both belong to $I \sem \ker (I)$, then both belong either to $L \sem \ker (I)$ or to $R \sem \ker (I)$,
by Basic Lemma~\ref{basic-lemma} (ii). Niceness then means that any such pair $(i, i')$ is pure by (2)(b). (2)(c) should be 
construed as saying that the pure pairs define ``small" intervals in the partial order $\la I, \leq \ra$.

\begin{obs} \label{pure-obs}
Assume $\la I, \leq \ra$ is a nice distributive almost-lattice. Also assume $I$ has no maximal element, let $\ell \notin I$, $J = I \cup \{ \ell \}$
and stipulate $i \leq \ell$ for all $i \in I$. Then $\pure (J) = \pure (I) \cup \{ (i, \ell ) : i \in I \sem \ker (I) \} \cup \{ ( \ell,\ell) \}$.
\end{obs}

Note that by Observation~\ref{almost-lattice-basic}, $J$ is a lattice, in particular $L (J) = R(J) = \ker (J) = J$.

\begin{proof}
Obviously $\pure (I) \sub \pure(J)$, and if $i \in I \sem \ker (I)$, say $i \in L(I) \sem \ker (I)$, then there is $j \in R(I) \sem \ker (I)$ such that
$i \lor j$ does not exist in $I$ by Basic Lemma~\ref{basic-lemma}. Thus $i \lor j = \ell$ in $J$ and $(i , \ell) = (i , i \lor j)$ belongs to $\pure (J)$.
On the other hand, suppose $(i ,\ell) \in \pure (J)$ according to (b). Then there must be $j \in I$ incomparable with $i$ such that $\ell = i \lor j$.
Therefore $i \in L \sem \ker (I)$ and $j \in R \sem \ker (I)$ or vice-versa. Next, if $(i , \ell) \in \pure (J)$ according to (c), then, in fact,
$(i,\ell)$ already belongs to $\pure (J)$ according to (b). Finally notice that if $(i,\ell) \in \pure (J)$ and $i \leq j < j' < \ell$ then $(j,j') \in 
\pure (I)$ by niceness.
\end{proof}

We are ready for the main definition of this section.

\begin{defin} \label{defin4}
Let $I$ be a nice distributive almost-lattice, and assume $\bar \AA = \la \AA_i : i\in I \ra $
is a system of complete Boolean algebras with complete embeddings and correct projections.
We then have a system of names for quotient algebras $\la \dot \AA_{ij} : i \leq j \in I \ra$,
that is, $\AA_j$ factors as a two-step iteration, $\AA_j \cong \AA_i \star \dot \AA_{ij}$ (where $\dot \AA_{ij}$ is the
name for the trivial algebra $\{ \zero,\one \}$ if $i = j$). Say that $\la \AA_i : \; i\in I \ra$ is equipped with a 
{\em coherent system} $\la \dot \mu_{ij} : (i,j) \in \pure (I) \ra$ {\em of measures} if
\begin{enumerate}
\item[(A)] $\dot \mu_{ij}$ is forced  (by $\AA_i$)  to be a finitely additive strictly positive
   measure on $\dot \AA_{ij}$ (in particular, $\dot \mu_{ii}$  is the trivial measure on $\dot \AA_{ii}$),
\item[(B)] ({\em coherence}) given $i \leq j\leq k$ with $(i,k) \in \pure (I)$ and an $\AA_i$-name $\dot p_{ik}$ for a condition in
   $\dot \AA_{ik}$ which is forced to factor as $(\dot p_{ij}, \dot p_{jk})$ when considered
   as a condition in the two-step iteration $\dot \AA_{ij} \star \dot \AA_{jk}$, we have
   \[ \forces_i \dot \mu_{ik} (\dot p_{ik}) = \int_{\dot p_{ij}} \dot \mu_{jk} (\dot p_{jk}) 
   d\dot \mu_{ij}, \]
\item[(C)] ({\em correctness}) given incomparable $i,j \in I$ such that $i \lor j$ exists (recall that, by Basic Lemma~\ref{basic-lemma} and by
   Definition~\ref{defin3a}, this means either that $i, j, i\lor j \in L$ or $i,j, i \lor j \in R$ and that $(i\land j, i),
   (i\land j,j), (i, i \lor j), (j, i \lor j)$ all belong to $\pure (I)$) and $p_i = (p_{i \land j} , \dot p_{i \land j, i}) \in
   \AA_{ i \land j} \star \dot \AA_{i \land j, i} \cong \AA_i$ which can be rewritten as
   $p_i = (p_{i \land j} , \dot p_{j, i \lor j}) \in \AA_j \star \dot \AA_{j, i \lor j} \cong \AA_{i \lor j}$, we have
    \[ p_{i \land j} \forces_j \dot \mu_{j, i \lor j} ( \dot p_{j, i \lor j})  = \dot \mu_{i \land j, i} (\dot p_{i \land j, i } ). 
    \;\;\;\;\; \dashv \]
\end{enumerate}
\end{defin}

One may wonder why we only consider measures on ``small" quotients $\dot \AA_{ij}$ for pure pairs $(i,j)$ and not on larger
ones or on the whole algebra. The reason for this is that in this more general case we cannot prove the preservation
of coherence in the amalgamated limit when we extend the measures from the Boolean algebras generated by the corresponding
dense sets $\dot D_i$ to the cBa's $\dot \AA_{i\ell}$. This is explained in detail in the last step of the proof of
Main Lemma~\ref{amal-coherent} and the discussion thereafter (see, in particular, Counterexample~\ref{ctrex-product-measure}). On the other hand,
the ``small" quotient measures considered here are sufficient for us as we shall see in Section~\ref{Cohen}; namely, they will allow
us to define measures on the whole shattered iteration by a recursive procedure, whenever all the iterands carry a measure, see
Lemma~\ref{shattered-measures}.

A few comments concerning the representation of elements and their measures in quotient algebras are in order. 
If $i < j \in I$, and $p_j \in \AA_j$, we may rewrite this condition as $p_j  = (p_i , \dot p_{ij}) \in \AA_i \star \dot \AA_{ij} = \AA_j$
(and we have used this already here and there). This clearly means that $h_{ji} (p_j) = p_i$. Further, if $G_i$ is $\AA_i$-generic
with $p_i \in G_i$, then, in $V[G_i]$, we may identify the condition $p_{ij} =\dot p_{ij} [G_i] \in \dot \AA_{ij} [G_i]$ with the condition $p_j \in \AA_j / G_i$
(while we mostly use the former notation here, we used the latter when introducing the sets $D_i$, dense subsets of
quotients in the amalgamated limit, in Section~\ref{amalgamated}). Thus $\mu_{ij} (\dot p_{ij} [G_i])$ and $\mu_{ij} (p_j)$ will be the same thing
(if $(i,j) \in \pure (I)$).

Next, letting $i,j \in I$ such that $i \lor j$ exists and $p_i = (p_{i \land j}, \dot p_{i \land j, i} ) \in \AA_i = \AA_{i \land j} \star \dot \AA_{i \land j, i}$,
where $p_{i\land j} = h_{i,i\land j} (p_i)$, we have $h_{i,i\land j} (p_i) = h_{i \lor j, j} (p_i)$ by correctness, and, thinking of $p_i$ as a condition
in $\AA_{i \lor j}$, we may write $p_i = ( p_{i \land j}, \dot p_{i \land j, i} , \dot \one_{i , i \lor j} ) = ( p_{i\land j}, \dot \one_{i \land j, j}, \dot p_{j, i \lor j})
= (p_{i\land j}, \dot p_{j, i \lor j})$ where for the last item we think of $p_{i\land j}$ as a condition in $\AA_j$. Note, in particular,
that $ \dot p_{j, i \lor j}$, while formally an $\AA_j$-name, really is an $\AA_{i\land j}$-name.

Using these comments, we see that coherence and correctness imply the following basic properties of the system of measures:

\begin{obs}   \label{coherentbasic}
Assume $\bar \AA = \la \AA_i : i\in I\ra$ carries a coherent system $\la \dot \mu_{ij} : (i,j) \in \pure (I) \ra$ of measures. 
\begin{enumerate}
\item  If $(i,k) \in \pure (I)$, $i \leq j \leq k$, and $\forces_i \dot p_{ij} \in \dot\AA_{ij} \sub \dot\AA_{ik}$, then $\forces_i \dot\mu_{ij} (\dot p_{ij} ) = \dot\mu_{ik} (\dot p_{ij})$.
\item Let $i, j \in I$ be incomparable such that $i \lor j$ exists. Also assume $(p_i, p_j) $ belongs to $D$
as witnessed by $(i,j)$, let $p_{i \land j} = h_{i, i \land j} (p_i) = h_{j, i \land j} (p_j)$, and write 
$p_i = (p_{i \land j}, \dot p_{i \land j, i})$, $p_j = (p_{i\land j}, \dot p_{i \land j, j})$ as elements of the
corresponding two-step iterations. Then
\[ p_{i \land j} \forces_{i \land j}  \dot \mu_{i \land j, i \lor j} \left( \dot p_{i \land j, i} \cdot \dot p_{i \land j, j} \right) 
   = \dot \mu_{i \land j, i } (\dot p_{i \land j, i})  \cdot \dot \mu_{i \land j, j} (\dot p_{i \land j, j}). \]
\item {\em (coherence for integrals)} If $(i,k) \in \pure (I)$, $i \leq j \leq k$, $\dot x$ is an $\AA_k$-name for a real in $[0,1]$, and $\dot p_{ik} = (\dot p_{ij}, \dot p_{jk})$ is an 
   $\AA_i$-name for a condition in $\dot \AA_{ik}$, then 
   \[ \forces_i \int_{\dot p_{ik}} \dot x d \dot\mu_{ik} = \int_{\dot p_{ij}} \left( \int_{\dot p_{jk}} \dot x d\dot \mu_{jk} \right) d \dot \mu_{ij} . \]
\item {\em (correctness for integrals)} If $i,j \in I$ are incomparable such that $i \lor j$ exists, $\dot x$ is an $\AA_i$-name for a real in $[0,1]$, and $p_i = (p_{i \land j} , \dot p_{i \land j, i}) \in
   \AA_i$ is identified with $p_i = (p_{i \land j} , \dot p_{j, i \lor j}) \in  \AA_{i \lor j}$, then
   \[ p_{i \land j} \forces_j  \int_{\dot p_{j, i \lor j}} \dot x d\dot \mu_{j, i \lor j} = \int_{\dot p_{i\land j, i}} \dot x d \dot \mu_{i \land j, i} . \]
\end{enumerate}
\end{obs}

\begin{proof}
(i) Indeed, by coherence (B), $\forces_i  \dot\mu_{ik} (\dot p_{ij}) = \int_{\dot p_{ij}} \dot \mu_{jk} (\dot \one_{jk}) 
d\dot \mu_{ij} = \int_{\dot p_{ij}}   d\dot \mu_{ij} = \dot\mu_{ij} (\dot p_{ij})$.

(ii) We may identify the $\AA_{i \land j}$-name $\dot p_{i \land j, i}$ with an equivalent $\AA_j$-name $\dot p_{j, i \lor j}$.
We then get
\begin{eqnarray*}  p_{i \land j} \forces_{i \land j}  \dot \mu_{i \land j, i \lor j} \left( \dot p_{i \land j, j} \cdot \dot p_{j, i \lor j} \right) 
   &=& \int_{\dot p_{i\land j, j}} \dot \mu_{j, i \lor j} (\dot p_{j, i \lor j}) \; d \dot \mu_{i \land j, j} \\
   & = & \int_{\dot p_{i\land j, j}} \dot \mu_{i \land j, i } (\dot p_{i \land j, i}) \; d \dot \mu_{i \land j, j} 
   = \dot \mu_{i \land j, i } (\dot p_{i \land j, i})  \cdot \dot \mu_{i \land j, j} (\dot p_{i \land j, j}) 
\end{eqnarray*}
where the first equality is by coherence (B) and the second, by correctness (C).

(iii) Work in $V [G_i]$ where $G_i$ is $\AA_i$-generic over $V$. Then
\begin{eqnarray*} \int_{ p_{ik}} \dot x d \mu_{ik} & = &  \sup_{m \to \infty} \sum_{\ell = 0}^{2^m-1} \left[ \mu_{ik} \left( p_{ik} \cdot \BOOL {\ell + 1 \over 2^m } \geq \dot x > {\ell \over 2^m}
   \BOOR \right) \cdot {\ell \over 2^m} \right] \\ & = & \sup_{m \to \infty} \sum_{\ell = 0}^{2^m-1} \left[ \int_{p_{ij}} \dot\mu_{jk} \left( \dot p_{jk} \cdot \BOOL {\ell + 1 \over 2^m } \geq \dot x > {\ell \over 2^m}
   \BOOR \right) d\mu_{ij} \cdot {\ell \over 2^m} \right]  \\
   & = & \sup_{m \to \infty}   \int_{p_{ij}} \left[ \sum_{\ell = 0}^{2^m-1}\dot\mu_{jk} \left( \dot p_{jk} \cdot \BOOL {\ell + 1 \over 2^m } \geq \dot x > {\ell \over 2^m}
   \BOOR \right)  \cdot {\ell \over 2^m} \right]  d\mu_{ij} \\ & \leq &  \int_{p_{ij}} \sup_{m \to \infty}  \sum_{\ell = 0}^{2^m-1}\left[ \dot\mu_{jk} \left( \dot p_{jk} \cdot 
   \BOOL {\ell + 1 \over 2^m } \geq \dot x > {\ell \over 2^m} \BOOR\right)  \cdot {\ell \over 2^m} \right] d\mu_{ij}  \\
   & = & \int_{p_{ij}} \left( \int_{\dot p_{jk}} \dot x d \dot \mu_{jk} \right) d\mu_{ij}  = \int_{p_{ij}} \inf_{m \to \infty}  \sum_{\ell = 0}^{2^m-1}\left[ \dot\mu_{jk} \left( \dot p_{jk} \cdot 
   \BOOL {\ell + 1 \over 2^m } \geq \dot x > {\ell \over 2^m}\BOOR \right)  \cdot {\ell + 1\over 2^m} \right] d\mu_{ij}  \\
   & \leq & \inf_{m \to \infty}   \int_{p_{ij}} \left[ \sum_{\ell = 0}^{2^m-1}\dot\mu_{jk} \left( \dot p_{jk} \cdot  \BOOL {\ell + 1 \over 2^m } \geq \dot x > {\ell \over 2^m}
   \BOOR \right)  \cdot {\ell + 1\over 2^m} \right]  d\mu_{ij} = \int_{p_{ik}} \dot x d \mu_{ik} 
\end{eqnarray*}
where we use the definition of the integral in the first and last equality as well as those between the $\leq$, and coherence (B) in the second and last equality. The remaining
(in)equalities are obvious.

(iv) This has a proof similar to, but much simpler than, (iii), using correctness (C) instead. We omit the details.
\end{proof}

\begin{mainlem}[coherence in the amalgamated limit]   \label{amal-coherent}
Assume $\bar \AA = \la \AA_i : i\in I\ra$ carries a coherent system $\la \dot \mu_{ij}:
(i,j) \in \pure (I)  \ra$ of measures. Then $\bar \AA = \la \AA_i : i\in I \cup \{ \ell \} \ra$ carries a coherent system $\la \dot \mu_{ij}:
(i,j) \in \pure ( I \cup \{\ell \} )\ra$ of measures where $\ell > i$ for all $i \in I$ and $\AA_\ell  = \lim \amal_{i \in I} \AA_i$.
\end{mainlem}

\begin{proof}
Recall from Observation~\ref{pure-obs} that the only new elements of $\pure (I \cup \{ \ell \})$ are of the form $(i,\ell)$ where $i \in I \sem \ker (I)$ (except for the
trivial $(\ell,\ell)$). Thus it suffices to define the $\dot \mu_{i\ell}$ for such $i$ (in particular, if $I = \ker (I)$ there is nothing to prove).

Fix $i_0 \in I \sem \ker (I)$ and work in the $\AA_{i_0}$-generic extension $V[G_{i_0}]$. Without loss of generality $i_0 \in L$ ($i_0 \in R$ is analogous).
Let $(p_i , p_j) \in D_{i_0}$ where $i \in L$, $j \in R$, and $i \geq i_0$. This means that
$p_{i_0}: = h_{ii_0} (p_i)  \in G_{i_0}$ and, a fortiori, $p_{i_0 \land j} := h_{j, i_0 \land j} (p_j) \in G_{i_0}$. Decomposing $p_i$ and $p_j$ as
conditions in iterations as in earlier arguments, we may write (in the ground model $V$) $p_i = ( p_{i_0}, \dot p_{i_0i} ) = (p_{i_0} , \dot p_{i_0, i_0 \lor (i \land j)} , \dot p_{ i_0 \lor (i \land j), i} )$
and $p_j = (p_{i_0 \land j} , \dot p_{i_0 \land j, j} ) = (p_{i_0 \land j} , \dot p_{i_0 \land j, i \land j} , \dot p_{i \land j, j} )$. 
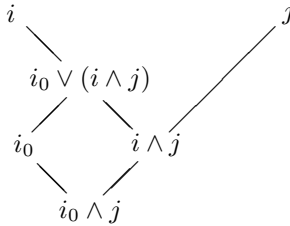
\[
\begin{picture}(105,75)(0,0)
\put(30,0){\makebox(0,0){$i_0 \land j$}}
\put(5,25){\makebox(0,0){$ i_0$}}
\put(55,25){\makebox(0,0){$i \land j$}}
\put(30,50){\makebox(0,0){$i_0 \lor (i \land j)$}}
\put(0,75){\makebox(0,0){$i$}}
\put(105,75){\makebox(0,0){$j$}}
\put(8,30){\line(1,1){13}}
\put(10,18){\line(1,-1){11}}
\put(35,6){\line(1,1){13}}
\put(35,42){\line(1,-1){11}}
\put(6,70){\line(1,-1){13}}
\put(60,30){\line(1,1){40}}
\end{picture}
\]
In $V [G_{i_0}]$, $p_i$ gets identified with
$p_{i_0i}  = (p_{i_0, i_0 \lor (i \land j)} , \dot p_{ i_0 \lor (i \land j), i} )$ and $p_j$, with $p_{i_0 \land j, j} = (p_{i_0 \land j, i \land j} , \dot p_{i \land j, j} )$. 
Thus we define
\[ \mu_{i_0 \ell} (p_i,p_j) = \mu_{i_0 \ell} (p_{i_0i}, p_{i_0 \land j, j}) : =  \int_{p_{i_0, i_0 \lor ( i \land j)}} \dot\mu_{i_0 \lor (i\land j), i} \left(\dot p_{i_0 \lor (i\land j), i }\right) 
\cdot \dot \mu_{i \land j, j} \left( \dot p_{i \land j, j}\right) \; d\mu_{i_0, i_0 \lor (i \land j)} \]
Note that the pairs $(i_0, i_0 \lor (i \land j))$, $(i_0 \lor (i \land j), i)$, and $(i \land j , j)$ are indeed pure, the former two because $i_0$ and $i$ both belong to
$L \sem \ker (I)$, and the latter, by (2)(a) or (2)(b) in Definition~\ref{defin3a} (note that either $j \leq i$ or $i$ and $j$ are incomparable). 
The first thing we need to verify is that this expression is \underline{\em well-defined}, for different members of $D_{i_0}$ may be equivalent
as conditions in $\AA_{i_0 \ell}$.

Assume that $(q_{i'}, q_{j'})$  is a condition in $D_{i_0}$, as witnessed by $i' \geq i_0$ in $L$  and $j'$ in $R$, 
which is equivalent to $(p_i, p_j)$. By Observation~\ref{basicobsquotient} (iii), we may assume that $i' \geq i$ and $j' \geq j$. By part (ii) of this Observation, there is
$r \in G_{i_0}$ such that $q_{i'} \cdot r = p_i \cdot h_{j, i' \land j} (p_j) \cdot r$ and $q_{j'} \cdot h_{i_0, i_0 \land j'} (r) = p_j \cdot h_{i,i \land j'} (p_i) \cdot h_{i_0, i_0 \land j'} (r)$.
Considering the corresponding $\AA_{ i \land j}$-names $\dot q_{i \land j, i' \land j}$ and $\dot q_{i \land j, i \land j'}$, that is,
$h_{j, i' \land j} (p_j) = (p_{i \land j} , \dot q_{i \land j, i' \land j} )$ and $h_{i, i\land j'} (p_i) = (p_{i \land j} , \dot q_{i \land j, i \land j'} )$,
this means that, in $V[G_{i_0}]$, $q_{i_0 i'} = ( p_{i_0i}, \dot q_{i\land j, i' \land j})$ and $q_{i_0 \land j', j'} = (p_{i_0 \land j , j}, \dot q_{i \land j , i \land j'})$. 
Since we may argue in two steps, to simplify things first consider the case $j = j'$. Thus $q_{i_0 \land j, j} = p_{i_0 \land j , j}$.
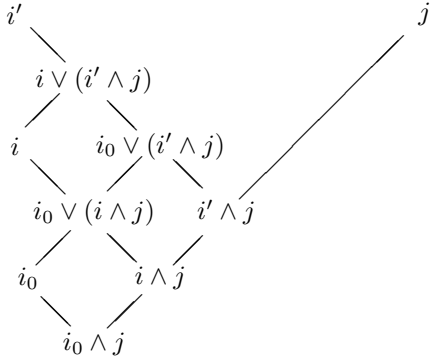
\[
\begin{picture}(155,125)(0,0)
\put(30,0){\makebox(0,0){$i_0 \land j$}}
\put(5,25){\makebox(0,0){$ i_0$}}
\put(55,25){\makebox(0,0){$i \land j$}}
\put(30,50){\makebox(0,0){$i_0 \lor (i \land j)$}}
\put(80,50){\makebox(0,0){$i' \land j$}}
\put(0,75){\makebox(0,0){$i$}}
\put(55,75){\makebox(0,0){$i_0 \lor (i' \land j)$}}
\put(30,100){\makebox(0,0){$i \lor (i' \land j)$}}
\put(0,125){\makebox(0,0){$i'$}}
\put(155,125){\makebox(0,0){$j$}}
\put(8,30){\line(1,1){13}}
\put(10,18){\line(1,-1){11}}
\put(35,6){\line(1,1){13}}
\put(35,42){\line(1,-1){11}}
\put(6,70){\line(1,-1){13}}
\put(60,30){\line(1,1){11}}
\put(4,82){\line(1,1){13}}
\put(6,120){\line(1,-1){13}}
\put(35,58){\line(1,1){13}}
\put(35,92){\line(1,-1){11}}
\put(85,55){\line(1,1){62}}
\put(60,70){\line(1,-1){11}}
\end{picture}
\]
Note $p_{i_0 \land j, j} = ( p_{i_0 \land j, i \land j }, \dot q_{i \land j, i' \land j}, \dot q_{i' \land j, j})$. 
We may identify $\dot q_{i \land j, i' \land j}$ with an $\AA_{i_0, i_0 \lor (i \land j)}$-name $\dot q_{i_0 \lor ( i \land j), i_0 \lor ( i' \land j)}$, and
$\dot p_{i_0 \lor (i \land j), i}$,  with an $\AA_{i_0,i_0  \lor (i ' \land j)}$-name $\dot q_{i_0 \lor (i' \land j), i'}$. Then 
$q_{i_0 i'} = (p_{i_0, i_0 \lor (i \land j)}, \dot p_{i_0 \lor (i \land j), i} , \dot q_{i \land j, i' \land j}) = $  \\
$(p_{i_0, i_0 \lor (i\land j)} , \dot q_{i_0 \lor ( i \land j), i_0 \lor ( i' \land j)}, \dot q_{i_0 \lor (i' \land j), i'})$, and
we compute 
\begin{eqnarray*}  
&&   \mu_{i_0 \ell} (p_{i_0i}, p_{i_0 \land j, j})  =   \int_{p_{i_0, i_0 \lor ( i \land j)}} \dot\mu_{i_0 \lor (i\land j), i} \left(\dot p_{i_0 \lor (i\land j), i }\right) 
\cdot \dot \mu_{i \land j, j} \left( \dot p_{i \land j, j}\right) \; d\mu_{i_0, i_0 \lor (i \land j)} \\
   & = &  \int_{p_{i_0, i_0 \lor ( i \land j)}} \dot\mu_{i_0 \lor (i\land j), i} \left(\dot p_{i_0 \lor (i\land j), i }\right)  \cdot \left( \int_{\dot q_{i \land j, i' \land j} }  \dot \mu_{i' \land j, j} \left(\dot q_{i' \land j,j}\right)  
      \; d \dot \mu_{i \land j,i' \land  j}  \right)  d\mu_{i_0, i_0 \lor (i \land j)}  \\
   & = &  \int_{p_{i_0, i_0 \lor ( i \land j)}}    \left(    \int_{\dot q_{i \land j, i' \land j} }     \dot\mu_{i_0 \lor (i\land j), i} \left(\dot p_{i_0 \lor (i\land j), i }\right)  \cdot    \dot \mu_{i' \land j, j} 
   \left(\dot q_{i' \land j,j}\right)  \; d \dot \mu_{i \land j,i' \land  j}  \right)  d\mu_{i_0, i_0 \lor (i \land j)} \\
   & = & \int_{p_{i_0, i_0 \lor ( i \land j)}}    \left(    \int_{\dot q_{i_0 \lor (i \land j),i_0 \lor ( i' \land j)} }     \dot \mu_{i_0 \lor (i' \land j), i'} \left(\dot q_{i_0 \lor (i '\land j),i'}\right)  \cdot    
        \dot \mu_{i' \land j, j} \left(\dot q_{i' \land j,j}\right)  \;    d \dot \mu_{i_0 \lor (i \land j), i_0 \lor (i' \land  j)}  \right)  d\mu_{i_0, i_0 \lor (i \land j)} \\
    & = &  \int_{q_{i_0, i_0 \lor (i '\land j)}}        \dot \mu_{i_0 \lor (i' \land j), i'} \left(\dot q_{i_0 \lor(i '\land j),i'}\right)  \cdot    \dot \mu_{i' \land j, j} \left(\dot q_{i' \land j,j}\right)  \; 
        d\mu_{i_0, i_0 \lor (i '\land j)}   \;\; = \;\; \mu_{i_0 \ell} (q_{i_0i'}, q_{i_0 \land j,j}) 
 \end{eqnarray*}
where the second and fifth equal signs use  coherence ((B) and Observation~\ref{coherentbasic} (iii), respectively), while the fourth is by correctness ((C) and Observation~\ref{coherentbasic} (iv)). 
(We leave it to the reader to verify that all respective pairs are again pure.)

Next consider the case $i = i'$. Thus $q_{i_0i} = p_{i_0i}$.
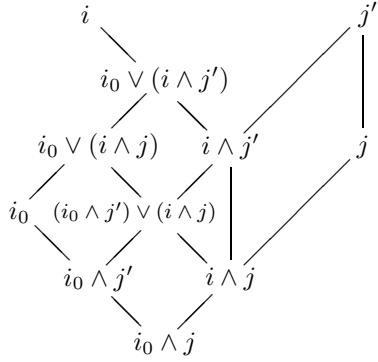
\[
\begin{picture}(132,125)(0,0)
\put(55,0){\makebox(0,0){$i_0 \land j$}}
\put(30,25){\makebox(0,0){$ i_0 \land j'$}}
\put(80,25){\makebox(0,0){$i \land j$}}
\put(0,50){\makebox(0,0){$i_0$}}
\put(55,50){\makebox(-20,0){{\footnotesize $(i_0 \land j') \lor (i \land j)$ }}}
\put(80,75){\makebox(0,0){$i \land j'$}}
\put(30,75){\makebox(0,0){$i_0 \lor (i \land j)$}}
\put(130,75){\makebox(0,0){$j$}}
\put(55,100){\makebox(0,0){$i_0 \lor (i \land j')$}}
\put(25,125){\makebox(0,0){$i$}}
\put(132,125){\makebox(0,0){$j'$}}
\put(33,30){\line(1,1){13}}
\put(35,18){\line(1,-1){11}}
\put(60,6){\line(1,1){13}}
\put(60,42){\line(1,-1){11}}
\put(6,45){\line(1,-1){13}}
\put(4,55){\line(1,1){13}}
\put(31,120){\line(1,-1){13}}
\put(35,82){\line(1,1){13}}
\put(35,67){\line(1,-1){11}}
\put(60,92){\line(1,-1){11}}
\put(60,55){\line(1,1){13}}
\put(80,32){\line(0,1){35}}
\put(130,82){\line(0,1){35}}
\put(85,30){\line(1,1){40}}
\put(85,80){\line(1,1){40}}
\end{picture}
\]
Note that $p_{i_0i} = (p_{i_0, i_0 \lor ( i \land j)}, \dot q_{i_0 \lor (i \land j) , i_0 \lor (i \land j') } , \dot q_{i_0 \lor (i \land j') , i } )$. We may identify
$\dot p_{i \land j, j}$ with an $\AA_{i_0 \land j' , i \land j'}$-name $\dot q_{i \land j' , j'}$ and then can write
$q_{i_0 \land j' , j'} = (p_{i_0 \land j, i \land j}, \dot p_{i \land j, j}, \dot q_{i \land j, i \land j'}) =
( p_{i_0 \land j, i\land j} , \dot q_{i \land j, i \land j'}, \dot q_{i \land j' , j'} )$. We compute
\begin{eqnarray*}
&&   \mu_{i_0 \ell} (p_{i_0i}, p_{i_0 \land j, j})  =   \int_{p_{i_0, i_0 \lor ( i \land j)}} \dot\mu_{i_0 \lor (i\land j), i} \left(\dot p_{i_0 \lor (i\land j), i }\right) 
\cdot \dot \mu_{i \land j, j} \left( \dot p_{i \land j, j}\right) \; d\mu_{i_0, i_0 \lor (i \land j)} \\
   & = &  \int_{p_{i_0, i_0 \lor ( i \land j)}} \left(   \int_{\dot q_{i_0 \lor (i \land j), i_0 \lor (i \land j')}} \dot\mu_{i_0 \lor (i\land j'), i} \left(\dot q_{i_0 \lor (i\land j'), i }\right)  
   \; d\dot\mu_{i_0 \lor (i \land j), i_0 \lor (i \land j')} \right)  \cdot  \dot  \mu_{i \land j, j} \left( \dot p_{i \land j, j}\right) \; d\mu_{i_0, i_0 \lor (i \land j)} \\
   & = &  \int_{p_{i_0, i_0 \lor ( i \land j)}} \left(   \int_{\dot q_{i_0 \lor (i \land j), i_0 \lor (i \land j')}} \dot\mu_{i_0 \lor (i\land j'), i} \left(\dot q_{i_0 \lor (i\land j'), i }\right)  
   \cdot   \dot \mu_{i \land j, j} \left( \dot p_{i \land j, j}\right)  \; d\dot\mu_{i_0 \lor (i \land j), i_0 \lor (i \land j')} \right)   \; d\mu_{i_0, i_0 \lor (i \land j)} \\
   & = &  \int_{p_{i_0, i_0 \lor ( i \land j)}} \left(   \int_{\dot q_{i_0 \lor (i \land j), i_0 \lor (i \land j')}} \dot\mu_{i_0 \lor (i\land j'), i} \left(\dot q_{i_0 \lor (i\land j'), i }\right)  
   \cdot  \dot  \mu_{i \land j', j'} \left( \dot q_{i \land j', j'}\right)  \; d\dot\mu_{i_0 \lor (i \land j), i_0 \lor (i \land j')} \right)   \; d\mu_{i_0, i_0 \lor (i \land j)} \\
    & = &  \int_{q_{i_0, i_0 \lor (i \land j')}}    \dot\mu_{i_0 \lor (i\land j'), i} \left(\dot q_{i_0 \lor (i\land j'), i }\right)  
   \cdot  \dot  \mu_{i \land j', j'} \left( \dot q_{i \land j', j'}\right)  \; d\mu_{i_0, i_0 \lor (i \land j')}   \;\; = \;\; \mu_{i_0 \ell} (q_{i_0i}, q_{i_0 \land j' ,j'}) 
 \end{eqnarray*}
where the second and fifth equal signs use  coherence ((B) and Observation~\ref{coherentbasic} (iii), respectively), while the fourth is by correctness (C). 
This completes the proof of well-definedness.

Obviously $\mu_{i_0 \ell}$ is \underline{\em strictly positive} on $D_{i_0}$ and $\mu_{i_0 \ell } (\one) = 1$. For \underline{\em finite additivity}, we need to consider finite sums of pairwise
incompatible elements of $D_{i_0}$ as in Preliminary Lemma~\ref{firstextension} (iv) because the $D_{i_0}$ are not closed under sums of
two incompatible elements. By Observation~\ref{basicobs} (vi), we may assume these finitely many conditions
have the same witnesses $i \geq i_0$, $i \in L$, and $j\in R$. So we have $(p_i, p_j) , (p_i^m, p_j^m) \in D$, $m < n$, such that
\[ r \forces `` (p_i^m , p_j^m) \mbox{ are pairwise incompatible  and } \sum_{m < n} (p^m_i,p^m_j) = (p_i,p_j) "\] 
for some $r \in \AA_{i_0}$.  In $V[G_{i_0}]$ with $r \in G_{i_0}$, 
this means that the $(p^m_{i_0i} , p^m_{ i_0 \land j , j } ) $ are pairwise incompatible and $\sum_{m < n} (p^m_{i_0i} , p^m_{ i_0 \land j , j } ) = ( p_{i_0i} , p_{i_0 \land j, j} )$.
Since $\one - p^m_{ i_0, i_0 \lor (i \land j)} \forces \dot p^m_{i_0 \lor (i \land j),i }  = \zero$, we get
\begin{eqnarray*}
 \mu_{i_0 \ell} (p^m_{i_0i}, p^m_{i_0 \land j, j})  & =  &  \int_{p_{i_0, i_0 \lor ( i \land j)}^m} \dot\mu_{i_0 \lor (i\land j), i} \left(\dot p^m_{i_0 \lor (i\land j), i }\right) 
\cdot \dot \mu_{i \land j, j} \left( \dot p^m_{i \land j, j}\right) \; d\mu_{i_0, i_0 \lor (i \land j)} \\
&=& \int_{p_{i_0, i_0 \lor ( i \land j)}} \dot\mu_{i_0 \lor (i\land j), i} \left(\dot p^m_{i_0 \lor (i\land j), i }\right) 
\cdot \dot \mu_{i \land j, j} \left( \dot p^m_{i \land j, j}\right) \; d\mu_{i_0, i_0 \lor (i \land j)}. 
\end{eqnarray*}

\begin{sclaim}   \label{amal-coherent-claim1}
\[ p_{i_0, i_0 \lor (i \land j)} \forces \sum_{m < n} \dot \mu_{i_0 \lor (i\land j), i} \left(\dot p^m_{i_0 \lor (i\land j), i }\right)  \cdot \dot \mu_{i \land j, j} \left( \dot p^m_{i \land j, j}\right)
= \dot \mu_{i_0 \lor (i\land j), i} \left(\dot p_{i_0 \lor (i\land j), i }\right)  \cdot \dot \mu_{i \land j, j} \left( \dot p_{i \land j, j}\right). \]
\end{sclaim}

\begin{proof}
To see this, let $G_{i_0, i_0 \lor (i\land j) }$ be $\AA_{i_0, i_0 \lor (i\land j) }$-generic over $V[G_{i_0}]$ with $p_{i_0, i_0 \lor (i \land j)} \in G_{i_0, i_0 \lor (i\land j) }$  and step into the generic extension
$V [G_{i_0}] [G_{i_0, i_0 \lor (i\land j) }] = V [ G_{i_0 \lor (i \land j)}] $. Then the $( p^m_{ i_0 \lor (i\land j), i} , p^m_{i \land j, j})$ are pairwise incompatible,
with the possibility of some of them being $\zero$, and their sum is $(p_{i_0 \lor (i\land j), i} , p_{i\land j, j})$. Note that incompatibility means that any two of them  are
either incompatible in the first coordinate or in the second. This means that we can find conditions $q^k_{i_0 \lor (i \land j) , i} \in \AA_{i_0 \lor (i \land j) , i}$, $k < \ell$, and $q^{k'}_{i \land j, j}
\in \AA_{i \land j, j}$, $k' < \ell '$, and sets $A_m \sub \ell, B_m \sub \ell '$, $m < n$, such that
\begin{enumerate}
\item the $q^k_{i_0 \lor (i \land j) , i}$ are pairwise incompatible, and so are the $q^{k'}_{i \land j, j}$,
\item $\sum_{k < \ell } q^k_{i_0 \lor (i \land j) , i} = p_{i_0 \lor( i \land j) , i}$ and $\sum_{k' < \ell '} q^{k'}_{i \land j, j} = p_{i\land j, j}$,
\item for $m \neq m'$, either $A_m \cap A_{m'} = \emptyset$ or $B_m \cap B_{m'} = \emptyset$,
\item either $A_m = B_m = \emptyset$ or both are non-empty,
\item $\sum_{k \in A_m } q^k_{i_0 \lor (i \land j) , i} = p^m_{i_0 \lor( i \land j) , i}$ and $\sum_{k' \in B_m} q^{k'}_{i \land j, j} = p^m_{i\land j, j}$,
\item for each pair $k , k'$ there is a (necessarily unique) $m$ such that $k \in A_m$ and $k' \in B_m$.
\end{enumerate}
Now, using the finite additivity of the measures $\mu_{i_0 \lor (i \land j) , i}$ and $\mu_{i \land j, j}$, we see
\begin{eqnarray*}
\mu_{i_0 \lor (i \land j), i} \left(p_{i_0 \lor (i \land j), i} \right) \cdot \mu_{i\land j, j} \left( p_{i \land j, j } \right) & \stackrel{\textstyle{\mathrm{(i),(ii)}}}{=}  & 
   \left( \sum_{k < \ell} \mu_{i_0 \lor (i \land j), i} \left(q^k_{i_0 \lor (i \land j), i} \right)  \right) \cdot \left( \sum_{k' < \ell'}  \mu_{i\land j, j} \left( q^{k'}_{i \land j, j } \right) \right) \\
  & = & \sum_{k < \ell} \sum_{k' < \ell'} \left( \mu_{i_0 \lor (i \land j), i} \left(q^k_{i_0 \lor (i \land j), i} \right)   \cdot \mu_{i\land j, j} \left( q^{k'}_{i \land j, j } \right) \right) \\
  & \stackrel{\textstyle{\mathrm {(vi)}}}{=} & \sum_{m < n} \sum_{k \in A_m} \sum_{k' \in B_m} \left( \mu_{i_0 \lor (i \land j), i} \left(q^k_{i_0 \lor (i \land j), i} \right)  
      \cdot \mu_{i\land j, j} \left( q^{k'}_{i \land j, j } \right) \right) \\
  & \stackrel{\textstyle{\mathrm{ (v)}}}{=} & \sum_{m < n} \mu_{i_0 \lor (i \land j), i} \left(p^m_{i_0 \lor (i \land j), i} \right) \cdot \mu_{i\land j, j} \left( p^m_{i \land j, j } \right) .
\end{eqnarray*}
This ends the proof of the claim.
\end{proof}

Thus we obtain 
\begin{eqnarray*}
\mu_{i_0 \ell} ( p_{i_0 i} , p_{i_0 \land j, j} ) & = & \int_{p_{i_0, i_0 \lor ( i \land j)}} \dot\mu_{i_0 \lor (i\land j), i} \left(\dot p_{i_0 \lor (i\land j), i }\right) 
\cdot \dot \mu_{i \land j, j} \left( \dot p_{i \land j, j}\right) \; d\mu_{i_0, i_0 \lor (i \land j)} \\
& = &  \int_{p_{i_0, i_0 \lor ( i \land j)}} \left( \sum_{m < n} \dot\mu_{i_0 \lor (i\land j), i} \left(\dot p^m_{i_0 \lor (i\land j), i }\right) 
\cdot \dot \mu_{i \land j, j} \left( \dot p^m_{i \land j, j}\right)\right) \; d\mu_{i_0, i_0 \lor (i \land j)} \\
& = & \sum_{m < n} \left(  \int_{p_{i_0, i_0 \lor ( i \land j)}}  \dot\mu_{i_0 \lor (i\land j), i} \left(\dot p^m_{i_0 \lor (i\land j), i }\right) 
\cdot \dot \mu_{i \land j, j} \left( \dot p^m_{i \land j, j}\right) \; d\mu_{i_0, i_0 \lor (i \land j)} \right) = \sum_{m < n} \mu_{i_0 \ell} ( p^m_{i_0 i} , p^m_{i_0 \land j, j} ) 
\end{eqnarray*}
and the proof of finite additivity is complete.

To show \underline{\em coherence} (B) of the measures, let $i_0 < i_1$ be elements of $L \sem \ker (I)$, and work in $V [ G_{i_0}]$ (the
case $i_0, i_1 \in R \sem \ker (I)$ is analogous). We have a condition $(p_{i_0i} , p_{i_0 \land j, j} ) \in D_{i_0}$, witnessed
by $i \geq i_1$, $i \in L$, and $j\in R$. We may write $p_{i_0i} = ( p_{i_0i_1}, \dot p_{i_1 i})$ and $p_{i_0 \land j,j} = ( p_{i_0 \land j , i_1 \land j} , \dot p_{i_1 \land j, j} )$.
Further decompositions are defined similarly in a straightforward manner.
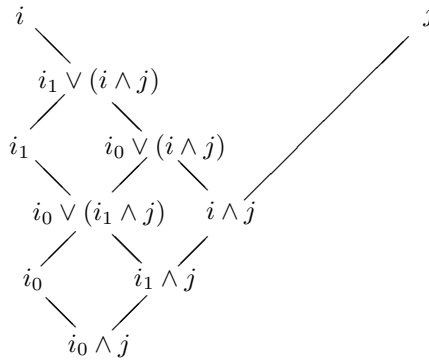
\[
\begin{picture}(155,125)(0,0)
\put(30,0){\makebox(0,0){$i_0 \land j$}}
\put(5,25){\makebox(0,0){$ i_0$}}
\put(55,25){\makebox(0,0){$i_1 \land j$}}
\put(30,50){\makebox(0,0){$i_0 \lor (i_1 \land j)$}}
\put(80,50){\makebox(0,0){$i \land j$}}
\put(0,75){\makebox(0,0){$i_1$}}
\put(55,75){\makebox(0,0){$i_0 \lor (i \land j)$}}
\put(30,100){\makebox(0,0){$i_1 \lor (i \land j)$}}
\put(0,125){\makebox(0,0){$i$}}
\put(155,125){\makebox(0,0){$j$}}
\put(8,30){\line(1,1){13}}
\put(10,18){\line(1,-1){11}}
\put(35,6){\line(1,1){13}}
\put(35,42){\line(1,-1){11}}
\put(6,70){\line(1,-1){13}}
\put(60,30){\line(1,1){11}}
\put(4,82){\line(1,1){13}}
\put(6,120){\line(1,-1){13}}
\put(35,58){\line(1,1){13}}
\put(35,92){\line(1,-1){11}}
\put(85,55){\line(1,1){62}}
\put(60,70){\line(1,-1){11}}
\end{picture}
\]
We need to establish $\mu_{i_0\ell} (p_{i_0 i} , p_{i_0\land j, j}) = \int_{p_{i_0i_1}} \dot \mu_{i_1 \ell} ( \dot p_{i_1i} , \dot p_{i_1 \land j, j} ) \; d\mu_{i_0i_1}$. But indeed,
applying coherence of the original system ((B) and Observation~\ref{coherentbasic} (iii)) several times, we obtain
\begin{eqnarray*}
  \mu_{i_0 \ell} (p_{i_0i}, p_{i_0 \land j, j})  & = &  \int_{p_{i_0, i_0 \lor ( i \land j)}} \dot\mu_{i_0 \lor (i\land j), i} \left(\dot p_{i_0 \lor (i\land j), i }\right) 
\cdot \dot \mu_{i \land j, j} \left( \dot p_{i \land j, j}\right) \; d\mu_{i_0, i_0 \lor (i \land j)} \\
& = &  \int_{p_{i_0, i_0 \lor ( i_1 \land j)}} \left[  \int_{\dot p_{i_0 \lor (i_1 \land j), i_0 \lor ( i \land j)}} \left(    \int_{\dot p_{i_0 \lor (i \land j) , i_1 \lor ( i \land j)}}
   \dot\mu_{i_1 \lor (i\land j), i} \left(\dot p_{i_1 \lor (i\land j), i }\right) \; d\dot\mu_{i_0 \lor (i \land j) , i_1 \lor ( i \land j)} \right) \right. \\ &&\hskip 3truecm \dot \mu_{i \land j, j} \left( \dot p_{i \land j, j}\right) \; 
   d\dot\mu_{i_0\lor (i_1 \land j), i_0 \lor (i \land j)} \Bigg]  d\mu_{i_0, i_0 \lor (i_1 \land j)}  \\
& = &  \int_{p_{i_0, i_0 \lor ( i_1 \land j)}}  \left(    \int_{\dot p_{i_0 \lor (i_1 \land j), i_1}}   \int_{\dot p_{i_1 , i_1 \lor ( i \land j)}} \dot\mu_{i_1 \lor (i\land j), i} \left(\dot p_{i_1 \lor (i\land j), i }\right) \right. \\
   &&\hskip 3truecm \dot \mu_{i \land j, j} \left( \dot p_{i \land j, j}\right) \; d\dot\mu_{i_1 , i_1 \lor ( i \land j)} \; d\dot\mu_{i_0\lor (i_1 \land j), i_1}  \Bigg)  \; d\mu_{i_0, i_0 \lor (i_1 \land j)}  \\
& = &  \int_{p_{i_0 i_1}}  \int_{\dot p_{i_1 , i_1 \lor ( i \land j)}} \dot\mu_{i_1 \lor (i\land j), i} \left(\dot p_{i_1 \lor (i\land j), i }\right) \cdot
   \dot \mu_{i \land j, j} \left( \dot p_{i \land j, j}\right) \; d\dot\mu_{i_1 , i_1 \lor ( i \land j)}   \; d\mu_{i_0i_1} \\
& = & \int_{p_{i_0i_1}} \dot \mu_{i_1 \ell} ( \dot p_{i_1 i} , \dot p_{i_1 \land j, j} ) \; d\mu_{i_0i_1} 
\end{eqnarray*}
as required.

To see that \underline{\em correctness} (C) is preserved, it suffices to consider the case where $i,j \in I$ such that $i \lor j$ does not exist in $I$ so that
$i \lor j = \ell$ in $I \cup \{ \ell \}$ (cf. Observation~\ref{almost-lattice-basic}). Thus, by Basic Lemma~\ref{basic-lemma}, without loss of generality, $i \in L \sem \ker (I)$
and $j \in R \sem \ker (I)$. We need to show that if we rewrite $p_j = (p_{i \land j}, \dot p_{i \land j, j})$ as $(p_{i \land j}, \dot p_{i \ell})$, we have
$p_{i \land j} \forces_i \dot \mu_{i \ell} (\dot p_{i \ell} ) = \dot \mu_{i \land j, j} (\dot p_{i \land j, j})$. Work in $V[G_i]$ with $p_{i \land j} \in G_i$.
Note that $p_{i \ell}$ is identified with the pair $(\one_{ii}, p_{i \land j, j})$ and, given that we are in the special case $i_0 = i$ in the definition of $\mu_{i\ell}$,
we see
\[ \mu_{i\ell} (p_{i \ell}) = \mu_{i\ell} ( \one_{ii}, p_{i \land j, j} ) = \int_{\one_{ii}} \dot \mu_{ii} (\dot \one_{ii}) \cdot \dot \mu_{i\land j, j} (\dot p_{i \land j, j}) \; d \mu_{ii} =  \mu_{i\land j, j} ( p_{i \land j, j}), \]
as required, because the names are actual objects and  $\mu_{ii} ( \one_{ii}) = 1$. Note that we do not need to deal with correctness anymore when
we extend the measures because the elements it refers to necessarily belong to  the $D_i$.

We come to \underline{\em measure extension}.
First, the measures satisfy all the required properties on  the $D_{i}$, and we can {\em extend them to the Boolean algebras  $\la D_{i} \ra$} generated 
by these sets. Indeed, by Preliminary Lemma~\ref{properties-of-D} (iii) and by the finite additivity proved above, the $D_i$ satisfy the assumptions of Preliminary Lemma~\ref{firstextension}
in $V[G_{i}]$. Thus we obtain finitely additive measures on  the $\la D_{i} \ra$, which we also denote by
$\mu_{i \ell}$. We need to show coherence (B) is preserved. 

Assume $i < j < \ell$ belong to $I \sem \ker (I)$ 
and work in $V[G_i]$. Let $a_{i\ell} = ( a_{ij} , \dot a_{j\ell}) \in \la D_i \ra = \AA_{ij} \star \la \dot D_j \ra$. By the proof of Preliminary Lemma~\ref{firstextension},
we know $a_{i\ell}$ is a disjoint sum of finitely many elements of $D_i$, $a_{i\ell } = \sum_{m < n} b_{i\ell}^m$ with $b_{i\ell}^m \in D_i$. Write 
$b_{i\ell}^m = (b_{ij}^m, \dot b_{j\ell}^m )$. Then
\begin{eqnarray*}
\mu_{i\ell} (a_{i\ell}) & = & \sum_{m < n} \mu_{i\ell} \left(b_{i\ell}^m \right)=  \sum_{m < n} \left( \int_{b_{ij}^m} \dot \mu_{j\ell} \left( \dot b_{j\ell}^m \right) d\mu_{ij} \right) \\
   & = &      \int_{a_{ij}} \left(   \sum_{m<n}  \dot \mu_{j\ell} \left( \dot b_{j\ell}^m \right)   \right) d\mu_{ij} =  \int_{a_{ij}}  \dot \mu_{j\ell} \left(\dot a_{j\ell} \right)   d\mu_{ij} 
\end{eqnarray*}
where the second equal sign comes from the coherence of the $D_i$ established above, and the third uses $\one - b_{ij}^m \forces \dot b_{j\ell}^m = \zero$. 

We finally need to {\em extend the measures  to the $\AA_{i \ell}$}. Let $\EE_\ell$ be a Boolean algebra with
$\la D \ra \sub \EE_\ell \sub \ro (D) = \AA_\amal$. We then put $\EE_i = \AA_i$ and $\dot \EE_{ij} = \dot \AA_{ij}$ for $i < j \in I$. Furthermore, if $G_i$ is $\EE_i$-generic
over $V$ for $i \in I$, define the quotient Boolean algebra $\EE_{i\ell} = \EE_\ell / G_i = \{  [a_\ell ] : a_\ell \in \EE_\ell \}$ where $a_\ell$ and $b_\ell \in \EE_\ell$
are equivalent if $\one - h_{\ell i} (a_\ell \triangle b_\ell) \in G_i$ and $[a_\ell]$ is the equivalence class of $a_\ell$. Here $\triangle$ denotes symmetric difference:
$a_\ell \triangle b_\ell = ( a_\ell - b_\ell) + (b_\ell - a_\ell)$. Note that $h_{\ell i} (a_\ell \triangle b_\ell) \in \EE_i$ because the latter is a cBa.
Also note that $[a_\ell] = \zero$ in $\EE_{i\ell}$ iff $h_{\ell i} (a_\ell) \notin G_i$.
Writing $a_\ell = (a_i, \dot a_{i\ell})$ with $h_{\ell i} (a_\ell) = a_i $, we may identify $[a_\ell] \in \EE_{i\ell}$ with $a_{i \ell} = \dot a_{i\ell} [ G_i]$. 
Thus, from $\EE_\ell$, we obtain a system $\la \dot \EE_{ij} : i < j \leq \ell \ra$ of Boolean algebras.

By a standard application of Zorn's Lemma (see also Lemma~\ref{secondextension} and its proof), it suffices to prove that if 
$\bar \EE = \la \EE_i : i \in  I \cup \{ \ell \} \ra$ is as in the preceding paragraph and carries a coherent system $\la \dot \mu_{ij} = \dot \mu_{ij}^\EE : (i,j) \in \pure ( I \cup \{ \ell \}) \ra$ of measures,
and if $b_\ell \in \AA_\ell \sem \EE_\ell$, then, letting $\EE_\ell ' = \la \EE_\ell \cup \{ b_\ell \} \ra$, $\bar \EE ' = \la \EE_i ' : i \in I \cup \{ \ell \} \ra$
carries a coherent system $\la \dot \mu_{ij} ' = \dot \mu_{ij}^{\EE '} : (i,j) \in \pure (I \cup \{ \ell\} ) \ra$ of measures extending the original system. (Note here that $\EE_i ' = \EE_i$ for $i < \ell$.)

Fix $i \in I \sem \ker (I)$ and assume without loss of generality that $i \in L \sem \ker (I)$ (the case $i \in R \sem \ker (I)$ is analogous). 
For $\varepsilon = \pm 1$, write $\varepsilon b_\ell =  ( c_i^\varepsilon , \varepsilon \dot b_{i \ell} ) \in \EE_i \star \dot \EE_{i \ell} '$ with $h_{\ell i} ( \varepsilon b_\ell) = c_i^\varepsilon$.
Clearly $c_i^1 + c_i^{-1} = \one$ but $c_i^1$ and $c_i^{-1}$ are compatible (otherwise $b_\ell$ would belong to $\EE_i ' = \EE_i \sub \EE_\ell$, contradiction). Let $G_i$ be
$\EE_i$-generic. (We may assume $c_i^1 \cdot c_i^{-1} \in G_i$  because outside of $c_i^1 \cdot c_i^{-1}$, $\dot b_{i\ell}$ is either $\zero$ or $\one$ and thus $\dot \EE_{i\ell}' =
\dot \EE_{i \ell}$ will be forced, though this is not really relevant). Let $a_{i\ell} \in \EE_{i \ell}$ be arbitrary. 
Following the proof of Lemma~\ref{secondextension}, in $V [G_i]$ define
\[ \tilde \mu_{i\ell} ( a_{i\ell} \cdot \varepsilon b_{i\ell} ) = \sup \{ \mu_{i\ell} ( a') : a' \in \EE_{i\ell} , a' \leq a_{i\ell} \cdot \varepsilon b_{i\ell} \} . \]

Next let $i < j < \ell$, $i,j \in \ker (I)$, and still work in $V[G_i]$. Decompose $a_{i\ell} \cdot \varepsilon b_{i\ell} \in \EE_{i\ell} '$ as $(c_{ij}^\varepsilon , \dot a_{j\ell} \cdot 
\varepsilon \dot b_{j\ell} ) \in \EE_{ij} \cdot \dot \EE_{j\ell}'$ where $c_{ij}^\varepsilon = h_{\ell j} ( a_{i\ell } \cdot \varepsilon b_{i \ell}) \in \EE_{ij}$. We first establish:

\begin{sclaim}    \label{amal-coherent-claim2}
\[  \tilde \mu_{i\ell} (a_{i\ell} \cdot \varepsilon b_{i\ell} ) \leq \int_{c_{ij}^\varepsilon } \dot{\tilde\mu}_{j\ell} (\dot a_{j\ell} \cdot \varepsilon \dot b_{j\ell}) d \mu_{ij} . \]
\end{sclaim}

\begin{proof}
Let $\delta > 0$, choose $a'_{i\ell} = (a_{ij} ' , \dot a_{j\ell} ' ) \in \EE_{i\ell} $ with $a'_{i\ell} \leq a_{i\ell} \cdot \varepsilon b_{i\ell}$ such that
$\tilde \mu_{i\ell} (a_{i\ell} \cdot \varepsilon b_{i\ell} )  < \mu_{i\ell} ( a'_{i\ell}) + \delta$. Thus $a_{ij} ' \leq c_{ij}^\varepsilon$.
Since $a'_{ij} \forces \dot a_{j\ell} ' \leq \dot a_{j\ell} \cdot \varepsilon \dot b_{j\ell}$,
we see $a'_{ij} \forces \dot \mu_{j\ell} ( \dot a_{j\ell} ')  \leq \dot{\tilde\mu}_{j\ell} (\dot a_{j\ell} \cdot \varepsilon \dot b_{j\ell})$. Therefore, using the coherence of $\bar \EE$, we obtain
\begin{eqnarray*}
\tilde \mu_{i\ell} (a_{i\ell} \cdot \varepsilon b_{i\ell} ) - \delta & < & \mu_{i\ell} ( a'_{i\ell}) = \int_{a'_{ij}} \dot \mu_{j\ell} ( \dot a_{j\ell} ')   d\mu_{ij} \leq 
   \int_{a'_{ij}} \dot{\tilde\mu}_{j\ell} (\dot a_{j\ell} \cdot \varepsilon \dot b_{j\ell})   d\mu_{ij} \\
   & \leq & \int_{c_{ij}^\varepsilon }  \dot{\tilde\mu}_{j\ell} (\dot a_{j\ell} \cdot \varepsilon \dot b_{j\ell})   d\mu_{ij} .
\end{eqnarray*}
Since this holds for every $\delta$, the claim is proved.
\end{proof}

Unfortunately equality need not hold here (see Counterexample~\ref{ctrex-product-measure} below right after the end of this proof). Therefore, by recursion on $\alpha < \omega_1$,
define
\begin{eqnarray*}
\mu^0_{i\ell} (a_{i\ell} \cdot \varepsilon b_{i\ell}) & =  & \tilde \mu_{i\ell} ( a_{i\ell} \cdot \varepsilon b_{i\ell})  \\
\mu^{\alpha + 1}_{i\ell} (a_{i\ell} \cdot \varepsilon b_{i\ell}) & =  & \sup_{i \leq j < \ell} \int_{c_{ij}^\varepsilon}  \dot \mu^\alpha_{j\ell} ( \dot a_{j\ell} \cdot \varepsilon \dot b_{j\ell})  d \mu_{ij} \\
\mu^\alpha_{i\ell} (a_{i\ell} \cdot \varepsilon b_{i\ell}) & =  & \lim_{\gamma < \alpha}  \mu_{i\ell}^\gamma ( a_{i\ell} \cdot \varepsilon b_{i\ell})  \;\;\;\;\; \mbox{ for limit } \alpha.
\end{eqnarray*}

\begin{sclaim}    \label{amal-coherent-claim3}
For all $\alpha < \omega_1$:
\begin{enumerate}
\item for $i \leq j < \ell$,  \[ \mu^\alpha_{i \ell} ( a_{i\ell} \cdot \varepsilon b_{i\ell})  \leq    \int_{c_{ij}^\varepsilon}  \dot \mu^\alpha_{j\ell} ( \dot a_{j\ell} \cdot \varepsilon \dot b_{j\ell})  d \mu_{ij} 
   \leq \mu^{\alpha + 1}_{i \ell} ( a_{i\ell} \cdot \varepsilon b_{i\ell}), \]
\item {\rm (monotonicity)}  for $i \leq j \leq j' < \ell$, \[ \int_{c_{ij}^\varepsilon}  \dot \mu^\alpha_{j\ell} ( \dot a_{j\ell} \cdot \varepsilon \dot b_{j\ell})  d \mu_{ij}  \leq
   \int_{c_{ij'}^\varepsilon}  \dot \mu^\alpha_{j' \ell} ( \dot a_{j' \ell} \cdot \varepsilon \dot b_{j' \ell})  d \mu_{ij'} . \]
\end{enumerate}
\end{sclaim}

\begin{proof}
By induction on $\alpha$. The second inequality in (i) is clear by definition for any $\alpha$. Monotonicity (ii) for any $\alpha$ is established by
\[    \int_{c_{ij}^\varepsilon}  \dot \mu^\alpha_{j\ell} ( \dot a_{j\ell} \cdot \varepsilon \dot b_{j\ell})  d \mu_{ij}  \leq
    \int_{c_{ij}^\varepsilon}  \int_{\dot c_{jj'}^\varepsilon}  \dot \mu^\alpha_{j' \ell} ( \dot a_{j' \ell} \cdot \varepsilon \dot b_{j' \ell})  d \dot \mu_{j j'} \; d \mu_{ij} =
   \int_{c_{ij'}^\varepsilon}  \dot \mu^\alpha_{j' \ell} ( \dot a_{j' \ell} \cdot \varepsilon \dot b_{j' \ell})  d \mu_{ij'} \]
where the inequality follows from (i) for $\alpha$ and the equality is coherence for integrals (Observation~\ref{coherentbasic} (iii)). The first
inequality of (i) for $\alpha = 0$ is Claim~\ref{amal-coherent-claim2}. For limit $\alpha$, it follows by induction hypothesis from
\[  \lim_{\gamma < \alpha}  \mu_{i\ell}^\gamma ( a_{i\ell} \cdot \varepsilon b_{i\ell})  = \lim_{\gamma < \alpha} \int_{c_{ij}^\varepsilon} 
    \dot \mu^\gamma_{j\ell} ( \dot a_{j\ell} \cdot \varepsilon \dot b_{j\ell})  d \mu_{ij}  \leq \int_{c_{ij}^\varepsilon}  \lim_{\gamma < \alpha} 
    \dot \mu^\gamma_{j\ell} ( \dot a_{j\ell} \cdot \varepsilon \dot b_{j\ell})  d \mu_{ij}. \]
So consider a successor $\alpha + 1$. Let $i \leq j < \ell$, and fix $\delta > 0$. Fix $j' \geq i$ such that 
\[ \mu_{i \ell}^{\alpha + 1} (a_{i \ell} \cdot \varepsilon b_{i \ell} ) <  \int_{c_{ij'}^\varepsilon}  \dot \mu^\alpha_{j'\ell} ( \dot a_{j'\ell} \cdot \varepsilon \dot b_{j'\ell})  d \mu_{ij'}  + \delta. \]
Recall that since $j , j' \in L$ they have a common upper bound $j \lor j' \in L$. By monotonicity for $\alpha$ (ii), we may thus assume $j' = j \lor j'$, i.e. $j' \geq j$.
By definition and coherence for integrals (Observation~\ref{coherentbasic} (iii)) we then see
\begin{eqnarray*}
\int_{c_{ij}^\varepsilon}  \dot \mu^{\alpha+1}_{j\ell} ( \dot a_{j\ell} \cdot \varepsilon \dot b_{j\ell})  d \mu_{ij} & \geq & 
   \int_{c_{ij}^\varepsilon}  \int_{c_{jj'}^\varepsilon}  \dot \mu^\alpha_{j'\ell} ( \dot a_{j'\ell} \cdot \varepsilon \dot b_{j'\ell})  d \mu_{j j'} \; d \mu_{ij} \\
   & = &   \int_{c_{ij'}^\varepsilon}  \dot \mu^\alpha_{j'\ell} ( \dot a_{j'\ell} \cdot \varepsilon \dot b_{j'\ell})  d \mu_{ij'}  > \mu^{\alpha + 1}_{i\ell} (a_{i\ell} \cdot \varepsilon b_{i\ell}) - \delta 
\end{eqnarray*}
Since this holds for every $\delta$ we obtain the inequality, and the proof of the claim is complete.
\end{proof}

So there is an ordinal $\alpha_0 = \alpha_0 (a_{i\ell}, \varepsilon ) < \omega_1$ such that for all $\alpha \geq \alpha_0$, we have 
$ \mu^\alpha_{i\ell} (a_{i\ell} \cdot \varepsilon b_{i\ell}) = \mu^{\alpha_0}_{i\ell} (a_{i\ell} \cdot \varepsilon b_{i\ell}) $. By the claim $c_{ij}^\varepsilon$
must force that $\dot \alpha_0 (\dot a_{j\ell} , \varepsilon) \leq \alpha_0$. Define
\[ \bar \mu_{i\ell} (a_{i\ell} \cdot \varepsilon b_{i \ell} ) =  \mu^{\alpha_0}_{i\ell} (a_{i\ell} \cdot \varepsilon b_{i\ell}). \]
Then by the claim
\[ \bar \mu_{i\ell} (a_{i\ell} \cdot \varepsilon b_{i \ell} ) = \int_{c_{ij}^\varepsilon} \dot{\bar \mu}_{j\ell} ( \dot a_{j\ell} \cdot \varepsilon b_{j \ell} ) d\mu_{ij} \tag{$\star$}  \]
for all $j$ with $i \leq j < \ell$.  

Next let 
\[ \mu_{i\ell} ' (a_{i\ell} \cdot \varepsilon b_{i\ell}) = {1 \over 2} \Big( \mu_{i\ell} (a_{i\ell}) +  \bar \mu_{i\ell} ( a_{i\ell} \cdot \varepsilon b_{i\ell} ) -  
   \bar \mu_{i\ell} ( a_{i\ell} \cdot (- \varepsilon b_{i\ell}) )  \Big) \]
and
\[ \mu_{i\ell} ' \left( a_{i\ell}^0 + a_{i\ell}^1 b + a_{i\ell}^2 (-b_{i\ell}) \right) = \mu_{i\ell} \left(a_{i\ell}^0\right) + \mu_{i\ell} ' \left(a_{i\ell}^1 b_{i\ell}  \right) + \mu_{i\ell} ' \left(a_{i\ell}^2 (-b_{i\ell})\right) \]
where $a_{i\ell}^0 + a_{i\ell}^1 b_{i\ell} + a_{i\ell}^2 (-b_{i\ell})$ is an arbitrary element of $\EE_{i\ell} '$ with the $a_{i\ell}^k \in \EE_{i\ell}$ being pairwise incompatible.
Note that if $a_{i\ell} \cdot \varepsilon b_{i \ell} \in \EE_{i \ell}$ then
\[  \mu_{i\ell} ' ( a_{i\ell} \cdot \varepsilon b_{i\ell} ) = \bar  \mu_{i\ell}  ( a_{i\ell} \cdot \varepsilon b_{i\ell} ) = \tilde  \mu_{i\ell} ( a_{i\ell} \cdot \varepsilon b_{i\ell} ) = 
    \mu_{i\ell}  ( a_{i\ell} \cdot \varepsilon b_{i\ell} ) \]
where the first and last equalities are obvious while the second follows from the coherence of the original system $\bar \EE$.  
Furthermore
\[  \mu_{i\ell} ' ( a_{i\ell} \cdot  b_{i\ell} )  + \mu_{i\ell} ' ( a_{i\ell} (- b_{i\ell} )) = \mu_{i \ell} (a_{i\ell}) \;\;\;   \mbox{ and } \;\;\; 
    \mu_{i\ell} ' ( a_{i\ell} \cdot \varepsilon b_{i\ell} ) +  \mu_{i\ell} ' ( \bar a_{i\ell} \cdot \varepsilon b_{i\ell} )  =  \mu_{i\ell} ' (( a_{i\ell} + \bar a_{i\ell}) \cdot \varepsilon b_{i\ell} ) \]
for incompatible $a_{i\ell} , \bar a_{i\ell}  \in \EE_{i\ell}$ where the second equality is first established for $\tilde \mu_{i\ell}$ and $\bar \mu_{i \ell}$. 
The argument of the proof of Lemma~\ref{secondextension} now shows that $\mu_{i\ell} '$ is a finitely additive strictly positive measure on $\EE_{i\ell} '$ extending the measure $\mu_{i\ell}$.

Finally, using $(\star)$ above and coherence of $\bar \EE$, we compute
\begin{eqnarray*}
\mu_{i\ell} ' ( a_{i\ell} \cdot \varepsilon b_{i\ell}    ) & = & {1 \over 2} \Big( \mu_{i\ell} 
   ( a_{i\ell}     ) + \bar \mu_{i \ell} (a_{i\ell} \cdot \varepsilon b_{i\ell}    ) - \bar \mu_{i \ell} (a_{i\ell} \cdot (- \varepsilon b_{i\ell} )    ) \Big) \\
   & = &  {1 \over 2} \left( \int_{a_{ij} } \dot \mu_{j\ell} 
   (  \dot a_{j\ell} ) d\mu_{ij}  +  \int_{ c_{ij}^\varepsilon } \dot{\bar \mu}_{j \ell} (    \dot a_{j\ell} \cdot \varepsilon \dot b_{j\ell} )  d\mu_{ij}
   - \int_{ c_{ij}^{-\varepsilon }  } \dot{\bar \mu}_{j \ell} (    \dot a_{j\ell} \cdot (-\varepsilon \dot b_{j\ell} )) d\mu_{ij}   \right)   \\
    & = &  {1 \over 2} \left( \int_{a_{ij} } \dot \mu_{j\ell} 
   (  \dot a_{j\ell} ) d\mu_{ij}  +  \int_{ a_{ij} } \dot{\bar \mu}_{j \ell} (    \dot a_{j\ell} \cdot \varepsilon \dot b_{j\ell} )  d\mu_{ij}
   - \int_{ a_{ij} }  \dot{\bar \mu}_{j \ell} (    \dot a_{j\ell} \cdot (-\varepsilon \dot b_{j\ell} )) d\mu_{ij}   \right)   \\
& = &  \int_{a_{ij}  } {1\over 2} \Big(      \dot \mu_{j\ell}  (  \dot a_{j\ell} )  +   \dot{\bar \mu}_{j \ell} (    \dot a_{j\ell} \cdot \varepsilon \dot b_{j\ell} )
   - \dot{\bar \mu}_{j \ell} (    \dot a_{j\ell} \cdot (-\varepsilon \dot b_{j\ell} ))  \Big) d\mu_{ij} \\
   & = & \int_{a_{ij} }   \dot \mu_{j\ell} ' ( \dot a_{j\ell} \cdot \varepsilon \dot b_{j\ell} )  d\mu_{ij} = \int_{c_{ij}^\varepsilon }   \dot \mu_{j\ell} ' ( \dot a_{j\ell} \cdot \varepsilon \dot b_{j\ell} )  d\mu_{ij}.
\end{eqnarray*}
This completes the proof of coherence of $\bar \EE'$ and of Main Lemma~\ref{amal-coherent}.
\end{proof}

The following explains why Claim~\ref{amal-coherent-claim2} is optimal.

\begin{counterexam}  \label{ctrex-product-measure}
Assume $\AA_0$ and $\AA_1$ are cBa's carrying finitely additive strictly positive measures $\mu_0$ and $\mu_1$, respectively. 
Let $\AA = \la \AA_0 \sem \{ \zero \} \times \AA_1 \sem \{ \zero \} \ra$, the algebra generated by the product, and let $\bar \AA = \ro ( \AA_0 \sem \{ \zero \} \times \AA_1 \sem \{ \zero \} )$,
the cBa generated by the product. Also assume there is a maximal antichain $(a_n : n \in \omega )$ in $\AA_0$ such that $\sum_n \mu_0 (a_n) \leq {1 \over 4}$ 
(this exists e.g. if $\AA_0$ is the Cohen algebra $\CC$), and that there are $(b_n : n \in \omega )$ with $\mu_1 (b_n) \geq {1 \over 2}$ such that whenever
$c' \leq \sum_n (a_n , b_n) \in \bar \AA$ belongs to $\AA$, then there is $n_0$ such that $c' \leq \sum_{n < n_0} (a_n , b_n)$ (this holds, e.g., if the $b_n$ are independent
elements of measure ${1 \over 2}$ in the random algebra $\AA_1 = \BB$). In particular $c:=\sum_n (a_n , b_n) \in \bar \AA \sem \AA$. Let $\mu$ be the (unique) measure
on $\AA$ induced by $\mu_0$ and $\mu_1$. Then $\tilde \mu (c) \leq {1 \over 4}$ where $\tilde \mu$ is defined from $\mu$ as in the proofs of Lemma~\ref{secondextension}
or Main Lemma~\ref{amal-coherent}. 

On the other hand, if we think of $\bar \AA$ as a two-step iteration first adding the $\AA_0$-generic and then the $\AA_1$-generic, then $c$ can be identified with a
condition of the form $(\one, \dot b)$ (where $a_n \forces \dot b = b_n$) and clearly $\int_\one \dot{\tilde \mu}_1 (\dot b) d \mu_0 \geq {1 \over 2}$. So $\tilde \mu (c)
< \int_\one \dot{\tilde \mu}_1 (\dot b) d \mu_0 $.

Furthermore, if we think of $\bar \AA$ as a two-step iteration the other way round, first adding the $\AA_1$-generic and then the $\AA_0$-generic, then $c$ is of the form
$(b, \dot a)$ where $b = \sum_n b_n$ and clearly $\int_b \dot{\tilde \mu}_0 (\dot a) d \mu_1 \leq {1 \over 4}$. In particular, if we adjust $\tilde \mu (c)$ to
$\bar \mu (c) = \int_\one \dot{\tilde \mu}_1 (\dot b) d \mu_0$ (as we did in the proof of Main Lemma~\ref{amal-coherent}), we will lose coherence on the other side.
\hfill $\dashv$
\end{counterexam}

This last argument explains why in the proof of Main Lemma~\ref{amal-coherent} we can only hope to preserve coherence for pairs $i,j \in I$ such that
$i \lor j$ exists in $I$, e.g. for $i,j$ both from $L \sem \ker (I)$ or for $i,j$ both from $R \sem \ker (I)$ (see the monotonicity argument in the proof of Claim~\ref{amal-coherent-claim3}),
and this is the reason why we consider quotient measures only for pure pairs.

\begin{mainlem}[coherence in iterations with $\BB$]    \label{random-coherent}
Assume $\bar \AA = \la \AA_i : i\in I\ra$ carries a coherent system $\la \dot \mu_{ij}^\AA :
(i,j) \in \pure (I) \ra$ of measures. Then $\bar \EE = \la \EE_i : i \in I \ra$ carries a coherent system 
$\la \dot \mu_{ij}^\EE : (i,j) \in \pure (I)  \ra$ of measures where we put
$\EE_i = \AA_i \star \dot \BB$. This is still true if $\dot \BB$ is replaced by $\dot \BB_\kappa$.
\end{mainlem}

\begin{proof}
Let $i < j$ be elements of $I$ with $(i,j) \in \pure (I)$. Assume we have conditions $q_i = (p_i , \dot b_i) \in \AA_i \star \dot \BB =\EE_i$ and
$q_j = (p_j , \dot b_j) \in \AA_j \star \dot \BB =\EE_j$ such that $h_{ji}^\EE (q_j) = q_i$,
in particular $h_{ji}^\AA (p_j) = p_i$. (As usual, we require here that 
$p_i \forces \mu (\dot b_i) > 0$ so that we have $h_{\EE_i \AA_i } (q_i) = p_i$. Similarly for $j$ instead of $i$.) We then also have
$\dot p_{ij} \in \dot \AA_{ij}$ and $\dot q_{ij} \in \dot\EE_{ij}$ such that
$p_j = (p_i , \dot p_{ij})$ and $q_j = (q_i , \dot q_{ij})$. We define 
\[ \dot\mu_{ij}^\EE \left(\dot q_{ij} \right) = 
\lim_{\dot c_i \leq \dot b_i, \dot \mu^\BB (\dot c_i) \to 0, \dot c_i \in \dot G_i^\BB}
\int_{\dot p_{ij}} { \dot \mu^\BB \left(\dot b_j \cap \dot c_i \right) \over \dot\mu^\BB \left(\dot c_i\right) }
d \dot \mu_{ij}^\AA \]
where $\dot \mu^\BB$ denotes (the name of) the standard measure on the random algebra.
That is, we take the limit along the $\dot\BB$-generic filter $\dot G_i^\BB$ in the
$\AA_i$-generic extension. Of course, this expression becomes clearer when working in this extension,
$p_i$ being in the generic, and we shall do so from now on.

We first need to argue this is \underline{\em well-defined}, that is, the limit exists.

\begin{sclaim}   \label{random-coherent-claim}
Given any $c_i \leq b_i$, $\epsilon > 0$, and $r \in [0,1]$,
there is $d_i \leq c_i$ such that \underline{either} for all $e_i \leq d_i$,
\[ \int_{p_{ij}} { \dot \mu^\BB \left(\dot b_j \cap e_i \right) \over \mu^\BB (e_i) } d\mu_{ij}^\AA \leq r 
  + \epsilon \]
\underline{or} for all $e_i \leq d_i$,
\[ \int_{p_{ij}} { \dot \mu^\BB \left(\dot b_j \cap e_i \right) \over \mu^\BB (e_i) } d\mu_{ij}^\AA \geq r 
  - \epsilon. \]
\end{sclaim}

\begin{proof}
Suppose this was false. Then there are $c_i \leq b_i$, $\epsilon > 0$, and
$r \in [0,1]$ such that for all $d_i \leq c_i$ there are $e_i , e_i ' \leq d_i$ such that
\[ \int_{p_{ij}} { \dot \mu^\BB \left(\dot b_j \cap e_i \right) \over \mu^\BB (e_i) } d\mu_{ij}^\AA > r 
  + \epsilon \]
and
\[ \int_{p_{ij}} { \dot \mu^\BB \left(\dot b_j \cap e_i '\right) \over \mu^\BB (e_i ') } d\mu_{ij}^\AA < r 
  - \epsilon. \]
So we get maximal antichains $\{ e_{i,n} :  n \in\omega\}$ and
$\{ e_{i,n} ' : n\in\omega\}$ below $c_i$ satisfying these inequalities. Thus we obtain
\begin{eqnarray*} \int_{p_{ij}}  \dot\mu^\BB \left(\dot b_j \cap c_i\right) d\mu_{ij}^\AA & = &
   \int_{p_{ij}}\left(  \sum_n \dot\mu^\BB \left(\dot b_j \cap e_{i,n}\right) \right) d\mu_{ij}^\AA \\
   & \stackrel{\textstyle (1)}{=} & \sum_n \left( \int_{p_{ij}} \dot\mu^\BB \left(\dot b_j \cap e_{i,n}\right) d\mu_{ij}^\AA \right) \\
   & > & \sum_n \mu^\BB (e_{i,n}) (r + \epsilon) = \mu^\BB (c_i) (r + \epsilon) 
\end{eqnarray*}
and
\[  \int_{p_{ij}}  \dot\mu^\BB \left(\dot b_j \cap c_i\right) d\mu_{ij}^\AA  = 
   \int_{p_{ij}} \sum_n \dot\mu^\BB \left(\dot b_j \cap e_{i,n} ' \right) d\mu_{ij}^\AA 
   < \mu^\BB (c_i) (r - \epsilon) , \]
a contradiction.

A comment concerning equality (1) about interchanging $\sum$ and $\int$ 
in the previous paragraph is in order, however. Of course, this is true when dealing
with $\sigma$-additive measures, but it fails for finitely additive measures in general, and
only the inequality $\geq$ can be proved\footnote{To see this note that if $a \in \AA$ and 
$\dot x_n$ are $\AA$-names for reals such that $a \forces ``0 \leq \dot x_n , \sum_n \dot x_n \leq 1"$, then
\[ \sum_{k < 2^m}  \left[ \mu \left(a \cdot \BOOL {k + 1 \over 2^m} \geq \sum_n \dot x_n > {k \over 2^m } \BOOR\right) \cdot {k \over 2^m} \right]
\geq \sum_n  \sum_{k < 2^m} \left[ \mu \left(a \cdot \BOOL {k+1 \over 2^m } \geq \dot x_n > {k \over 2^m} \BOOR \right)
\cdot {k \over 2^m } \right] \]
for all $m$.}. 
To see the reverse inequality $\leq$ also holds above,
proceed as follows.

Fix $m_0$ and $k_0 < 2^{m_0}$ and consider 
\[p = p_{ij} \cdot \BOOL { k_0 + 1 \over 2^{m_0} } \geq \sum_n
\dot\mu^\BB \left(\dot b_j \cap e_{i,n}\right) > {k_0 \over 2^{m_0}} \BOOR . \] 
It suffices to prove that
\begin{equation} \sum_n \int_p \dot\mu^\BB \left(\dot b_j \cap e_{i,n} \right) d\mu_{ij}^\AA \geq \mu^\AA_{ij} (p) \cdot {k_0 \over 2^{m_0}} . \end{equation}
Let $\epsilon > 0$. Then there is $n_0$ such that 
\begin{equation} p \FORCES_{ij} \sum_{n < n_0} \dot\mu^\BB \left(\dot b_j \cap e_{i,n}\right) > {k_0 \over 2^{m_0}} - {\epsilon \over 2} .  \end{equation}
This follows from the $\sigma$-additivity of the measure $\mu^\BB$ on $\BB$ and the fact
that the $e_{i,n}$ are not names. Choose $m_1 \geq m_0$ such that ${n_0 \over 2^{m_1}} < {\epsilon \over 2}$. We
shall argue that
\begin{equation} \sum_{n < n_0} \sum_{k < 2^{m_1}} \left[ \mu^\AA_{ij} \left(p \cdot \BOOL {k + 1 \over 2^{m_1}} \geq \dot\mu^\BB \left(\dot b_j \cap
e_{i,n}\right) > {k \over 2^{m_1}} \BOOR \right) \cdot {k \over 2^{m_1}} \right] \geq \left( \mu^\AA_{ij} (p) \cdot {k_0 \over 2^{m_0}} 
\right) - \epsilon \end{equation}
which implies (2) because $\epsilon$ is arbitrary. For $f : n_0 \to 2^{m_1}$ let
\begin{equation} p_f = p \cdot \prod_{n < n_0} \BOOL {f(n) + 1 \over 2^{m_1}} \geq \dot\mu^\BB \left(\dot b_j \cap e_{i,n}\right) >
{f(n) \over 2^{m_1}} \BOOR .  \end{equation}
Clearly 
\begin{equation} \sum_{n < n_0} \sum_{k < 2^{m_1}} \left[ \mu^\AA_{ij} \left(p \cdot \BOOL {k + 1 \over 2^{m_1}} \geq \dot\mu^\BB \left(\dot b_j \cap
e_{i,n}\right) > {k \over 2^{m_1}} \BOOR \right) \cdot {k \over 2^{m_1}} \right] = \sum_{f: n_0 \to 2^{m_1}} \left( \mu^\AA_{ij} (p_f) \sum_{n < n_0}
{f(n) \over 2^{m_1}} \right) . \end{equation}
However, by (3) and (5),
\[ p_f \FORCES_{ij} { \left(\sum_{n < n_0} f(n) \right) + n_0 \over 2^{m_1}} \geq \sum_{n < n_0} \dot\mu^\BB \left(\dot b_j \cap
e_{i,n} \right) > {k_0 \over 2^{m_0}} - {\epsilon \over 2}.  \]
Thus
\[ \sum_{n < n_0} {f(n) \over 2^{m_1}} > {k_0 \over 2^{m_0}} - {\epsilon \over 2} - {n_0 \over 2^{m_1}}
> {k_0 \over 2^{m_0}} - \epsilon \]
and
\[ \sum_{f: n_0 \to 2^{m_1}} \left( \mu^\AA_{ij} (p_f) \sum_{n<n_0} {f(n) \over 2^{m_1}} \right) > \sum_{f: n_0 \to 2^{m_1}} \mu^\AA_{ij} (p_f) \cdot
\left( {k_0 \over 2^{m_0}} - \epsilon\right) = \mu^\AA_{ij} (p) \left({k_0 \over 2^{m_0}} - \epsilon \right) \]
as required (see (4) and (6)).

This completes the proof of Claim~\ref{random-coherent-claim}.
\end{proof}

We continue with the proof of Main Lemma~\ref{random-coherent}.

Partitioning $[0,1]$ into finitely many small enough intervals and applying the claim repeatedly,
we see that for all $c_i \leq b_i$ and all $\epsilon > 0$ there are $d_i \leq c_i$ and
$r \in [0,1]$ such that for all $e_i \leq d_i$,
\[ \left|   \int_{p_{ij}}  {\dot\mu^\BB \left(\dot b_j \cap e_i \right) \over \mu^\BB (e_i) } d\mu^\AA_{ij} - r 
   \right| < \epsilon. \]
However, by density, this just means that the limit along the generic filter $\dot G_i^\BB$ exists.

Before proceeding with the proof, let us check that we have \underline{\em commutativity}
in the sense that
\[   \int_{p_{ij}} \dot\mu^\BB \left(\dot b_j \right) d\mu^\AA_{ij} = \int_{b_i} \dot\mu_{ij}^\EE \left(\dot q_{ij} \right)
   d\mu^\BB \]
where $(p_{ij} , \dot b_j) = (b_i , \dot q_{ij})$. Note that $\AA_{ij} \star \dot\BB
\cong \BB \star \dot \EE_{ij}$.
(We still work in the $\AA_i$-extension.)

By the above, given $\epsilon > 0$ and $c_i \leq b_i$, there is $d_i \leq c_i$ such that
\[   d_i \forces_\BB  \left|  \dot\mu_{ij}^\EE \left(\dot q_{ij}\right) - \int_{p_{ij}}   { \dot\mu^\BB \left(\dot b_j
   \cap d_i\right) \over \mu^\BB (d_i) } d\mu_{ij}^\AA \right| < \epsilon \]
which means that
\[   \left|   \int_{d_i}  \dot\mu_{ij}^\EE   \left(\dot q_{ij}\right) d\mu^\BB - \int_{p_{ij}}   \dot \mu^\BB \left(\dot b_j
   \cap d_i\right)  d\mu_{ij}^\AA   \right|   < \epsilon \cdot  \mu^\BB (d_i).  \]
Let $\{ d_{i,n} : n\in\omega\}$ be a maximal antichain of conditions below $b_i$ with this
property. Then
\begin{eqnarray*}   \left|  \int_{b_i}   \dot\mu^\EE_{ij}   \left(\dot q_{ij}\right) d\mu^\BB  \right.  -  \left. \int_{p_{ij}}  \dot \mu^\BB 
   \left(\dot b_j\right)  d\mu_{ij}^\AA  \right|   & = & \left| \sum_n \int_{d_{i,n}}  \dot \mu_{ij}^\EE
   \left(\dot q_{ij}\right) d\mu^\BB - \int_{p_{ij}} \left( \sum_n \dot\mu^\BB \left(\dot b_j \cap d_{i,n} \right) \right) d\mu^\AA_{ij} \right| \\
   &  \leq & \sum_n \left| \int_{d_{i,n}} \dot\mu^\EE_{ij}  \left(\dot q_{ij}\right) d\mu^\BB  -  \int_{p_{ij}}
   \dot\mu^\BB   \left(\dot b_j \cap d_{i,n} \right) d\mu^\AA_{ij} \right|   \\
   & < & \sum_n \epsilon \cdot \mu^\BB (d_{i,n}) = \epsilon \cdot \mu^\BB (b_i) \leq \epsilon
\end{eqnarray*}
where we use again the fact, established in the proof of the claim, that we can interchange $\sum$ and $\int$ in the right-hand term.
Since this is true for every $\epsilon$, commutativity follows.

It is easy to see that $\dot \mu^\EE_{ij} (\zero) = 0$ and $\dot \mu^\EE_{ij} (\one) = 1$ are forced, and 
we argue the measure $\dot \mu^\EE_{ij}$ is \underline{\em finitely additive}. (We still work in the $\AA_i$-extension and
use $\BB \star \dot \EE_{ij} \cong \AA_{ij} \star \dot \BB$ as before.)
Let $b_i \in \BB$ and $\dot q_{ij}^\ell \in \dot \EE_{ij}$, $\ell = 0,1$, such that
$b_i \forces_\BB ``\dot q_{ij}^0$ and $\dot q_{ij}^1$ are incompatible."
Let $p_{ij}^\ell \in \AA_{ij}$ and $\dot b_j^\ell \in \dot \BB$ such that
$(b_i,\dot q_{ij}^\ell) = (p_{ij}^\ell, \dot b_j^\ell)$.
Since the $(b_i,\dot q_{ij}^\ell)$ are incompatible, so are the $(p_{ij}^\ell, \dot b_j^\ell)$.
Therefore $p_{ij}^0 \cdot p_{ij}^1 \forces_{ij} \dot\mu^\BB (\dot b_j^0 \cap \dot b_j^1)
= 0$. Also we naturally have that $\one - p_{ij}^\ell \forces_{ij} \dot \mu^\BB (\dot b_j^\ell)
= 0$. Thus 
\[ p_{ij}^0 + p_{ij}^1 \forces_{ij} \dot\mu^\BB \left(\dot b^0_j \cup \dot b_j^1\right) =
\dot\mu^\BB \left(\dot b_j^0\right) + \dot\mu^\BB \left(\dot b_j^1\right) \]
and, clearly, $(b_i , \dot q_{ij}^0 + \dot q_{ij}^1) = ( p_{ij}^0 + p_{ij}^1 , \dot b^0_j \cup \dot b_j^1
)$. Let $p_{ij} = p_{ij}^0 + p_{ij}^1$ and $\dot b_j = \dot b_j^0 \cup \dot b_j^1$. Therefore
\begin{eqnarray*} \dot\mu_{ij}^\EE \left(\dot q_{ij}^0 + \dot q_{ij}^1 \right) & = & \lim_{c_i \leq
b_i , \mu^\BB (c_i) \to 0 , c_i \in \dot G^\BB_i } \int_{p_{ij}} { \dot \mu^\BB \left(\dot b_j \cap c_i\right)
\over \mu^\BB (c_i) } d\mu^\AA_{ij} \\
& = & \lim \int_{p_{ij}} { \dot \mu^\BB \left(\dot b_j^0 \cap c_i\right) + \dot \mu^\BB \left(\dot b_j^1 \cap c_i\right)
\over \mu^\BB (c_i) } d\mu^\AA_{ij} \\
& = & \lim \int_{p_{ij}^0} { \dot \mu^\BB \left(\dot b_j^0 \cap c_i\right)
\over \mu^\BB (c_i) } d\mu^\AA_{ij} +   \lim \int_{p_{ij}^1} { \dot \mu^\BB \left(\dot b_j^1 \cap c_i\right)
\over \mu^\BB (c_i) } d\mu^\AA_{ij} \\
& = & \dot\mu_{ij}^\EE \left(\dot q_{ij}^0\right) + \dot\mu_{ij}^\EE \left(\dot q_{ij}^1\right). 
\end{eqnarray*}
To see \underline{\em strict positivity}, let us take again a condition $(b_i , \dot q_{ij}) = (p_{ij} , \dot b_j) \in \BB \star \dot \EE_{ij} = \AA_{ij} \star \dot \BB$.
This representation means, in particular, that we assume $h^\EE_{ji} (p_{ij}, \dot b_j) = h^\EE_{ji}  ( b_i, \dot q_{ij} ) = b_i$
because $b_i \forces_\BB \dot q_{ij} > \zero$. Assume by way of contradiction that $c_i \forces_\BB \dot \mu^\EE_{ij} ( \dot q_{ij} ) = 0$  for some $c_i \leq b_i$. 
By commutativity,
\[ \int_{p_{ij}} \dot \mu^\BB \left( \dot b_j \cap c_i \right) d\mu_{ij}^\AA = \int_{c_i} \dot\mu_{ij}^\EE \left(\dot q_{ij} \right) d\mu^\BB = 0 \]
so that $p_{ij} \forces_{ij} \dot\mu^\BB ( \dot b_j \cap c_i) = 0$,  contradicting $h^\EE_{ji} (p_{ij}, \dot b_j) = b_i$.

Next we show \underline{\em coherence} (B). Assume $i < j <k$ are elements of $I$ with $(i,k) \in \pure (I)$. (We still work in the $\AA_i$-extension.) We have $b_i \forces_\BB \dot q_{ik}
= (\dot q_{ij} , \dot q_{jk}) \in \dot\EE_{ik} = \dot \EE_{ij} \star \dot\EE_{jk}$. There
are $p_{ik}$, $p_{ij}$ and $\dot b_k$, $\dot b_j$ such that $(p_{ik}, \dot b_k) = (b_i , \dot q_{ik})$ and
$(p_{ij} , \dot b_j) = (b_i , \dot q_{ij})$.
There is also $\dot p_{jk}$ such that $p_{ik} = (p_{ij} , \dot p_{jk})$. Now,
for any $c_i \leq b_i$, by commutativity,
\[ \int_{p_{ik}} \dot\mu^\BB \left(\dot b_k \cap c_i \right) d\mu_{ik}^\AA = \int_{c_i} \dot\mu_{ik}^\EE \left(\dot q_{ik}\right) d\mu^\BB \]
as well as
\begin{eqnarray*}
\int_{p_{ik}} \dot\mu^\BB \left(\dot b_k \cap c_i \right) d\mu_{ik}^\AA & = & 
\int_{p_{ij}} \int_{\dot p_{jk}} \dot\mu^\BB \left(\dot b_k \cap c_i \right) d\dot\mu_{jk}^\AA d\mu_{ij}^\AA \\
&=& \int_{p_{ij}} \int_{\dot b_j \cap c_i} \dot\mu^\EE_{jk} \left(\dot q_{jk} \right) d\dot\mu^\BB d\mu_{ij}^\AA =
\int_{c_i} \int_{\dot q_{ij}} \dot\mu_{jk}^\EE \left(\dot q_{jk}\right) d\dot\mu^\EE_{ij} d\mu^\BB
\end{eqnarray*}
where the first equality uses coherence of the system $\bar \AA$, as applied to integrals (Observation~\ref{coherentbasic} (iii)) while the other two equalities are again by commutativity.
Since this holds for any $c_i \leq b_i$, we get 
\[ b_i \forces_\BB \dot\mu_{ik}^\EE (\dot q_{ik}) = \int_{ \dot q_{ij}} \dot\mu_{jk}^\EE (\dot q_{jk}) d\dot\mu_{ij}^\EE \]
as required. 

Finally we show \underline{\em correctness} (C). Let $i, j \in I$ be incomparable such that $i \lor j $ exists, and work in the $\AA_{i \land j}$-extension. 
Let $b_{i \land j} \in \BB$ and let $\dot q_{i \land j, i}$ be a $\BB$-name for a condition in $\dot \EE_{i \land j, i}$. 
Put $b_j = b_{i \land j}$ and think of this as a condition in $\AA_{i \land j, j} \star \dot \BB = \BB \star \dot \EE_{i \land j, j}$,
and identify $\dot q_{i \land j, i}$ with an $\AA_{i \land j, j} \star \dot \BB$-name $\dot q_{j, i\lor j} \in \Dot \EE_{j, i \lor j}$.
There are a condition $p_{i \land j, i} = p_{j, i \lor j} \in \AA_{i \land j, i}$ and an $\AA_{i \land j, i}$-name $\dot b_i = \dot b_{i \lor j}$
such that $( p_{i \land j, i} , \dot b_i) = (b_{i\land j} , \dot q_{i \land j, i})$ and $( p_{j, i \lor j} , \dot b_{i\lor j } ) = (b_{ j} , \dot q_{j, i \lor j})$.
For any $c_{i \land j} \leq b_{i \land j}$, by commutativity,
\[ \int_{ p_{ i \land j, i}} \dot\mu^\BB \left( \dot b_i \cap c_{i \land j} \right) d \mu_{i \land j, i}^\AA = \int_{c_{i\land j}} 
   \dot \mu^\EE_{i \land j, i} \left(\dot q_{i \land j, i }\right)  d \mu^\BB \]
and
\[ \int_{ p_{ j , i \lor j}} \dot\mu^\BB \left( \dot b_{i \lor j} \cap c_{i \land j} \right) d \mu_{j, i \lor j}^\AA = \int_{c_{i\land j}} 
   \dot \mu^\EE_{j, i \lor j}  \left(\dot q_{j, i \lor j}\right)  d \mu^\BB . \]
However, the two left-hand integrals are the same by correctness for integrals (Observation~\ref{coherentbasic} (iv)) of the system $\bar \AA$ and, thus, so are the two 
right-hand integrals. Since this holds for any $c_{i \land j} \leq b_{i \land j}$, we get
\[ b_{i \land j} \forces_\BB \dot \mu^\EE_{j , i \lor j} (\dot q_{i , i \lor j} ) = \dot \mu^\EE_{i \land j, i} (  \dot q_{i \land j, i } ) \]
as required.

This completes the proof of Main Lemma~\ref{random-coherent}.
\end{proof}



\section{The shattered iteration of Cohen forcing}   \label{Cohen}

In this section we prove Main Theorem A.
Let $\kappa < \lambda$ be uncountable regular cardinals. We think of them
as being disjoint and use  $\kappa \dot\cup \lambda$ to denote this disjoint union.
For ordinals $\zeta < \xi$, $[\zeta,\xi) = \{ \zeta' : \zeta \leq \zeta' < \xi \}$ denotes the
corresponding interval of ordinals.

First some heuristics which motivate the use of the amalgamation (and the amalgamated limit)
in our context. To guarantee $\Cov N = \lambda$ we plan to add $\lambda$ random reals, adjoined by
$\BB_{\kappa \dot\cup \lambda}$ for technical reasons, and this may be thought of as the basic step of our iteration.
To increase $\Cov M$ to $\kappa$ we will then  start adjoining Cohen reals $c_{\zeta,\eta}$ for any pair $(\zeta,\eta) \in
\kappa \times \lambda$, and this will be done in such a fashion that $c_{\zeta,\eta}$ added by $\dot \CC_{\zeta,\eta}$
is Cohen over $V^{\BB_{\zeta\dot\cup\eta}}$, that is, over an initial segment of the randoms, while the remaining
randoms, namely, those adjoined by $\BB_{[\zeta,\kappa) \dot\cup [\eta,\lambda)}$ are still random over 
$V^{\BB_{\zeta\dot\cup\eta} \star \dot \CC_{\zeta,\eta}}$. Notice that sets of the form $\zeta \dot\cup\eta$ and $\zeta' \dot\cup\eta'$
are in general not contained one in the other  so that the basic building
block we shall encounter is this: we have random reals $b_0 , b_1$ and Cohen reals $c_0, c_1$
such that $c_i$ is Cohen over $b_i$ while $b_i$ is random over $c_{1-i}$. Writing $\BB_i$
and $\CC_i$ for the corresponding forcing notions we get the diagram
\[
\begin{picture}(180,60)(0,0)
\put(5,5){\makebox(0,0){$\CC_0$}}
\put(175,5){\makebox(0,0){$\CC_1$}}
\put(95,5){\makebox(0,0){$\BB_0 \star \dot \BB_1 \cong \BB_1 \star \dot \BB_0$}}
\put(50,30){\makebox(0,0){$(\BB_0 \times \CC_0) \star \dot \BB_1$}}
\put(135,30){\makebox(0,0){$(\BB_1 \times \CC_1) \star \dot \BB_0$}}
\put(95,55){\makebox(0,0){$\AA$}}
\put(55,37){\line(2,1){25}}
\put(55,20){\line(2,-1){20}}
\put(110,10){\line(2,1){20}}
\put(108,50){\line(2,-1){25}}
\put(10,10){\line(2,1){20}}
\put(148,20){\line(2,-1){20}}
\end{picture}
\]
in which the amalgamation $\AA = \amal_{\BB_0 \star \dot \BB_1} ((\BB_0 \times \CC_0) \star
\dot \BB_1 , (\BB_1 \times \CC_1) \star \dot \BB_0 )$ arises naturally. To appreciate this diagram, recall that
being random is commutative, i.e., if $b_1$ is random over $V[b_0]$, then
$b_0$ is random over $V[b_1]$~\cite{Ku84}, and that forcing with an iteration $\PP \star \dot \CC$
is equivalent to forcing with the product $\PP \times \CC$ since Cohen forcing is forcing
with finite partial functions.

The product $\kappa \times \lambda$ is naturally partially ordered by $(\zeta , \eta) \leq (\xi,\theta)$ if
$\zeta \leq \xi$ and $\eta \leq \theta$. Thus $ (\zeta,\eta) < (\xi,\theta)$ means that $\zeta \leq \xi$,
$\eta \leq \theta$, and either $\zeta < \xi$ or $\eta < \theta$. Note that this partial order has meets and joins
given by $(\zeta,\eta) \land (\xi,\theta) = (\zeta \cap \xi, \eta \cap \theta)$ and $(\zeta,\eta) \lor (\xi,\theta) 
= (\zeta \cup \xi, \eta \cup \theta)$. Also note that $\leq$ is wellfounded on
$\kappa \times \lambda$. This is crucial because we use this ordering for the recursive definition of an iteration.
Say a pair $((\zeta,\eta),(\xi,\theta))$ with $(\zeta,\eta) \leq (\xi,\theta)$ is {\em pure} if either $\zeta = \xi$ or $\eta = \theta$.
For $(\xi,\theta) \in \kappa \times \lambda$ let
\begin{eqnarray*}
I_{\xi,\theta} & = & \{ (\zeta,\eta) : (\zeta,\eta) < (\xi,\theta) \}, \\
J_{\xi,\theta} & = & \{ (\zeta,\eta) : (\zeta,\eta) \leq (\xi,\theta) \}.
\end{eqnarray*}

\begin{obs}   \label{ordinal-lattice}
\begin{enumerate}
\item All $J_{\xi,\theta}$ as well as the $I_{\xi,\theta}$ with $\xi = 0$ or $\theta = 0$ are lattices. Furthermore,
   $\pure (J_{\xi,\theta}) = \{ (( \zeta ,\eta),(\zeta',\eta')):  (\zeta,\eta), (\zeta',\eta') \in J_{\xi,\theta}$ and the pair is pure$\}$ 
   for $\xi >0$ and $\theta >0$.
\item For $\xi > 0$ and $\theta > 0$, $I_{\xi,\theta}$ is a nice distributive almost-lattice with $\ker (I_{\xi,\theta} ) =
   \{ (\zeta,\eta) : \zeta < \xi$ and $\eta < \theta \}$, $L (I_{\xi,\theta} ) =
   \{ (\zeta,\eta) : \zeta < \xi$ and $\eta \leq \theta \}$, and $R (I_{\xi,\theta} ) =
   \{ (\zeta,\eta) : \zeta \leq  \xi$ and $\eta < \theta \}$. Furthermore, $\pure (I_{\xi,\theta}) = \{ (( \zeta ,\eta),(\zeta',\eta')):  (\zeta,\eta), (\zeta',\eta') \in I_{\xi,\theta}$ and the pair is pure$\}$.
\end{enumerate}
\end{obs}

\begin{proof}
(i) It is easy to see that the sets are lattices. To see any pure pair from $J_{\xi,\theta}$ belongs to $\pure (J_{\xi,\theta})$, notice this is clear
for pairs of the form $((\zeta,\eta), (\zeta,\eta))$ by (2)(a) in Definition~\ref{defin3a}. Assume $((\zeta,\eta),(\zeta',\eta))$ is pure with $\zeta < \zeta '$. If $\eta = 0$,
notice that $(\zeta,\theta)$ and $(\zeta ', 0)$ are incomparable and $(\zeta,0) = (\zeta,\theta) \land (\zeta',0)$. If $\eta > 0$, notice that
$(\zeta',0)$ and $(\zeta,\eta)$ are incomparable and $(\zeta',\eta) = (\zeta',0) \lor (\zeta,\eta)$. So $((\zeta,\eta),(\zeta',\eta)) \in \pure (J_{\xi,\theta})$
by (2)(b). Similarly for pure pairs of the form $((\zeta,\eta),(\zeta,\eta '))$. On the other hand, if $((\zeta,\eta), (\zeta',\eta')) \in \pure (J_{\xi,\theta})$,
according to (2)(a), (b), or (c), we must have either $\zeta = \zeta'$ or $\eta = \eta '$ or both.

(ii) Clearly $I_{\xi,\theta}$ is a distributive almost-lattice, and the characterizations of  $\ker$, $L$, and $R$ are also obvious.
To see niceness, let $(\zeta,\eta) < (\zeta',\eta')$ both from $I_{\xi,\theta} \sem \ker (I_{\xi,\theta})$, without loss of generality $(\zeta,\eta) , (\zeta',\eta') \in
L (I_{\xi,\theta}) \sem \ker (I_{\xi,\theta})$. Thus $\eta = \eta ' = \theta$ and $\zeta < \zeta '< \xi$. We then see that $(\zeta', 0)$ is incomparable with
$(\zeta,\theta)$ and $(\zeta,\theta) \lor (\zeta' , 0) = (\zeta ', \theta)$. The argument for $\pure (I_{\xi,\theta})$ is exactly like the one for $\pure (J_{\xi,\theta})$.
\end{proof}

We remark that for the $I_{0,\theta}$ and the $I_{\xi,0}$, the set of pure pairs is trivial, and similarly for the corresponding $J$'s.

\begin{defin} \label{defin5}
By recursion on $(\zeta,\eta) \in \kappa \times \lambda$ we define the {\em shattered iteration} 
$\left\la \AA_{\zeta,\eta}^{\xi,\theta} : \zeta \leq \xi < \kappa, \eta \leq \theta < \lambda \right\ra$ of the
forcing notions $\left\la \dot \EE_{\zeta, \eta} : \zeta < \kappa, \eta < \lambda \right\ra$ as follows:
\begin{enumerate}
\item $\EE_{0,0}$ is a cBa in the ground model and $\AA_{0,0}^{\xi,\theta} = \EE_{0,0} \star \dot\BB_{ \xi \dot\cup \theta}$.
\item For $(\zeta,\eta) > (0,0)$ we first let $\AA_{\amal}^{\zeta,\eta} = \lim\amal_{(\zeta_0, \eta_0) < (\zeta,\eta)} \AA^{\zeta,\eta}_{\zeta_0,\eta_0}$;
   $\dot \EE_{\zeta,\eta}$ is an $\AA_{\amal}^{\zeta,\eta}$-name for a cBa and we let $\AA^{\xi,\theta}_{\zeta,\eta} = \AA_{\amal}^{\zeta,\eta}
   \star \dot \EE_{\zeta,\eta} \star \dot \BB_{[\zeta,\xi) \dot\cup [\eta,\theta)}$.
\end{enumerate}
$\AA = \AA_{\lim} = \lim \dir_{(\zeta,\eta) \in \kappa \times \lambda} \AA^{\zeta,\eta}_{\zeta,\eta}$ is the direct limit of the shattered iteration.
\hfill $\dashv$
\end{defin}

We first need to see that this definition makes sense. For indeed, for the definition of the amalgamated limit we need
that the corresponding systems of cBa's have complete embeddings and correct projections. Therefore, the above
definition goes by recursion-induction, with the simultaneous inductive proof of the following lemma.

\begin{lem}   \label{shattered-embed}
\begin{enumerate}
\item For any $\zeta_0 \leq \zeta \leq \xi$, $\zeta_0 \leq \xi_0 \leq \xi$, $\eta_0 \leq \eta \leq \theta$, and $\eta_0 \leq \theta_0 \leq \theta$,
   we have $\AA_{\zeta_0,\eta_0}^{\xi_0,\eta_0} \embed \AA_{\zeta,\eta}^{\xi,\theta}$.
\item The systems $\left\la \AA_{\zeta_0,\eta_0}^{\zeta,\eta} : (\zeta_0,\eta_0 ) < (\zeta,\eta) \right\ra$ and $\left\la \AA_{\zeta_0,\eta_0}^{\zeta,\eta} : (\zeta_0,\eta_0 ) \leq (\zeta,\eta) \right\ra$
   have correct projections.
\end{enumerate}
\end{lem}

\begin{proof}
(i) and (ii) are proved by simultaneous induction on $(\zeta,\eta) \in \kappa\times \lambda$. Each stage $(\zeta,\eta)$ here corresponds 
to stage $(\zeta,\eta)$ in Definition~\ref{defin5}.

\underline{$(\zeta,\eta) = (0,0)$}. Nothing to prove.

\underline{$(\zeta,\eta) > (0,0)$}. By induction hypothesis we assume all instances of (i) and (ii) hold for any $(\zeta', \eta') < (\zeta,\eta)$.
We first prove the first half of (ii). Note that by induction hypothesis for (i), the system $\left\la \AA_{\zeta_0,\eta_0}^{\zeta,\eta} : (\zeta_0,\eta_0 ) < (\zeta,\eta) \right\ra$
has complete embeddings. Now let $(\zeta_0,\eta_0), (\zeta_1,\eta_1)$ be such that $(\zeta_0 \cup \zeta_1, \eta_0 \cup \eta_1) < (\zeta,\eta)$ and consider the
diagram
\[
\begin{picture}(50,60)(0,0)
\put(28,3){\makebox(0,0){$\AA_{\zeta_0 \cap \zeta_1, \eta_0 \cap \eta_1}^{\zeta,\eta}$}}
\put(-7,28){\makebox(0,0){$ \AA_{\zeta_0,\eta_0}^{\zeta,\eta} $}}
\put(58,28){\makebox(0,0){$\AA_{\zeta_1,\eta_1}^{\zeta,\eta}$}}
\put(28,58){\makebox(0,0){$\AA_{\zeta_0 \cup \zeta_1, \eta_0 \cup \eta_1}^{\zeta,\eta}$}}
\put(6,36){\line(1,1){13}}
\put(8,23){\line(1,-1){11}}
\put(31,11){\line(1,1){13}}
\put(33,48){\line(1,-1){11}}
\end{picture}
\]
If we replace the superscript $(\zeta,\eta)$ by $(\zeta_0 \cup \zeta_1, \eta_0 \cup \eta_1) $, then we have correct projections by induction hypothesis for (ii).
Note, however, that for any subscript $(\zeta',\eta') \leq ( \zeta_0 \cup \zeta_1, \eta_0 \cup \eta_1)$, $\AA_{\zeta',\eta'}^{\zeta,\eta} = \AA_{\zeta',\eta'}^{\zeta_0 \cup \zeta_1, \eta_0 \cup \eta_1} \star
\dot \BB_{[ \zeta_0 \cup \zeta_1, \zeta) \dot\cup [  \eta_0 \cup \eta_1, \eta) }$. Thus the correctness of the above diagram follows from Lemma~\ref{random-correct}.

Using this, $\AA^{\zeta,\eta}_{\amal}$ and all $\AA^{\xi,\theta}_{\zeta,\eta}$ can indeed be defined according to Definition~\ref{defin5}.

We next prove (i). First notice that since $\AA^{\zeta,\eta}_{\zeta,\eta} =  \lim\amal_{(\zeta_0, \eta_0) < (\zeta,\eta)} \AA^{\zeta,\eta}_{\zeta_0,\eta_0}
 \star \dot \EE_{\zeta,\eta}$, by Main Lemma~\ref{amal-embed}, 
\begin{equation} 
\AA^{\zeta,\eta}_{\zeta_0,\eta_0} \embed \AA^{\zeta,\eta}_{\zeta,\eta} 
\end{equation}
holds for all $(\zeta_0,\eta_0) < (\zeta, \eta)$. For the arbitrary case notice that 
\begin{eqnarray*}
\AA_{\zeta_0 , \eta_0}^{\xi_0 , \theta_0} 
& \stackrel{\rm (def)}{=} & \AA_{\zeta_0 , \eta_0}^{\zeta_0 , \eta_0} \star \dot \BB_{[\zeta_0,\xi_0) \dot\cup
   [\eta_0,\theta_0)} \embed
\AA_{\zeta_0 , \eta_0}^{\zeta_0 , \eta_0} \star \dot \BB_{[\zeta_0,\xi) \dot\cup
   [\eta_0,\theta)}\\
& = &
\left(\AA_{\zeta_0 , \eta_0}^{\zeta_0 , \eta_0} \star \dot \BB_{[\zeta_0,\zeta) \dot\cup
   [\eta_0,\eta)} \right) \star \dot \BB_{[\zeta,\xi) \dot\cup
   [\eta,\theta)}   \stackrel{\rm (def)}{=}
\AA_{\zeta_0 , \eta_0}^{\zeta , \eta} \star \dot \BB_{[\zeta,\xi) \dot\cup
   [\eta,\theta)}  \\
& \stackrel{\rm (7)}{\embed} &
\AA_{\zeta , \eta}^{\zeta , \eta} \star \dot \BB_{[\zeta,\xi) \dot\cup
   [\eta,\theta)} \stackrel{\rm (def)}{=}
\AA_{\zeta , \eta}^{\xi , \theta},
\end{eqnarray*} 
as required.

Finally note that the second half of (ii) now follows from Lemma~\ref{amal-correct} and part (iii) of
Observation~\ref{correct-basic}.
\end{proof}

\begin{lem}   \label{shattered-measures}
Assume that, in Definition~\ref{defin5}, for all $(\zeta,\eta) \in \kappa \times \lambda$, $\AA_{\amal}^{\zeta,\eta}$ forces that
$\dot \EE_{\zeta,\eta}$ carries a finitely additive strictly positive measure (where  $\AA_{\amal}^{0,0}$ is the trivial forcing).
Then  the systems $\left\la \AA_{\zeta_0,\eta_0}^{\xi,\theta} : (\zeta_0,\eta_0 ) < (\zeta,\eta) \right\ra$ and 
$\left\la \AA_{\zeta_0,\eta_0}^{\xi,\theta} : (\zeta_0,\eta_0 ) \leq (\zeta,\eta) \right\ra$ carry coherent systems of measures
$\left\la \dot \mu^{\xi,\theta}_{(\zeta_0,\eta_0),(\zeta_1,\eta_1)} : ((\zeta_0,\eta_0) ,(\zeta_1, \eta_1) ) \in \pure (I_{\zeta,\eta}) \right\ra$
and $\left\la \dot \mu^{\xi,\theta}_{(\zeta_0,\eta_0),(\zeta_1,\eta_1)} :  ((\zeta_0,\eta_0) ,(\zeta_1, \eta_1) ) \in \pure (J_{\zeta,\eta}) \right\ra$
for any $(\xi,\theta) \geq (\zeta,\eta)$. Furthermore all $\AA_{\zeta,\eta}^{\xi,\theta}$ carry finitely additive strictly positive measures.
\end{lem}

\begin{proof}
By induction on $(\zeta,\eta)$.

\underline{$(\zeta,\eta) = (0,0)$}. Since $\EE_{0,0}$ carries a finitely additive strictly positive measure, so does the two-step iteration
$\AA_{0,0}^{\xi,\theta} = \EE_{0,0} \star \dot \BB_{\xi \dot\cup \theta}$.

\underline{$(\zeta,0) $ and $(0,\eta)$, where $\zeta$, $\eta > 0$}. By the remark after the proof of Observation~\ref{ordinal-lattice}, the set of pure pairs is trivial.
Hence we only need to show the $\AA^{\xi,\theta}_{\zeta,0}$ and $\AA^{\xi,\theta}_{0,\eta}$ carry finitely additive strictly positive measures. This is done
by induction on $\zeta$ ($\eta$, respectively). $\zeta = 0$ is the previous case. If $\zeta = \zeta_0 + 1$ is successor then $\AA^{\xi,\theta}_{\zeta,0} =
\AA^{\zeta,0}_{\zeta_0,0} \star \dot \EE_{\zeta,0} \star \dot \BB_{[\zeta,\xi) \dot\cup \theta}$ and if $\zeta$ is limit then $\AA^{\xi,\theta}_{\zeta,0} = \left(\lim\dir_{\zeta_0 < \zeta}
\AA^{\zeta,0}_{\zeta_0, 0} \right) \star \dot \EE_{\zeta,0} \star \dot \BB_{[\zeta,\xi) \dot\cup \theta}$ and we obtain the measure from the induction hypothesis, the measures on
$\dot \EE_{\zeta, 0}$ and $\dot \BB_{[\zeta,\xi) \dot\cup \theta}$, and the fact that the direct limit of cBa's with a measure carries a measure \cite[2.6]{Ka89}.

\underline{$(\zeta,\eta) = (1,1)$}. The non-trivial pure pairs in $I_{1,1}$ are $((0,0), (0,1))$ and $((0,0), (1,0))$. Note that $\AA^{1,0}_{1,0} = \AA^{1,0}_{0,0} \star
\dot \EE_{1,0}$ so that $\dot \AA^{1,0}_{ (0,0), (1,0)} = \dot \EE_{1,0}$ which carries a measure $\dot\mu^{1,0}_{(0,0),(1,0)}$. Similarly for $\dot\mu^{0,1}_{(0,0),(0,1)}$.
Main Lemma~\ref{random-coherent} gives us measures $\dot\mu^{1,1}_{(0,0),(1,0)}$ and $\dot\mu^{1,1}_{(0,0),(0,1)}$. From Main Lemma~\ref{amal-coherent} we obtain
measures $\dot \mu^{1,1}_{(0,1), \amal}$ and $\dot \mu^{1,1}_{(1,0), \amal}$ satisfying correctness (C) and then measures $\dot \mu^{1,1}_{(0,1),(1,1)}$
and $\dot \mu^{1,1}_{(1,0),(1,1)}$ on $\dot \AA^{1,1}_{(0,1) , (1,1) }= \dot \AA^{1,1}_{(0,1), \amal} \star \dot \EE_{1,1}$ and 
$\dot \AA^{1,1}_{(1,0) , (1,1) }= \dot \AA^{1,1}_{(1,0), \amal} \star \dot \EE_{1,1}$ so that correctness (C) is trivially preserved. 
Finally, another application of Main Lemma~\ref{random-coherent} gives us measures $\dot \mu^{\xi,\theta}_{(0,1),(1,1)}$
and $\dot \mu^{\xi,\theta}_{(1,0),(1,1)}$ satisfying correctness (C). We obtain a measure on $\AA^{\xi,\theta}_{1,1}$ for example by iterating
a measure on $\AA^{\xi,\theta}_{1,0}$ with $\dot \mu^{\xi,\theta}_{(1,0),(1,1)}$.

\underline{$(\zeta,\eta) > (1,1)$}. By induction hypothesis, we already have the system 
$\left\la \dot \mu^{\xi,\theta}_{(\zeta_0,\eta_0),(\zeta_1,\eta_1)} : ((\zeta_0,\eta_0) ,(\zeta_1, \eta_1) ) \in \pure (I_{\zeta,\eta}) \right\ra$ of measures and,
since coherence (B) and correctness (C) are a local properties, this system is clearly coherent and corrrect. Next, for $(\zeta_0,\eta_0) < (\zeta,\eta)$ with $\zeta_0 = \zeta$ or $\eta_0 = \eta$
we obtain the measures $\dot \mu^{\zeta,\eta}_{(\zeta_0,\eta_0),\amal}$ on the quotient algebras $\dot \AA^{\zeta,\eta}_{(\zeta_0,\eta_0), \amal}$
(satisfying $\AA^{\zeta,\eta}_\amal = \AA^{\zeta,\eta}_{\zeta_0,\eta_0} \star \dot \AA^{\zeta,\eta}_{(\zeta_0,\eta_0), \amal}$, see Lemma~\ref{amal-remainder})
from Main Lemma~\ref{amal-coherent} so that coherence (B) and correctness (C) are preserved. Next, we have the measure $\dot \mu^\EE_{\zeta,\eta}$ on
$\dot \EE_{\zeta,\eta}$. Now, for a name for a condition $(\dot p, \dot q)$ in the two-step iteration $\dot \AA^{\zeta,\eta}_{(\zeta_0,\eta_0), \amal}
\star \dot \EE_{\zeta,\eta}$, $\zeta_0 = \zeta$ or $\eta_0 = \eta$, simply define 
\[ \dot \mu^{\zeta,\eta}_{(\zeta_0,\eta_0),(\zeta,\eta)} (\dot p , \dot q)  = \int_{\dot p} \dot \mu^\EE_{\zeta,\eta} ( \dot q) \; d\dot \mu^{\zeta,\eta}_{(\zeta_0,\eta_0),\amal}. \]
This is clearly a finitely additive strictly positive measure, and to see coherence (B) note that for $(\zeta_0,\eta) < (\zeta_1,\eta) < (\zeta,\eta)$, 
writing now $\dot p_{0, \amal} = ( \dot p_{0,1} , \dot p_{1, \amal} ) \in \dot \AA^{\zeta,\eta}_{(\zeta_0,\eta), \amal} = \dot \AA^{\zeta,\eta}_{(\zeta_0,\eta),(\zeta_1,\eta)}
\star \dot \AA^{\zeta,\eta}_{(\zeta_1,\eta), \amal}$, we have
\begin{eqnarray*} 
   \dot \mu^{\zeta,\eta}_{(\zeta_0,\eta),(\zeta,\eta)} (\dot p_{0,\amal} , \dot q) & = &\int_{\dot p_{0,\amal}} \dot \mu^\EE_{\zeta,\eta} ( \dot q) \; d\dot \mu^{\zeta,\eta}_{(\zeta_0,\eta),\amal}
   = \int_{\dot p_{0,1}} \left( \int_{\dot p_{1,\amal}} \dot \mu^\EE_{\zeta,\eta} ( \dot q) \; d\dot \mu^{\zeta,\eta}_{(\zeta_1,\eta),\amal} \right) d\dot \mu^{\zeta,\eta}_{(\zeta_0,\eta), (\zeta_1,\eta)} \\
   &= & \int_{\dot p_{0,1}}     \dot \mu^{\zeta,\eta}_{(\zeta_1,\eta),(\zeta,\eta)} (\dot p_{1,\amal} , \dot q)      \; d\dot \mu^{\zeta,\eta}_{(\zeta_0,\eta), (\zeta_1,\eta)},
\end{eqnarray*}
by coherence for integrals (Observation~\ref{coherentbasic} (iii)) as required. The case $(\zeta,\eta_0) < (\zeta,\eta_1) < (\zeta,\eta)$ is analogous. Preservation of correctness (C) is
trivial. Thus we have the coherent system $\left\la \dot \mu^{\zeta,\eta}_{(\zeta_0,\eta_0),(\zeta_1,\eta_1)} :  ((\zeta_0,\eta_0) ,(\zeta_1, \eta_1) ) \in \pure (J_{\zeta,\eta}) \right\ra$, 
and the coherent system $\left\la \dot \mu^{\xi,\theta}_{(\zeta_0,\eta_0),(\zeta_1,\eta_1)} :  ((\zeta_0,\eta_0) ,(\zeta_1, \eta_1) ) \in \pure (J_{\zeta,\eta}) \right\ra$ for any
$(\xi,\theta) > (\zeta,\eta)$ is now obtained by an application of Main Lemma~\ref{random-coherent}. We get a measure on $\AA^{\xi,\theta}_{\zeta,\eta}$
by iterating a measure on $\AA^{\xi,\theta}_{\zeta,0}$ with $\dot \mu^{\xi,\theta}_{(\zeta,0), (\zeta,\eta)}$.
\end{proof}

\begin{cor}   \label{ccc}
Under the assumption of Lemma~\ref{shattered-measures}, $\AA_{\lim} $ carries a finitely additive strictly positive measure and
thus is ccc.
\end{cor}

Let us briefly discuss how the part of the iteration with one of the coordinates being 0 looks like. Say $\eta = 0$.
We then have $\AA^{\zeta,0}_{\amal} = \lim \amal_{\zeta_0 < \zeta} \AA^{\zeta, 0}_{\zeta_0, 0}$. Now, if $\zeta $ is limit,
this is simply the direct limit,
\[ \AA^{\zeta,0}_{\amal} = \lim\dir_{\zeta_0 < \zeta} \AA^{\zeta, 0}_{\zeta_0, 0} = \lim\dir_{\zeta_0 < \zeta} \left( \AA^{\zeta_0 ,0}_{\zeta_0,0} \star \dot \BB_{[\zeta_0,\zeta)} \right) ,\]
and if $\zeta = \zeta_0 + 1$ is successor, then
\[ \AA^{\zeta,0}_{\amal} =  \AA^{\zeta, 0}_{\zeta_0, 0} = \AA^{\zeta_0,0}_{\amal} \star \dot \EE_{\zeta_0,0} \star \dot \BB_{ \{ \zeta_0 \} }. \]
Thus we see this is somewhat similar to a finite support iteration adding alternately $\dot \EE_{\zeta,0}$-generics and random reals, except for the fact that in limit stages we do not
take the direct limit of  the initial stages of the iteration but the direct limit of the two-step iteration of the initial stages composed with random algebras. This is a fundamental
difference because it means that in limit stages of countable cofinality a real random over all initial extensions will be added while a standard finite support iteration
of random forcing does not add such a real~\cite[Theorem 3.2.27]{BJ95}. Of course,
the iteration in two directions is more complicated, and there are still lots of random reals over this one-dimensional iteration, namely those
added in the second coordinate $\eta$.

We are ready for the proof of Main Theorem A. We use the shattered iteration with all iterands $\dot \EE_{\zeta,\eta}$ being
Cohen forcing $\CC$. Cohen forcing is countable, so in particular $\sigma$-centered, and thus carries a finitely additive
strictly positive measure so that all the previous results apply, and the iteration is indeed ccc (Corollary~\ref{ccc}). Furthermore,
since Cohen forcing is absolute, the iteration with $\CC$ is the same as the product with $\CC$, that is, we have
\[ \AA^{\zeta,\eta}_{\zeta,\eta} = \AA^{\zeta,\eta}_{\amal} \star \dot \CC = \AA^{\zeta,\eta}_{\amal} \times \CC  \]
for any $(\zeta,\eta) \in \kappa \times \lambda$. Let us denote the name of the Cohen real added at stage $(\zeta,\eta)$ by $\dot c_{\zeta,\eta}$.
Similarly, the name of the random real at stage $\zeta < \kappa$ ($\eta < \lambda$, respectively) is $\dot r_\zeta$ ($\dot r_\eta$, resp.).

We need to compute the cardinals related to $\M$ and $\N$, but, compared to the technical work in Sections~\ref{amalgamated} and~\ref{coherent},
this turns out to be rather simple; it is the contents of the following four facts, all of which have similar proofs.

\begin{fact}   \label{fact1}
$\cov (\M) \geq \kappa$ in the $\AA$-generic extension.
\end{fact}

\begin{proof}
Let $\nu < \kappa$, and assume $\dot B_\alpha, \alpha < \nu$, are $\AA$-names for Borel meager sets. 
Since $\AA = \lim \dir_{(\zeta,\eta) \in \kappa \times \lambda} \AA^{\zeta,\eta}_{\zeta,\eta}$, by the ccc,  for all $\alpha < \nu$
there is $(\zeta_\alpha ,\eta_\alpha) \in \kappa \times \lambda$ such that $\dot B_\alpha$ is an $\AA^{\zeta_\alpha, \eta_\alpha}_{\zeta_\alpha, \eta_\alpha}$-name.
Since $\nu < \kappa, \lambda$ and the latter are both regular, we may find $\zeta > \zeta_\alpha$ and $\eta > \eta_\alpha$ for all $\alpha < \nu$.
Thus all $\dot B_\alpha$ are $\AA_{\amal}^{\zeta,\eta}$-names. Since $\AA^{\zeta,\eta}_{\zeta,\eta} = \AA^{\zeta,\eta}_{\amal} \star \dot \CC$
adds the Cohen real $\dot c_{\zeta,\eta}$ on top of the $\AA^{\zeta,\eta}_{\amal}$-extension, we see that $\AA^{\zeta,\eta}_{\zeta,\eta} $ forces that
$\dot c_{\zeta,\eta} \notin \dot B_\alpha$ for all $\alpha < \nu$. Therefore no covering family for the meager ideal can have size $< \kappa$.
\end{proof}

\begin{fact}   \label{fact2}
$\non (\N) \leq \kappa$  in the $\AA$-generic extension. In fact $\AA$ forces $\{ \dot r_\zeta : \zeta < \kappa \} \notin \N$.
\end{fact}

\begin{proof}
Let $\dot B$ be an $\AA$-name for a Borel null set. As in the previous proof, there is $(\zeta,\eta) \in \kappa \times \lambda$ such that
$\dot B$ is a $\AA^{\zeta,\eta}_{\zeta,\eta}$-name. Since $\AA^{\zeta,\eta}_{\zeta,\eta} \star \dot \BB_{[\zeta,\kappa) \dot\cup [\eta,\lambda) } \embed \AA$,
we see that $\AA$ forces that $\dot r_\zeta$ is random over the $\AA^{\zeta,\eta}_{\zeta,\eta}$-extension and, thus, that $\dot r_\zeta \notin \dot B$.
\end{proof}

\begin{fact} \label{fact3}
$\cov (\N) \geq \lambda$  in the $\AA$-generic extension.
\end{fact}

\begin{proof}
Let $\nu < \lambda$, and assume $\dot B_\alpha, \alpha < \nu$, are $\AA$-names for Borel null sets. Again we can find $(\zeta_\alpha ,\eta_\alpha) \in \kappa \times \lambda$,
$\alpha < \nu$, such that $\dot B_\alpha$ is an $\AA^{\zeta_\alpha, \eta_\alpha}_{\zeta_\alpha, \eta_\alpha}$-name. By regularity of $\lambda$,
there is $\eta < \lambda$ with $\eta_\alpha < \eta$ for all $\alpha < \nu$. As in the previous proof this means that $\AA$ forces that
$\dot r_\eta$ is random over the $\AA^{\zeta_\alpha,\eta_\alpha}_{\zeta_\alpha,\eta_\alpha}$-extension and, thus, that $\dot r_\eta \notin \dot B_\alpha$,
for any $\alpha < \nu$. Therefore no covering family for the null ideal can have size $< \lambda$.
\end{proof}

\begin{fact}   \label{fact4}
$\non (\M) \leq \lambda$  in the $\AA$-generic extension. In fact $\AA$ forces $\{ \dot c_{\zeta,\eta} : \zeta < \kappa , \eta < \lambda \} \notin \M$.
\end{fact}

\begin{proof}
Let $\dot B$ be an $\AA$-name for a Borel meager set. Again there is $(\zeta,\eta) \in \kappa \times \lambda$ such that
$\dot B$ is a $\AA^{\zeta,\eta}_{\amal}$-name, and
we see that $\AA$ forces that $\dot c_{\zeta,\eta}$ is Cohen over the $\AA^{\zeta,\eta}_{\amal}$-extension and, thus, that $\dot c_{\zeta,\eta} \notin \dot B$.
\end{proof}

Using the ZFC-inequalities $\cov (\M) \leq \non (\N)$ and $\cov (\N) \leq \non (\M)$ (\cite[Theorem 2.1.7]{BJ95}, see also Introduction), we infer that $\cov (\M) = \non (\N) = \kappa$ and 
$\cov (\N) = \non (\M) = \lambda$, and the proof of Main Theorem A is complete.



\section{Preservation results for shattered iterations}   \label{preservation}

Preservation results for iterated forcing constructions are theorems saying that if the initial stages or the quotients of an iteration satisfy
a certain property then so does the whole iteration. There are a plethora of preservation results for classical iteration techniques like
finite support iteration and countable support iteration (see, in particular, \cite[Chapter 6]{BJ95}, \cite{Go93}, and~\cite{Sh98}),
starting with the preservation of the ccc for the former, and of properness for the latter. Unfortunately, very little is known with regard to
preservation for shattered iterations, and this has to do with the amalgamated limit construction (see the discussion at the beginning of
Section~\ref{coherent}). A major result is the preservation of finitely additive measurability, and this, like other preservation theorems, has impact
beyond the preservation of the ccc, as we shall now see.

Let $\I$ be a non-trivial $\sigma$-ideal on the real numbers (by ``non-trivial" we mean that $\{ x \} \in \I$ for all $x \in \RR$ and $\RR \notin \I$).
Define 
\begin{itemize}
\item $\add (\I) := \min \{ | \F | : \F \sub \I$ and $\bigcup \F \notin \I \}$, the {\em additivity} of $\I$ and
\item $\cof (\I) := \min \{ | \F | : \F \sub \I$ and for all $A \in \I$ there is $B \in \F$ with $A \sub B \}$, the {\em cofinality} of $\I$.
\end{itemize}
Clearly $\aleph_1 \leq \add (\I) \leq \cov (\I), \non (\I) \leq \cof (\I)$ for any non-trivial $\sigma$-ideal, and additionally 
$\add (\M) = \min \{ \cov (\M), \bb \}$ and $\cof (\M) = \max \{ \non (\M) , \dd \}$~\cite[Corollary 2.2.9 and Theorem 2.2.11]{BJ95},
as well as $\add (\N) \leq \add (\M)$ and $\cof (\M) \leq \cof (\N)$~\cite[Theorem 2.3.7]{BJ95}, the famous Bartoszy\'nski-Raisonnier-Stern Theorem.
The cardinal invariants related to measure and category as well as $\bb$ and $\dd$ are usually displayed in {\em Cicho\'n's diagram} 
(see~\cite[Chapter 2]{BJ95} or~\cite[Section 5]{Bl10} for details), where cardinals grow as one moves
up or right.
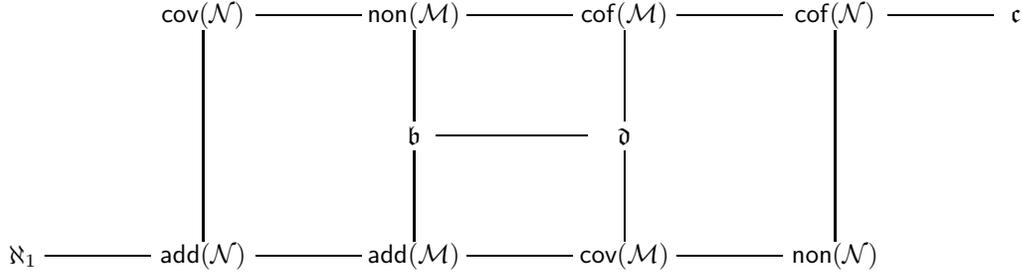
\[
\setlength{\unitlength}{0.2000mm}
\begin{picture}(750.0000,180.0000)(10,15)
\thinlines
\put(615,180){\line(1,0){70}}
\put(475,180){\line(1,0){70}}
\put(195,180){\line(1,0){70}}
\put(335,180){\line(1,0){70}}
\put(315,100){\line(1,0){100}}
\put(475,20){\line(1,0){70}}
\put(335,20){\line(1,0){70}}
\put(195,20){\line(1,0){70}}
\put(55,20){\line(1,0){70}}
\put(160,30){\line(0,1){140}}
\put(580,30){\line(0,1){140}}
\put(440,30){\line(0,1){60}}
\put(440,110){\line(0,1){60}}
\put(300,110){\line(0,1){60}}
\put(300,30){\line(0,1){60}}
\put(550,170){\makebox(60,20){$ \cof(\mathcal{N})$}}
\put(550,10){\makebox(60,20){$\non(\mathcal{N})$}}
\put(410,170){\makebox(60,20){$\cof(\mathcal{M})$}}
\put(410,90){\makebox(60,20){$\mathfrak{d}$}}
\put(410,10){\makebox(60,20){$ \cov(\mathcal{M})$}}
\put(270,10){\makebox(60,20){$\add(\mathcal{M})$}}
\put(270,90){\makebox(60,20){$\mathfrak{b}$}}
\put(270,170){\makebox(60,20){$ \non(\mathcal{M})$}}
\put(130,170){\makebox(60,20){$\cov(\mathcal{N})$}}
\put(130,10){\makebox(60,20){$ \add(\mathcal{N})$}}
\put(20,10){\makebox(40,20){$\aleph_1$}}
\put(680,170){\makebox(40,20){$\mathfrak{c}$}}
\end{picture} 
\]
Let $g \in \omom$. A function $\varphi \in ( \omloms )^\omega$ is a {\em $g$-slalom} if $| \varphi (n) | \leq g(n)$ for all $n$.
$\varphi$ {\em localizes} $f \in \omom$ if for all but finitely many $n \in \omega$, $f(n) \in \varphi (n)$. The following
is implicit in Kamburelis' work.

\begin{prop}[Kamburelis {\cite[Proposition 1.1]{Ka89}}]  \label{measure-localization}
Let $\AA$ be a cBa with a finite additive strictly positive measure $\mu$. Let $g \in \omom$, let $\dot\varphi$ be an $\AA$-name
for a $g$-slalom, and let $h(n) = n \cdot g(n)$ for all $n$. Then there is an $h$-slalom $\psi$ such that for all $f \in \omom$,
if $\psi$ does not localize $f$ then $\forces ``\dot\varphi$ does not localize $f"$.
\end{prop}

\begin{proof}
Clearly, for each $n$, $\sum_k \mu \left( \Bool k \in \dot \varphi (n) \Boor \right) \leq g(n)$ and thus, letting $\psi (n) = \left\{ k :
\mu \left( \Bool k \in \dot \varphi (n) \Boor \right) \geq {1 \over n } \right\}$, we see that $| \psi (n) | \leq n \cdot g(n) = h(n)$.
Assume $\psi$ does not localize $f \in\omom$. Then $\mu \left( \Bool f(n) \in \dot \varphi (n) \Boor \right) <{1 \over n } $
for infinitely many $n$, and therefore, given any $n_0$ and any $p \in \AA$, we may find $n > n_0$ with
$\mu \left( \Bool f(n) \in \dot \varphi (n) \Boor \right) < \mu (p)$ and thus $q: = p - \Bool f(n) \in \dot \varphi (n) \Boor > \zero$.
Clearly $q$ forces $f(n) \notin \dot\varphi (n)$. Hence $\forces ``\dot\varphi$ does not localize $f"$.
\end{proof}

We recall:

\begin{theorem}[Bartoszy\'nski, {see~\cite[Proposition 1.4]{Ka89} and~\cite[Theorem 2.3.9]{BJ95}}]   \label{Bartoszynski}
Fix $g \in \omom$ with $\lim_n g(n) = \infty$. Then
\begin{itemize}
\item $\add (\N) = \min \{ | \F| : \F \sub \omom $ and no $g$-slalom localizes all elements of $\F \} $ and
\item $\cof (\N) = \min \{ | \Phi |  : \Phi$ is a set of $g$-slaloms and for all $f \in\omom$ there is $\varphi \in \Phi$ such that
   $\varphi $ localizes $f \}$.
\end{itemize}
\end{theorem}

\begin{cor}  \label{measure-addN}
Assume $\add (\N) = \aleph_1$ and $\cof (\N) = \cc = \nu$. Let $\AA$ be a  cBa with a finitely additive strictly positive  measure
$\mu$ of size $\nu$. Then in the $\AA$-generic extension we still have $\add (\N) = \aleph_1$ and $\cof (\N) = \cc = \nu$.
\end{cor}

Note that an infinite cBa $\AA$ has size at least $\cc$. If it is larger then forcing with $\AA$ may increase the size of $\cc$. This accounts for the size restriction.
Also note that the assumption can be forced by first adding $\nu = \nu^\omega$ many Cohen reals over a model of, say, GCH~\cite[Model 7.5.8]{BJ95}.

\begin{proof}
Since $\AA$ is ccc of size $\cc$, $\cc = \nu$ is preserved. We use the characterization of $\add (\N)$ and $\cof (\N)$ of the
previous theorem. Fix any $g \in \omom$ with $\lim_n g(n) = \infty$ and let $h (n) = n \cdot g(n)$ for all $n$.

To see $\add (\N) = \aleph_1$, fix a family $\F \sub \omom$ of size $\aleph_1$ that is not localized by any $h$-slalom.
By Proposition~\ref{measure-localization}, $\F$ is not localized by any $g$-slalom of the extension, as required.

To see $\cof (\N) = \nu$, let $\nu' < \nu$ and let $\{ \dot \varphi_\alpha : \alpha < \nu ' \}$ be $\AA$-names for $g$-slaloms.
By  Proposition~\ref{measure-localization}, we find $h$-slaloms $\psi_\alpha$, $\alpha < \nu '$, satisfying the conclusion of the proposition.
By $\cof (\N) = \cc$ in the ground model, there is $f \in \omom$ that is not localized by any $\psi_\alpha$. Therefore
$f $ is forced not to be localized by any $\dot\varphi_\alpha$. 
\end{proof}

\begin{cor}   \label{shattered-Cohen}
Assume $\add (\N) = \aleph_1$, $\cof (\N) = \cc = \nu$ and let $\kappa,\lambda$ be regular cardinals with $\aleph_1 \leq \kappa < \lambda \leq \nu$.
In the extension by the shattered iteration of Cohen forcing of Section~\ref{Cohen}, we still have $\add (\N) = \aleph_1$ and $\cof (\N) = \cc = \nu$.
\end{cor}

In the original model of Section~\ref{Cohen}, we said nothing about the value of $\cc$ in the ground model. This is because its value
is irrelevant for the computation of the uniformities and covering numbers of the meager and null ideals. If GCH holds, or just $\lambda^\omega = \lambda$,
then $\cc = \lambda$ will be forced. But we may also start with a model where $\cc$ is larger than $\lambda$, and that's what we do here.

\begin{proof}
By Corollaries~\ref{ccc} and~\ref{measure-addN}. (Note here that $\AA$, being a ccc cBa generated by an iteration of length $\kappa \times \lambda$ of forcing
notions adding single reals, has size $\lambda^\omega = \nu$.)
\end{proof}

We next consider the values of $\bb$ and $\dd$ (and $\add (\M)$ and $\cof (\M)$) in the model of Section~\ref{Cohen}. We do believe
they are preserved as well, but we are unable to prove this.

\begin{conj}   \label{b-d-conjecture}
Assume $\bb = \aleph_1$, $\dd = \cc = \nu$ and let $\kappa,\lambda$ be regular cardinals with $\aleph_1 \leq \kappa < \lambda \leq \nu$.
In the extension by the shattered iteration of Cohen forcing, we still have $\bb= \add (\M) =\aleph_1$ and $\dd = \cof (\M) = \cc = \nu$.
\end{conj}

To show this two approaches are natural:
\begin{itemize}
\item Show the unbounded family of size $\aleph_1$ from the ground model is preserved. The techniques of~\cite[Section 1]{BJ93} may be relevant here.
\item Show that the ``first" $\aleph_1$ many Cohen reals, that is, $\{ c_{\zeta, 0} : \zeta < \aleph_1 \}$, remain unbounded in the whole extension.
\end{itemize}
We do know that $\bb \leq \kappa$ and $\dd \geq \lambda$ in our model, and, since we need the argument anyway for the results in the next section, we
provide it here. But this is clearly not optimal.

We start with a useful characterization of unbounded reals over intermediate extensions.

\begin{lem}   \label{unbounded-charac}
Assume $\AA_0 \embed \AA$ are cBa's. Let $\dot f \in \omom$ be an $\AA$-name. Then the following are equivalent:
\begin{enumerate}
\item The trivial condition of $\AA$ forces that $\dot f $ is unbounded over $V^{\AA_0} \cap \omom$.
\item For any $p \in \AA$ there are $q \leq_{\AA_0} h_{\AA \AA_0} (p)  $ and $n$ such that for all $k$ and all $q' \leq_{\AA_0} q$, 
   $q'$ and $p \cdot \Bool \dot f (n) \geq k \Boor$ are compatible. 
\end{enumerate}
\end{lem}

Note that the conclusion of (ii) can be rewritten as: ``for all $k$, $h_{\AA \AA_0} \left( p \cdot \Bool \dot f (n) \geq k \Boor \right) \geq_{\AA_0}  q$".

\begin{proof}
\underline{(i) $\Loriar$ (ii)}. Suppose there is $p \in \AA$ such that for all $q \leq_{\AA_0} h_{\AA \AA_0} (p)$ and all $n$, there are $k$ and  $q' \leq_{\AA_0} q$ 
such that $q'$ and $p \cdot \Bool \dot f (n) \geq k \Boor$ are incompatible. Note that the latter is equivalent to saying $q' \cdot p \forces_\AA \dot f (n) < k$.
For $n \in \omega$, let $D_n = \{ q' \leq_{\AA_0} h_{\AA\AA_0} (p) : \exists k \; ( q' \cdot p \forces_\AA \dot f (n) < k ) \}$, and note that
by assumption all $D_n$ are predense below $h_{\AA\AA_0} (p)$ in $\AA_0$. So, taking an $\AA_0$-generic filter $G_0$ with $h_{\AA\AA_0} (p) \in G_0$,
there are $q_n \in D_n \cap G_0$ and numbers $ h(n)$ witnessing $q_n \in D_n$ (i.e. $q_n \cdot p \forces_\AA \dot f (n) < h(n)$). 
Now, $h \in V [G_0] \cap \omom$ and $p \forces_{\AA / G_0} \dot f(n) < h(n)$ for all $n$, a contradiction.

\underline{(ii) $\Loriar$ (i)}. Assume there are an $\AA_0$-name $\dot h \in \omom$ and a condition $p \in \AA$ such that $p \forces_\AA \dot f(n) < \dot h(n)$
for all $n$. Let $q \leq_{\AA_0} h_{\AA\AA_0} (p)$ and $n$ be as in (ii). Choose $q' \leq_{\AA_0} q$ deciding the value of $\dot h (n)$,
say $q' \forces_{\AA_0} \dot h (n) = k$. By assumption, $q'$ is compatible with $p \cdot \Bool \dot f (n) \geq k \Boor$, and a common extension forces
$\dot f (n) \geq k = \dot h(n)$, a contradiction.
\end{proof}

\begin{lem}   \label{shattered-unbounded}
Let  $\left\la \AA^{\xi,\theta}_{\zeta,\eta} : \zeta \leq \xi < \kappa, \eta \leq \theta < \lambda \right\ra$ be a shattered iteration of
forcing notions $\left\la \dot \EE_{\zeta,\eta} : \zeta < \kappa, \eta < \lambda \right\ra$ such that each $\dot \EE_{\zeta,\eta}$ adds an unbounded
real $\dot e_{\zeta,\eta}$ over $\AA^{\zeta,\eta}_{\amal}$. Then:
\begin{enumerate}
\item  For all $\eta$ and all $\xi < \zeta$, $\forces_{\AA^{\zeta,\eta}_{\zeta,\eta}} `` \dot e_{\zeta,0}$ is unbounded over $V^{\AA^{\zeta,\eta}_{\xi,\eta}} \cap \omom"$.
\item For all $\zeta$ and all $\theta < \eta$, $\forces_{\AA^{\zeta,\eta}_{\zeta,\eta}} `` \dot e_{0,\eta}$ is unbounded over $V^{\AA^{\zeta,\eta}_{\zeta,\theta}} \cap \omom"$.
\end{enumerate}
\end{lem}

Notice that, by the discussion after Corollary~\ref{ccc} in Section~\ref{Cohen},
the $\dot e_{\zeta,0}$ are added by an iteration adding alternatively unbounded and random reals. Similarly for the $\dot e_{0,\eta}$.

\begin{proof}
The two statements are clearly symmetric, so it suffices to prove the first. This is done by induction on $\eta$, showing somewhat more generally that 
$\forces_{\AA^{\zeta,\theta}_{\zeta,\eta}} `` \dot e_{\zeta,0}$ is unbounded over $V^{\AA^{\zeta,\theta}_{\xi,\eta}} \cap \omom$" for any $\theta \geq \eta$.

\underline{$\eta = 0$}. Since $\dot e_{\zeta,0}$ is unbounded over $V^{\AA^{\zeta,0}_{\amal}} \cap \omom$ in $V^{\AA^{\zeta,0}_{\zeta,0}}$, it is
also unbounded over $V^{\AA^{\zeta,0}_{\xi, 0}} \cap \omom$. Furthermore, as $\AA^{\zeta,\theta}_{\zeta',0} = \AA^{\zeta,0}_{\zeta' , 0} \star \dot\BB_\theta$ for $\zeta' \leq \zeta$,
and random forcing is $\omom$-bounding, $\dot e_{\zeta,0}$ is still unbounded over $V^{\AA^{\zeta,\theta}_{\xi, 0}} \cap \omom$ in $V^{\AA^{\zeta,\theta}_{\zeta,0}}$.

\underline{$\eta > 0$}. We assume the statement has been proved for all $\eta' < \eta$ and prove it for $\eta$. As explained in the
previous paragraph, increasing the superscript $\eta$ to $\theta$ preserves the correctness of the statement by the $\omom$-bounding property of random forcing,
so it suffices to show it for $\theta = \eta$.   Moreover, since both $\AA^{\zeta,\eta}_{\zeta,0}$ and $\AA^{\zeta,\eta}_{\xi,\eta}$ completely embed into
$\AA^{\zeta,\eta}_{\amal} \embed \AA^{\zeta,\eta}_{\zeta,\eta}$, it suffices to see that the statement is forced by $\AA^{\zeta,\eta}_{\amal} $. 

We verify (ii) in Lemma~\ref{unbounded-charac}. To this end take $(p,q) \in D$, where $D$ is the dense set in the definition of the amalgamated limit,
that is, $\AA^{\zeta,\eta}_{\amal} = \ro (D)$. Without loss of generality (exchanging the role of $p$ and $q$, if necessary), we may assume that
$p \in \AA^{\zeta,\eta}_{\zeta,\eta'}$ for some $\eta' < \eta$ and $q \in \AA^{\zeta,\eta}_{\zeta', \eta}$ for some $\zeta' < \zeta$. We may also suppose that
$\xi \leq \zeta'$. We then have, by definition of $D$, that $h^{\zeta,\eta}_{(\zeta,\eta'), (\zeta',\eta')} (p) = h^{\zeta,\eta}_{(\zeta',\eta), (\zeta',\eta')} (q)
\in \AA^{\zeta,\eta}_{\zeta',\eta'}$. Note also that the projection of $(p,q)$ to $\AA^{\zeta,\eta}_{\xi,\eta}$ is $h^{\zeta,\eta}_{(\zeta',\eta),(\xi,\eta)} (q)$. 

By induction hypothesis for $\eta'$ applied to $\zeta'$, $\forces_{\AA^{\zeta,\eta}_{\zeta,\eta'}} `` \dot e_{\zeta,0}$ is unbounded over $V^{\AA^{\zeta,\eta}_{\zeta',\eta'}} \cap \omom$".
So, by Lemma~\ref{unbounded-charac}, there are $r \leq_{\AA^{\zeta,\eta}_{\zeta',\eta'}} h^{\zeta,\eta}_{(\zeta,\eta') , (\zeta ' ,\eta')} (p)$ and $n$ such that
for all $k$, we have $h^{\zeta,\eta}_{(\zeta,\eta'),(\zeta',\eta')} \left( p \cdot \Bool \dot e_{\zeta,0} (n) \geq k \Boor \right) \geq_{\AA^{\zeta,\eta}_{\zeta',\eta'}} r$.
Now let $s = h^{\zeta,\eta}_{(\zeta',\eta),(\xi,\eta)} (r \cdot q) \in \AA^{\zeta,\eta}_{\xi,\eta}$. Let $k$ and $s' \leq_{\AA^{\zeta,\eta}_{\xi,\eta}} s$ be arbitrary.
Since 
\[ h^{\zeta,\eta}_{(\zeta',\eta),(\zeta',\eta')} (s'\cdot r \cdot q) \leq_{\AA^{\zeta,\eta}_{\zeta',\eta'}} r \leq_{\AA^{\zeta,\eta}_{\zeta',\eta'}} h^{\zeta,\eta}_{(\zeta,\eta'),(\zeta',\eta')}
\left( p \cdot \Bool \dot e_{\zeta,0} (n) \geq k \Boor \right), \] 
$\left( h^{\zeta,\eta}_{(\zeta',\eta),(\zeta',\eta')} (s' \cdot r \cdot q ) \cdot p \cdot \Bool \dot e_{\zeta,0} (n) \geq k \Boor, 
s' \cdot r \cdot q \right)$ belongs to $D$ as witnessed by $h^{\zeta,\eta}_{(\zeta',\eta),(\zeta',\eta')} (s'\cdot r \cdot q)$, and is a common extension of $s'$,
$(p,q)$ and $\Bool \dot e_{\zeta,0} (n) \geq k \Boor$ in $\AA^{\zeta,\eta}_{\amal}$;
so the conditions are compatible, as required. This completes the induction step and the proof of the lemma.
\end{proof}

\begin{cor}  \label{shattered-bandd}
In the model of Corollary~\ref{shattered-Cohen}, $\bb = \add (\M) \leq \kappa$ and $\dd = \cof (\M) \geq \lambda$.
\end{cor}

\begin{proof}
By Lemma~\ref{shattered-unbounded} (i), $\{\dot c_{\zeta, 0} : \zeta < \kappa \}$ is forced to be unbounded in the $\AA$-generic extension,
and we obtain $\bb \leq \kappa$. By part (ii) of the same lemma, $\{ \dot c_{0,\eta} : \eta < \lambda \}$ witnesses $\dd \geq \lambda$, for if
$\nu ' < \lambda $ and $\{ \dot f_\alpha : \alpha < \nu ' \}  \sub \omom$ then for each $\alpha < \nu'$ we find $\eta_\alpha$ such that
all $\dot c_{0,\eta}$ for $\eta \geq \eta_\alpha$ are unbounded over $\dot f_\alpha$ and, letting $\eta = \sup \{ \eta_\alpha : \alpha < \nu ' \} < \lambda$,
$\dot c_{0,\eta}$ is unbounded over all $\dot f_\alpha$ so that $\{ \dot f_\alpha : \alpha < \nu ' \}$ is not dominating.
\end{proof}



\section{The shattered iteration of other forcing notions}   \label{Hechler}

We present the models obtained by the shattered iteration of two other well-known forcing notions. We start with the
proof of Main Theorem B. See Example~\ref{exam3} in Section~\ref{correct} for the definition of Hechler forcing.

\begin{theorem}   \label{shattered-Hechler}
Assume $\add (\N) = \aleph_1$, $\cof (\N) = \cc = \mu$ and let $\kappa, \lambda$ be regular with $\aleph_1 \leq \kappa < \lambda \leq \mu$.
Let $\left\la \AA_{\zeta,\eta}^{\xi,\theta} : \zeta \leq \xi < \kappa , \eta \leq \theta < \lambda \right\ra$ be the shattered iteration of Hechler forcing, that is,
for all $(\zeta,\eta) \in \kappa \times \lambda$, $\AA^{\zeta,\eta}_{\amal}$ forces that $\dot \EE_{\zeta,\eta} = \dot \DD$. Then, in the
$\AA$-generic extension, $\add (\N) = \aleph_1$, $\add (\M) = \bb = \cov (\M) = \non (\N) = \kappa$, $\cov (\N) = \non (\M) = \dd = \cof (\M) = \lambda$,
and $\cof (\N) = \cc = \mu$.
\end{theorem}

\begin{proof}
First note that, being $\sigma$-centered, Hechler forcing carries a finitely additive strictly positive measure by the discussion before Lemma~\ref{secondextension}
in Section~\ref{coherent}. In particular $\AA$ is ccc by Corollary~\ref{ccc}.

The proofs of $\non (\N) \leq \kappa$ and $\cov (\N) \geq \lambda$ carry over from Facts~\ref{fact2} and~\ref{fact3} in Section~\ref{Cohen} because they use
the random generics $r_\zeta, \zeta < \kappa,$ and $r_\eta, \eta < \lambda$. Next, recall that if $d \in \omom$ is $\DD$-generic over $V$, then
$c \in \twoom$ defined by $c (n) = d(n) \mod 2$ is $\CC$-generic over $V$. Thus, from the Hechler generics $\dot d_{\zeta, \eta}, (\zeta,\eta) \in \kappa \times \lambda$,
added by $\dot \EE_{\zeta,\eta} = \dot \DD$, we can decode Cohen generics $\dot c_{\zeta,\eta}$, which can be used to show that
$\cov (\M) \geq \kappa$ and $\non (\M) \leq \lambda$ as in Facts~\ref{fact1} and~\ref{fact4}.

A similar argument using the $d_{\zeta,\eta}$ shows $\bb \geq \kappa$ and $\dd \leq \lambda$: 
let $\nu < \kappa$ and let $\dot f_\alpha, \alpha < \nu,$ be $\AA$-names for reals in $\omom$. There is $(\zeta,\eta) \in \kappa \times \lambda$ such that
all $\dot f_\alpha$ are $\AA^{\zeta,\eta}_{\amal}$-names. Since $\AA^{\zeta,\eta}_{\zeta,\eta} = \AA^{\zeta,\eta}_{\amal} \star \dot \DD$ adds the Hechler real
$\dot d_{\zeta,\eta}$ on top of the $\AA^{\zeta,\eta}_{\amal}$-extension, $\AA^{\zeta,\eta}_{\zeta,\eta}$ forces $\dot f_\alpha \leq^* \dot d_{\zeta,\eta}$ for all
$\alpha < \nu$. Thus the $\{ \dot f_\alpha : \alpha < \nu \}$ are bounded, and $\bb \geq \kappa$ follows. The same argument shows that
$\{ \dot d_{\zeta,\eta} : \zeta < \kappa, \eta < \lambda \}$ is forced to be a dominating family, and we obtain $\dd \leq \lambda$.

Next, $\bb \leq \kappa$ and $\dd \geq \lambda$ are proved exactly as in Corollary~\ref{shattered-bandd} using Lemma~\ref{shattered-unbounded} in Section~\ref{preservation}.

Finally, $\add (\N) = \aleph_1$ and $\cof (\N) = \mu = \cc$ follows from Corollary~\ref{measure-addN} in Section~\ref{preservation} because $\AA$ carries a finitely additive measure
by Corollary~\ref{ccc}.
\end{proof}

Let us consider some more cardinal invariants of the continuum. For sets $A, B \in \omoms$, say that $A$ {\em splits} $B$ if 
$B \cap A$ and $B \sem A$ are both infinite. A family $\A \sub \omoms$ is {\em splitting} if every $B \in \omom$ is split by some member of $\A$.
Dually, $\B \sub \omoms$ is {\em unreaped} if no $A \in \omoms$ splits all members of $\B$. Using these notions we define:
\begin{itemize}
\item $\sss : = \min \{ | \A | : \A \sub \omoms$ is a splitting family$\}$, the {\em splitting number} and
\item $\rr : = \min \{ | \B | : \B \sub \omoms$ is an unreaped family$\}$, the {\em reaping number}.
\end{itemize}
It is well-known that $\sss \leq \dd, \non (\M), \non (\N)$ and, dually, $\rr \geq \bb, \cov (\M), \cov (\N)$~\cite[Theorems 3.3, 3.8, and 5.19]{Bl10}.

Let $\U$ be a non-principal ultrafilter on the natural numbers $\omega$. {\em Mathias forcing} $\MM_\U$ with $\U$ consists of pairs $(s, A)$ with 
$s \in \omlom$ and $A \in \U$ such that $\max (s) < \min (A)$. The ordering is given by $(t,B) \leq (s,A)$ if $t \supseteq s$, $B \sub A$, and $t \sem s \sub A$.
Mathias forcing with $\U$ is $\sigma$-centered and generically adds a set $M \in \omoms$ diagonalizing the ultrafilter $\U$, that is,
$M \sub^* A$ for all $A \in \U$. Indeed, $M = \bigcup \{ s : $ there is $A \in \U$ such that $(s,A) \in G \}$ where $G$ is the $\MM_\U$-generic filter over $V$.
Clearly, $M$ is not split by any ground model set, and $\MM_\U$ thus is a natural forcing for destroying splitting families and, when iterated, increasing the
splitting number $\sss$. Furthermore, if $\U$ is not a P-point then $\MM_\U$ adds a dominating real~\cite[Lemma 4]{Ca88}.

\begin{theorem}   \label{shattered-Mathias}
Assume $\add (\N) = \aleph_1$, $\cof (\N) = \cc = \mu$ and let $\kappa, \lambda$ be regular with $\aleph_1 \leq \kappa < \lambda \leq \mu$.
Let $\left\la \AA_{\zeta,\eta}^{\xi,\theta} : \zeta \leq \xi < \kappa , \eta \leq \theta < \lambda \right\ra$ be the shattered iteration of Mathias forcing with 
non-P-point ultrafilters $\dot \U_{\zeta,\eta}$, that is, for all $(\zeta,\eta) \in \kappa \times \lambda$, $\dot \U_{\zeta,\eta}$ is an $\AA^{\zeta,\eta}_{\amal}$-name and
$\AA^{\zeta,\eta}_{\amal}$ forces that $\dot \EE_{\zeta,\eta} = \MM_{\dot \U_{\zeta,\eta}}$. Then, in the
$\AA$-generic extension, $\add (\N) = \aleph_1$, $\sss = \add (\M) = \bb = \cov (\M) = \non (\N) = \kappa$, $\rr = \cov (\N) = \non (\M) = \dd = \cof (\M) = \lambda$,
and $\cof (\N) = \cc = \mu$.
\end{theorem}

\begin{proof}
Except for the values of $\sss$ and $\rr$, the proof is exactly like the proof of Theorem~\ref{shattered-Hechler}. In particular, by the $\sigma$-centeredness
of the $\MM_{\dot \U_{\zeta,\eta}}$, $\AA$ carries a finitely additive strictly positive measure and is ccc. Furthermore, since $\sss \leq \non (\N)$ and
$\rr \geq \cov (\N)$ in ZFC, $\sss  \leq \kappa$ and $\rr \geq \lambda$ follow.

So we are left with showing $\sss \geq \kappa$ and $\rr \leq \lambda$, and this is exactly what is achieved by the Mathias generics.
Indeed, assume $\nu < \kappa$ and $\dot A_\alpha , \alpha < \nu ,$ are $\AA$-names for infinite subsets of $\omega$. Then there is $(\zeta, \eta) \in \kappa \times \lambda$
such that all $\dot A_\alpha$ are $\AA^{\zeta,\eta}_{\amal}$-names. Since $\AA^{\zeta,\eta}_{\zeta,\eta} = \AA^{\zeta,\eta}_{\amal} \star \MM_{\dot \U_{\zeta,\eta}}$ 
diagonalizes the ultrafilter $\dot \U_{\zeta,\eta}$, the generic Mathias real $\dot M_{\zeta,\eta}$ is not split by any $\dot A_\alpha$ so that they cannot form
a splitting family in the final extension. This shows $\sss \geq \kappa$. For $\rr \leq \lambda$, simply note that the same argument shows that
$\{ \dot M_{\zeta,\eta} : \zeta < \kappa, \eta < \lambda \}$ is forced to be an unreaped family, and the proof of the theorem is complete.
\end{proof}

We do not know the values of $\sss$ and $\rr$ in the models of Corollary~\ref{shattered-Cohen} and Theorem~\ref{shattered-Hechler}, except for
the trivial $\sss \leq \kappa$ and $\rr \geq \lambda$.

\begin{conj}
In the models obtained by the shattered iteration of either Cohen or Hechler forcing, $\sss = \aleph_1$ and $\rr = \cc = \mu$.
\end{conj}

The reason for this conjecture is that these iterations only involve the Suslin ccc forcing notions $\BB$ and $\CC$ or $\DD$, and it is known
that a finite support iteration of Suslin ccc forcing over a model of CH preserves $\sss = \aleph_1$~\cite[Theorem 3.6.21]{BJ95}.



\section{Prospects and problems}   \label{problems}

Building on the work in Sections~\ref{Cohen} to~\ref{Hechler} we can make more precise in what sense our 
constructions {\em add random generics in limit stages} as opposed to fsi's
which add Cohen generics in limit stages.

Let $\kappa$ be an arbitrary infinite cardinal and let $\AA$ be a cBa.
Say $\AA$ is a {\em $\kappa$-Cohen iteration} if $\CC_\kappa \embed \AA$ and for all cBa's $\AA' \embed \AA$
generated by a name for a subset of $\mu < \kappa$, there is $\alpha < \kappa$ such that
$\AA ' \star \dot \CC_{[\alpha,\kappa)} \embed \AA$ (that is, a cofinal segment of the $\CC_\kappa$-generic
is still generic over the $\AA'$-extension). Of course $\AA' \star \dot\CC_{[\alpha,\kappa)} \cong
\AA' \times \CC_{[\alpha,\kappa)}$. Clearly we have:

\begin{obs}   \label{obs1}
Assume $\AA$ is a finite support iteration of ccc forcing of length $\kappa$ where $\kappa$
is regular. Then $\AA$ is a $\kappa$-Cohen iteration. 
\end{obs}

Similarly, $\kappa$-random iterations are obtained by replacing $\CC$ by $\BB$ in the previous definition:
$\AA$ is a {\em $\kappa$-random iteration} if $\BB_\kappa \embed \AA$ and for all cBa's $\AA' \embed \AA$
generated by a name for a subset of $\mu < \kappa$, there is $\alpha < \kappa$ such that
$\AA ' \star \dot \BB_{[\alpha,\kappa)} \embed \AA$. Obviously $\CC_\kappa$ ($\BB_\kappa$, respectively) is
a $\mu$-Cohen ($\mu$-random, resp.) iteration for every regular $\mu$ with $\mu \leq \kappa$.
Furthermore:

\begin{obs}   \label{obs2}
Assume $\kappa$ is uncountable.
\begin{enumerate}
\item If $\AA$ is a $\kappa$-Cohen iteration then $\Non M \leq \kappa$ and $\Cov M  \geq \mm_{\CC_\kappa}
   \geq \kappa$ in the $\AA$-generic extension.
\item If $\AA$ is a $\kappa$-random iteration then $\Non N \leq \kappa$ and $\Cov N \geq \mm_{\BB_\kappa}
   \geq \kappa$ in the $\AA$-generic extension.  
\end{enumerate}
\end{obs}

Here, for a partial order $\PP$, $\mm_\PP$ denotes the smallest cardinal $\mu$ such that
Martin's Axiom MA$_\mu$ fails for $\PP$. It is well-known (and easy to see) that $\Cov M = \mm_\CC$ and 
$\Cov N = \mm_\BB$~\cite[Theorem 3.1.8]{BJ95} and, more generally, $\mm_{\CC_\kappa} = \cov (\M_\kappa)$ and $\mm_{\BB_\kappa} =
\cov (\N_\kappa)$ where $\M_\kappa$ and $\N_\kappa$ denote the meager and null ideals on the
space $2^\kappa$, respectively. Obviously $\mm_{\CC_\lambda} \leq \mm_{\CC_\kappa}$ and $\mm_{\BB_\lambda} \leq \mm_{\BB_\kappa}$ for $\lambda \geq \kappa$,
and the inequality consistently is strict (see~\cite[Section 5 and Appendix B]{Br06} for models for this).

\begin{proof}
By analogy of Cohen and random forcings, it suffices to prove (i). First notice that the Cohen reals $\{ c_\alpha : \alpha < \kappa \}$
added by $\CC_\kappa$ are a non-meager set and thus a witness for $\non (\M) \leq \kappa$. For indeed, if $\dot A$ is an 
$\AA$-name for a Borel meager set, there is a countably generated cBa $\AA ' \embed \AA$ such that $\dot A$ is an $\AA '$-name,
and since $\AA' \star \dot \CC_{[\alpha,\kappa)} \embed\AA$ for some $\alpha < \kappa$, the $c_\beta$ with $\alpha \leq \beta < \kappa$
are still Cohen over the $\AA'$-generic extension and thus avoid the interpretation of $\dot A$.

On the other hand, let $\mu < \kappa$ and let $A_\alpha, \alpha < \mu$, be closed nowhere dense Baire subsets of the space $2^\kappa$.
Being Baire, the $A_\alpha$ have countable support, that is, there are $C_\alpha \sub \kappa$ countable and
$A'_\alpha \sub 2^{C_\alpha}$ closed nowhere dense such that $A_\alpha = A_\alpha ' \times 2^{\kappa \sem C_\alpha}$
(see~\cite{Ku84}). Therefore there is a $\mu$-generated cBa $\AA' \embed \AA$ such that the $A_\alpha '$, and thus also the $A_\alpha$,
belong to the $\AA'$-generic extension. As $\AA' \star \dot \CC_{[\alpha,\kappa)} \embed\AA$ for some $\alpha < \kappa$,
the canonical bijection between $[\alpha,\kappa)$ and $\kappa$ yields a generic Cohen function $f \in 2^\kappa$ in the
$\AA$-generic extension avoiding all $A_\alpha$. Hence $\mm_{\CC_\kappa} = \cov (\M_\kappa) \geq \kappa$ as required.
\end{proof}

\begin{lem} \label{random-iteration}
Let $\kappa \leq \lambda$ be regular uncountable. The shattered iterations $\AA$ considered in Sections~\ref{Cohen}
to~\ref{Hechler} are both  $\kappa$-random iterations and  $\lambda$-random iterations.
\end{lem}

\begin{proof}
By construction, $\BB_\kappa \embed \AA$ and if $\AA ' \embed \AA$ is generated by a name for a subset of
$\mu < \kappa$ then, by ccc-ness, there are $\zeta < \kappa$ and $\eta < \lambda$ such that
$\AA ' \embed \AA^{\zeta, \eta}_{\zeta, \eta}$. Since $\AA ' \star \dot \BB_{[ \zeta, \kappa)} \embed
\AA^{\zeta,\eta}_{\zeta,\eta} \star \dot \BB_{[\zeta,\kappa)} \embed \AA$ by construction,
$\AA$ is indeed a $\kappa$-random iteration.

$\BB_\lambda \embed \AA$ is also obvious. Assume $\AA ' \embed \AA$ is generated by a name for
a subset of $\mu < \lambda$. By ccc-ness, $\AA '$ is (contained in the cBa) generated by some
$\AA_{\zeta_\alpha, \eta_\alpha}^{\zeta_\alpha, \eta_\alpha}$, $\alpha < \mu$. Let 
$\eta < \lambda$ be such that $\eta_\alpha < \eta$ for all $\alpha < \mu$. We claim
that $\AA ' \star \dot \BB_{[\eta, \lambda)} \embed \AA$.
To see this it suffices to show that the generic random function $f_{[\eta, \lambda)}$
avoids all Baire null sets coded in the $\AA'$-extension.
We may assume, however, that such a Baire null set is coded in an 
$\AA_{\zeta_\alpha, \eta_\alpha}^{\zeta_\alpha, \eta_\alpha}$-extension for some $\alpha < \mu$.
Since $\AA_{\zeta_\alpha, \eta_\alpha}^{\zeta_\alpha, \eta_\alpha} \star \dot \BB_{[\eta,\lambda)}
\embed \AA$, the claim follows, and we proved $\AA$ is a $\lambda$-random iteration as well. 
\end{proof}

Note that in view of Observation~\ref{obs2} (ii), this subsumes Facts~\ref{fact2} and~\ref{fact3}
in Section~\ref{Cohen}.

Fsi's usually are only $\kappa$-Cohen iterations for one regular cardinal $\kappa$, and one may
wonder what $\kappa$-random iterations for one $\kappa$ might be. In fact, we have already encountered such
iterations as one-dimensional subiterations of our two-dimensional shattered iterations, see the discussion after
Corollary~\ref{ccc} in Section~\ref{Cohen}. Unfortunately we do not know whether these one-dimensional shattered iterations are useful for any
consistency results. 

On the other hand, two-dimensional fsi's, usually referred to as {\em matrix iterations}, have been considered a lot in the literature
and can be traced back at least to work of Blass and Shelah~\cite{BS89} in the eighties. It should be mentioned, however, that most of them
are not $\mu$-Cohen iterations for two distinct $\mu$, but some are. We leave it to the reader to show the consistency of 
$\cov (\N) = \non (\M) = \kappa$ and $\cov (\M) = \non (\N) = \lambda$ for arbitrary regular $\kappa < \lambda$ using a fsi
along a $\kappa \times \lambda$-matrix, which is both a $\kappa$- and a $\lambda$-Cohen iteration (see the Introduction for an alternative proof sketch of this). The question of constructing
higher-dimensional (at least three-dimensional) matrix iterations has been considered, but one runs into problems very similar to Problem~\ref{amal3-problem} in Section~\ref{amalgamated}:
it is not clear how to show complete embeddability in the limit step. This problem has been solved in a special case in~\cite{FFMM18}.
For shattered iterations we ask:

\begin{question}  \label{threedimensions}
Can we build a non-trivial shattered iteration $\AA$ that is a $\mu$-random iteration for three distinct $\mu$? 
\end{question}

By ``non-trivial" here we mean that $\AA$ is not a measure algebra.
A related problem is:

\begin{problem}  \label{conj2}
Show the consistency of  $\bb = \lambda$ and $\Cov M = \kappa$ for arbitrary  regular uncountable $\kappa < \lambda$!
\end{problem}

In all known models for $\bb > \cov (\M)$, $\cov (\M) = \aleph_1$ holds. To solve this problem one would have to somehow embed 
a $\lambda$-stage fsi of Hechler forcing into the shattered iteration framework, see in particular the proof of Main Theorem B in Section~\ref{Hechler}.
This would involve triple amalgamations which do not exist in general (see Counterexample~\ref{counterexam3} in Section~\ref{amalgamated}) but which might exist in the
case relevant for the construction.

\begin{question}   \label{shattered-maximum}
Can we separate more cardinals in shattered iteration models?
\end{question}

So far, we get at most four values, $\aleph_1 < \kappa < \lambda < \cc = \mu$. In the fsi context, the maximum number of distinct values,
referred to as {\em Cicho\'n's maximum}, has been achieved in recent breakthrough work of Goldstern, Kellner, Mej\'ia, Shelah, and others (see, in particular,
\cite{GKS19, GKMSta1, KST19, KTT18}). A particular instance of Question~\ref{shattered-maximum} is:

\begin{question}
Is $\aleph_1 < \cov (\M) < \non (\N) < \cov (\N) < \non (\M)$ consistent?
\end{question}

We do not know what happens if either $\kappa$ or $\lambda$ is singular.

\begin{question}   \label{conj3}
Assume either $cf (\kappa) < \kappa < \lambda$ or $\kappa \leq cf (\lambda) < \lambda$. Is $\Cov M = \Non N = \kappa$
and $\Cov N = \Non M = \lambda$ consistent?
\end{question}

Note that $\cof (\I) = \cov (\I)$ implies $\non (\I) \leq cf (\cof (\I))$ for any ideal $\I$ (this is an old observation of Fremlin). Therefore, if $\cof (\N) = \lambda$
in the second part of this question, we must have $\kappa \leq cf (\lambda)$. Making $\kappa$ or $\lambda$ singular also seems to involve amalgamations
over structures more complicated than distributive almost-lattices. 
A positive solution to the first part  may be a first step towards a positive solution of:

\begin{question}[Bartoszy\'nski-Judah, {see~\cite{BJ89} or~\cite[2.10]{Sh666}}]     \label{conj4}
Is $\Add M > cf (\Cov M)$  consistent?
\end{question}



\end{document}